\documentclass[12pt]{amsart}
\usepackage{amsmath}
\usepackage{amsbsy}
\usepackage{amssymb}
\usepackage{latexsym}
\newtheorem{theorem}{Theorem}[section]
\newtheorem{corollary}[theorem]{Corollary} 
\newtheorem{lemma}[theorem]{Lemma}
\newtheorem{proposition}[theorem]{Proposition}
\theoremstyle{definition}
\newtheorem{definition}[theorem]{Definition}
\theoremstyle{remark}
\newtheorem{remark}[theorem]{Remark}

\newtheorem{problem}[theorem]{Problem}

\numberwithin{equation}{section}

\def\bbC{{\mathbb C}}
\def\bbZ{{\mathbb Z}}
\def\bbN{{\mathbb N}}
\def\bbR{{\mathbb R}}
\def\bbQ{{\mathbb Q}}

\def\bfa{{\boldsymbol A}}
\def\bfb{{\boldsymbol B}}
\def\bfc{{\boldsymbol C}}
\def\bfd{{\boldsymbol D}}
\def\bfe{{\boldsymbol E}}
\def\bff{{\boldsymbol F}}
\def\bfh{{\boldsymbol H}}
\def\bfo{{\boldsymbol O}}
\def\bfr{{\boldsymbol R}}
\def\bfu{{\boldsymbol U}}
\def\bfx{{\boldsymbol X}}
\def\bfv{{\boldsymbol V}}

\def\bfzero{{\boldsymbol 0}}

\def\a{{\mathcal A}}
\def\b{{\mathcal B}}
\def\calc{{\mathcal C}}
\def\e{{\mathcal E}}
\def\f{{\mathcal F}}
\def\g{{\mathcal G}}
\def\h{{\mathcal H}}
\def\k{{\mathcal K}}
\def\l{{\mathcal L}}
\def\m{{\mathcal M}}
\def\n{{\mathcal N}}
\def\o{{\mathcal O}}
\def\p{{\mathcal P}}
\def\q{{\mathcal Q}}
\def\r{{\mathcal R}}
\def\calu{{\mathcal U}}
\def\calv{{\mathcal V}}
\def\forb{{\mathcal Forb}}
\def\norm#1{{\lVert{#1}\rVert}}
\begin{document}
\title[Fra\"{\i}ss\'{e} Limits, Ramsey Theory, and Topological
Dynamics]
{\bf Fra\"{\i}ss\'{e} Limits, Ramsey Theory, and Topological Dynamics of
Automorphism Groups}
\author{A. S. Kechris, V. G. Pestov, and S. Todorcevic}
\address{Department of Mathematics, Caltech 253-37,
Pasadena, CA 91125}
\email{kechris@caltech.edu}
\address{Department of Mathematics and Statistics,
University of Ottawa,
 585 King Edward Avenue,
 Ottawa, Ontario,
 Canada K1N6N5
}
\email{vpest283@uottawa.ca}
\address{ U.F.R de Math\'{e}matiques, U.M.R. 7056
 Universit\'{e} Paris 7,
 2, Pl. Jussieu, Case 7012,
 75251 Paris Cedex 05,
 France}
\email{stevo@math.jussieu.fr}

\date{}
\maketitle

\section{Introduction}

{\bf (A)} We study in this paper some connections between the Fra\"{\i}ss\'{e}
theory of amalgamation classes and ultrahomogeneous structures, Ramsey theory, and
topological dynamics of automorphism groups of countable structures.

A prime concern of topological dynamics is the study of continuous
actions of (Hausdorff) topological groups $G$ on (Hausdorff) compact
spaces $X$.  These are usually referred to as (compact) $G$-{\it flows}.
Of particular interest is the study of {\it minimal} $G$-{\it flows},
those for which every orbit is dense.  Every $G$-flow contains a minimal
subflow.  A general result of topological dynamics asserts that every
topological group $G$ has a {\it universal minimal flow} $M(G)$, a 
minimal $G$-flow which can be homomorphically mapped onto any other
minimal $G$-flow.  Moreover, this is uniquely determined, by this property,
up to isomorphism.  (As usual a {\it homomorphism} $\pi :X\rightarrow Y$
between $G$-flows is a continuous $G$-map and an {\it isomorphism} is a
bijective homomorphism.)  For separable, metrizable groups $G$, which are
the ones that we are interested in here, the universal minimal flow of
$G$ is an inverse limit of manageable, i.e., metrizable $G$-flows, but
itself may be very complicated, for example non-metrizable.  In fact, for
the ``simplest'' infinite $G$, i.e., the countable discrete ones, 
$M(G)$ is a very complicated compact $G$-invariant subset of the space $\beta G$ of
ultrafilters on $G$ and is always non-metrizable.

Rather remarkably, it turned out that there are non-trivial
topological groups $G$
for which $M(G)$ is actually trivial, i.e., a singleton.  This is
equivalent to saying that $G$ has a very strong fixed point property,
namely every $G$-flow has a fixed point (i.e., a point $x$ such that
$g\cdot x=x,\ \forall g\in G$).  (For separable, metrizable groups this is 
also equivalent to the fixed point property restricted to metrizable $G$-flows.)  
Such groups are said to have the
{\it fixed point on compacta property} or be {\it extremely amenable.}
The latter name comes from one of the standard characterizations of second
countable locally compact amenable groups.  A
second countable locally compact group $G$ is {\it amenable} iff every
metrizable $G$-flow has an invariant (Borel probability) measure.  However,
no non-trivial
locally compact group can be extremely amenable, because, by a theorem
of Veech \cite{V}, every such group admits a {\it free} $G$-flow (i.e., a flow
for which $g\cdot x=x\Rightarrow g=1_G$).  Nontriviality of the universal minimal flow for locally compact groups also follows from the
earlier results of Granirer-Lau \cite{GraL}.  This probably explains the rather
late emergence of extreme amenability. Note that the corresponding property for semigroups is much more 
common and easier to come by, and in fact the study of the fixed point 
on compacta property was initiated in the context of topological 
semigroups by Mitchell \cite{Mit}, followed by Granirer \cite{Gra}. In 1966, Mitchell \cite{Mit} asked the question of existence of 
non-trivial extremely amenable topological groups. The first 
examples of such groups were constructed by Herer-Christensen \cite{HC}.
They found Polish abelian so-called {\it pathological groups},
i.e., topological groups with no non-trivial unitary representations.  
Then they showed (see Theorem 4 in their paper) that every amenable pathological
group is extremely amenable.  Remarkably though it turned out that a lot
of important (non-locally compact) Polish groups are indeed extremely
amenable.  Gromov-Milman \cite{GroM} showed that the unitary group of infinite
dimensional separable Hilbert space is extremely amenable, Furstenberg-Weiss (unpublished) and independently Glasner \cite{Gl98} showed that the group of measurable maps from $I=[0,1]$
to the unit circle $\mathbb T$ is extremely amenable, Pestov \cite{P98a}
(see also \cite{P98}) showed that
the groups $H_+ (I),\ H_+(\bbR )$ of increasing homeomorphisms
of $I,\bbR$, resp., are extremely amenable, and Pestov \cite{P98a} showed that
the group Aut$(\langle\bbQ ,<\rangle )$ of automorphisms of the rationals
is extremely amenable.  More recently, Pestov \cite{P02} proved that the
universal Polish group Iso$(\bfu )$, of all isometries of the Urysohn space
$\bfu$, is extremely amenable, and Giordano-Pestov \cite{GiP,GiP2} showed that the group
Aut$(I,\lambda )$ (resp., Aut$^*(I,\lambda ))$ of measure preserving
automorphisms of $I$ with Lebesgue measure $\lambda$ (resp., measure-class
preserving automorphisms of $I,\lambda$) is extremely amenable.

In most known examples of extremely amenable groups, beginning with the
result by Gromov and Milman \cite{GroM} on the unitary group, extreme amenability
was established by using the phenomenon of concentration of measure on
high-dimensional structures, see, e.g., Milman and Schechtman \cite{MS} or 
Ledoux \cite{L}.  However,
we will not touch upon this subject here, referring the reader instead to
the introductory article \cite{P02b} or the most recent work in this direction
\cite{GiP2} and references therein.

Beyond the extremely amenable groups there were very few cases of metrizable
universal minimal flows that had been computed.  The first such example
is in Pestov \cite{P98a}, where the author shows that the universal minimal flow
of $H_+(\mathbb T )$, the group of orientation preserving homeomorphisms of
the circle, has as a universal minimal flow its natural (evaluation) action
on $\mathbb T$.  Then Glasner-Weiss \cite{GlW02} showed that the universal minimal flow
of $S_\infty$, the infinite symmetric group of all permutations of $\bbN$,
is its canonical action on the space of all linear orderings on $\bbN$.
Finally, Glasner-Weiss \cite{GlW03} showed that the universal minimal flow 
of $H(2^\bbN )$, the group of 
homeomorphisms of the Cantor space, is its canonical action on the space
of maximal chains of compact subsets of $2^\bbN$, a space introduced in
\cite{Usp00}.

\medskip
{\bf (B)} Motivated by Pestov's result that Aut$(\langle\bbQ ,<\rangle )$
is extremely amenable and the Glasner-Weiss computation of the universal
minimal flow of $S_\infty$, we develop in this paper a general framework,
in which such results can be viewed as special instances.  In particular, this
gives many new examples of automorphism groups that are extremely amenable
and calculations of universal minimal flows.  There are two main
ingredients that come into play here.  The first is the Fra\"{\i}ss\'{e} theory
of amalgamation classes and ultrahomogeneous structures, and the second
is the structural Ramsey theory that arises in the works of Graham, Leeb,
Rothschild, Ne\v set\v ril and R\"odl.  As repeatedly stressed 
by one of the present authors (see, e.g., Pestov \cite{P02b}), extreme amenability is
related to Ramsey-type phenomena.  For instance, Pestov's proof that 
Aut$(\langle\bbQ ,<\rangle )$ is extremely amenable depends on the
classical finite Ramsey theorem and in fact it is equivalent to it.
Generalizing this, we will see that, once things are put in the proper
context, extreme amenability of automorphism groups and calculation
of universal minimal flows turn out to have equivalent formulations in 
terms of concepts that have arisen in structural Ramsey theory.

\medskip
{\bf (C)} Let us first review some basic facts of the Fra\"{\i}ss\'{e} theory.
A (countable) {\it signature} consists of a set of symbols $L=\{R_i\}_{i\in I}
\cup\{f_j\}_{j\in J}$ ($I,J$ countable), to each of which there is an
associated {\it arity} $n(i)\in\{1,2,\dots\}\ (i\in I)$ and $m(j)\in\bbN\ 
(j\in J)$.  We call $R_i$ the {\it relation symbols} and $f_j$ the
{\it function symbols} of $L$.  A {\it structure} for $L$ is of the form
$\bfa =\langle A,\{R^\bfa_i\}_{i\in I},\ \{f^\bfa_j\}_{j\in J}\rangle$,
where $A \not= \emptyset$, $R^\bfa_i\subseteq A^{n(i)},\ f^\bfa_j:A^{m(j)}\rightarrow A$.
The set $A$ is called the {\it universe} of the structure.  An {\it embedding}
between structures $\bfa ,\bfb$ for $L$ is an injection $\pi :A\rightarrow B$
such that $R^\bfa_i(a_1,\dots , a_{n(i)})\Leftrightarrow R^\bfb_i(\pi (a_1),
\dots ,\pi (a_{n(i)}))$ and $\pi (f^\bfa_j(a_1,\dots ,a_{m(j)}))=f^\bfb_j
(\pi (a_1),\dots ,\pi (a_{m(j)}))$.  If $\pi$ is the identity, we say that
$\bf A$ is a {\it substructure} of $\bf B$.  An {\it isomorphism} is an onto 
embedding.  We write $\bfa\leq\bfb$ if $\bfa$ can be embedded in $\bfb$
and $\bfa\cong\bfb$ if $\bfa$ is isomorphic to $\bfb$.

A class $\k$ of finite structures for $L$ is {\it hereditary} if $\bfa
\leq\bfb\in\k$ implies $\bfa\in\k$.  It satisfies the {\it joint
embedding property} if for any $\bfa ,\bfb\in\k$ there is $\bfc\in\k$
with $\bfa\leq\bfc ,\ \bfb\leq\bfc$.  Finally, it satisfies the {\it
amalgamation property} if for any embeddings $f:\bfa\rightarrow\bfb ,\ 
g:\bfa\rightarrow\bfc$, with $\bfa ,\bfb ,\bfc\in\k$, there is $\bfd\in\k$
and embeddings $r:\bfb\rightarrow\bfd$ and $s:\bfc\rightarrow\bfd$, such
that $r\circ f=s\circ g$.  We call $\k$ a {\it Fra\"{\i}ss\'{e} class} if it is
hereditary, satisfies joint embedding and amalgamation, 
contains only countably many structures, up to 
isomorphism, and contains
structures of arbitrarily large (finite) cardinality.

If now $\bfa$ is a countable structure, which is locally finite (i.e., 
finitely generated substructures are finite), its {\it age}, Age$(\bfa )$,
is the class of all finite structures which can be embedded in $\bfa$.
We call $\bfa$ {\it ultrahomogeneous} if every isomorphism between
finite substructures of $\bfa$ can be extended to an automorphism of $\bfa$.
We call a locally finite, countably infinite, ultrahomogeneous structure a {\it Fra\"{\i}ss\'{e}
structure}.

There is a canonical 1-1 correspondence between Fra\"{\i}ss\'{e} classes and
structures, discovered by Fra\"{\i}ss\'{e}.  If $\bfa$ is a Fra\"{\i}ss\'{e} structure,
then Age$(\bfa )$ is a Fra\"{\i}ss\'{e} class.  Conversely, if $\k$ is a
Fra\"{\i}ss\'{e} class, then there is a unique Fra\"{\i}ss\'{e} structure, the {\it
Fra\"{\i}ss\'{e} limit of} $\k$, denoted by Flim$(\k )$, whose age is exactly
$\k$.  Here are a couple of examples:  the Fra\"{\i}ss\'{e} limit of the class of finite
linear orderings is $\langle\bbQ ,<\rangle$, and the Fra\"{\i}ss\'{e} limit of the class of
finite graphs is the {\it random graph}.

\medskip
{\bf (D)} We now come to structural Ramsey theory.  Let $\k$ be a hereditary class of
finite structures in a signature $L$.  For $\bfa\in\k ,\ \bfb\in\k$ with
$\bfa\leq\bfb$, we denote by $\begin{pmatrix}\bfb\\ \bfa\end{pmatrix}$
the set of all substructures of $\bfb$ isomorphic to $\bfa$.  If $\bfa
\leq\bfb\leq\bfc$ are in $\k$ and $k=2,3,\dots$, we write
\[
\bfc\rightarrow (\bfb )^\bfa_k,
\]
if for every coloring $c:\begin{pmatrix}\bfc \\ \bfa\end{pmatrix}
\rightarrow \{1,\dots ,k\}$, there is $\bfb'\in\begin{pmatrix}\bfc \\ \bfb
\end{pmatrix}$ which is homogeneous, i.e., $\begin{pmatrix}\bfb' \\ \bfa
\end{pmatrix}$ is monochromatic.  We say that $\k$ satisfies the {\it
Ramsey property} if for every $\bfa\leq\bfb$ in $\k$ and $k\geq 2$, there
is $\bfc\in\k$ with $\bfb\leq\bfc$ such that $\bfc\rightarrow (\bfb)^\bfa_k$.
For example, the classical finite Ramsey theorem is equivalent to the
statement that the class of finite linear orderings has the Ramsey
property.  Also Ne\v set\v ril and R\"odl showed that the class of 
finite ordered graphs has the Ramsey property, and Graham-Leeb-Rothschild
showed that the class of finite-dimensional vector spaces over a finite
field has the Ramsey property.

\medskip
{\bf (E)} Consider now automorphism groups Aut$(\bfa )$ of countably
infinite structures $\bfa$, for which we may as well assume that $A=\bbN$.
Thus, with the pointwise convergence topology, Aut$(\bfa )$ is a
closed subgroup of $S_\infty$, the infinite symmetric group.  Conversely,
given a closed subgroup $G\leq S_\infty$, $G$ is the automorphism group
of some structure on $A=\bbN$ (in some signature).

%
%
%

\medskip
Assume now $L$ is a signature containing a distinguished binary relation
symbol $<$.  An {\it order structure} $\bfa$ for $L$ is a structure $\bfa$ for which
$<^\bfa$ is a linear ordering.  An {\it order class} $\k$ for $L$ is one
for which all $\bfa\in\k$ are order structures.

%

We obtain the following result (Theorem \ref{4.7}):
\medskip

$\blacklozenge$ {\it The extremely amenable closed subgroups of
$S_\infty$ are exactly the groups of the form {\rm Aut}$(\bfa )$, where
$\bfa$ is the Fra\"{\i}ss\'{e} limit of a Fra\"{\i}ss\'{e} order class 
with the Ramsey property.}

\medskip
Another way to formulate this result is by saying that {\it the group 
{\rm Aut}$(\bfa )$ of automorphisms of the Fra\"{\i}ss\'{e} limit of
a Fra\"{\i}ss\'{e} order class $\k$ is
extremely amenable if and only if $\k$ has the Ramsey property.}

%
%
%

\medskip
{\bf (F)} We can now use this, and known results of structural Ramsey theory, to find
many new examples of extremely amenable automorphism groups (see Section 6).  Notice that,
by the preceding result, the extreme amenability of these groups is in fact
equivalent to the corresponding Ramsey theorem.

Consider the class of finite ordered graphs.  Its Fra\"{\i}ss\'{e} limit is
the random graph with an appropriate linear ordering.  We call it the
{\it random ordered graph}.  Let $K_n$ be the complete graph with $n$
elements, $n=3,4,\dots$.  Consider the class of $K_n$-free finite ordered
graphs, whose Fra\"{\i}ss\'{e} limit we call the {\it random} $K_n$-{\it free
ordered graph}. Finally consider the class of finite linear orderings.  Its
Fra\"{\i}ss\'{e} limit is $\langle\bbQ ,<\rangle$.

All of the above classes satisfy the Ramsey property.  This is due to
Ne\v set\v ril-R\"{o}dl \cite{NR77},\cite{NR83}
(see also Ne\v set\v ril \cite{N89} and \cite{N95}) for the graph cases, and it is of course
the classical Finite Ramsey Theorem for the last case.  So all the
corresponding automorphism groups of their Fra\"{\i}ss\'{e} limits are extremely
amenable.

This can be generalized to hypergraphs.  Let $L_0 = \{R_i\}_{i\in I}$ be a
finite relational signature with the arity of each $R_i$ at least 2.  
A {\it hypergraph of type} $L_0$ is a 
structure $\bfa_0=\langle A_0,\{R^{\bfa_0}_i\}_{i\in I}\rangle$ in which
$(a_1,\dots a_{n(i)})\in R^\bfa_i\Rightarrow a_1,\dots ,a_{n(i)}$ are
distinct, and $R^{\bfa_0}_i$ is closed under permutations.  Thus, essentially,
$R^{\bfa_0}_i\subseteq [A_0]^{n(i)}=$ the set of subsets of $A_0$ of
cardinality $n(i)$.  Consider the class of all finite ordered hypergraphs
of type $L_0$, whose Fra\"{\i}ss\'{e} limit we call the {\it random ordered
hypergraph of type} $L_0$.  More generally, for every class $\a$ of finite
irreducible hypergraphs of type $L_0$ (where $\bfa_0$ is {\it irreducible}
if it has at least two elements and for every $x\neq y$ in $A_0$ there is
$i\in I$ with $\{x,y\}\subseteq R^{\bfa_0}_i$), let 
${\mathcal OForb} (\a )$ be the
class of all finite ordered hypergraphs of type $L_0$ which omit $\a$ (i.e.,
no element of $\a$ can be embedded in them).  We call the Fra\"{\i}ss\'{e} limit
of ${\mathcal OForb}
 (\a )$ the {\it random $\a$-free ordered hypergraph of type} $L_0$.
Again Ne\v set\v ril-R\"odl \cite{NR77}, \cite{NR83}, and independently Abramson--Harrington
\cite{AH78} for $\a =\emptyset$,
showed that these classes have the
Ramsey property, so the corresponding automorphism groups are extremely
amenable.

There are similar results for metric spaces. Consider the class of 
finite ordered metric spaces with rational distances. 
Its Fra\"{\i}ss\'{e} limit is the so-called rational Urysohn space 
with an appropriate ordering. We call it the {\it ordered rational 
Urysohn space}. In response to an inquiry by the authors, 
Ne\v set\v ril \cite{N03} verified that the class of finite ordered 
metric spaces has the Ramsey property. Thus the automorphism 
group of the ordered rational Urysohn space is extremely amenable. 
We also show how this result can be used to give a new proof of 
the result of Pestov \cite{P02} that the isometry group of the Urysohn 
space is extremely amenable.

We next consider some other kinds of examples.  We first look at the
class of all finite {\it convexly ordered equivalence relations}, where
convexly ordered means that each equivalence class is convex (whenever
two elements are in it every element between them is also in it).  Their
Fra\"{\i}ss\'{e} limit is the rationals with the usual ordering and an equivalence
relation whose classes are convex, order isomorphic to the rationals,
and moreover the set of classes itself is ordered like the rationals.  We
show that the automorphism group of this structure is extremely amenable.
This implies that the corresponding class has the Ramsey property, a fact
that can also be proved directly.

Further we consider finite-dimensional vector spaces over a fixed finite
field $F$.  A {\it natural ordering} on such a vector space is one induced
antilexicographically by an ordering of a basis.  These were considered
in Thomas \cite{Th}, who showed that the class of naturally ordered
finite-dimensional spaces over $F$ is a Fra\"{\i}ss\'{e} class.  We call its
limit the $\aleph_0$-{\it dimensional vector space over $F$ with the canonical
ordering}.  The Ramsey property for the class of naturally ordered
finite-dimensional vector spaces over $F$ is easily seen to be equivalent
to the Ramsey property for the class of finite-dimensional vector spaces
over $F$, which was established in \cite{GLR}.  So the
corresponding automorphism group of the Fra\"{\i}ss\'{e} limit is extremely
amenable.

Finally, we consider the class of naturally ordered finite Boolean algebras,
where a {\it natural ordering} on a finite Boolean algebra is one
antilexicographically induced by an ordering of its atoms.  By analogy
with Thomas' result, we show that this is also a Fra\"{\i}ss\'{e} class, and
we call its limit the {\it countable atomless Boolean algebra with the 
canonical ordering}.  The Ramsey property for the class of naturally
ordered finite Boolean algebras is again easily seen to be equivalent
to the Ramsey property for the class of finite Boolean algebras and this
is trivially equivalent to the Dual Ramsey Theorem of Graham-Rothschild
\cite{GR}.  Thus the corresponding automorphism group is extremely amenable.
\medskip

{\bf (G)} Finally we use the results in {\bf (E)}, and some additional
considerations, to compute universal minimal flows.  In {\bf (E)} we have
seen a host of examples of Fra\"{\i}ss\'{e} order classes $\k$ in a signature
$L\supseteq\{<\}$.  Let $L_0=L\setminus\{ < \}$, the signature without the
distinguished symbol for the ordering.  For any structure $\bfa$ for $L$,
we denote by $\bfa_0=\bfa |L_0$ its {\it reduct} to $L_0$, i.e., $\bfa_0$
is the structure $\bfa$ with $<^\bfa$ dropped.  Denote also by $\k_0=\k
|L$ the class of all reducts $\bfa_0=\bfa |L_0$ for $\bfa\in\k$.  When
$\k$ satisfies a mild (and easily verified in every case we are
interested in) condition, in which case we call $\k$ {\it reasonable}
(see \ref{5.1} below for the precise definition), then $\k_0$ is a Fra\"{\i}ss\'{e}
class, whose limit is the reduct of the Fra\"{\i}ss\'{e} limit of $\k$.  Put
$\bff_0={\rm Flim} (\k_0),\ \bff ={\rm Flim} (\k )$, so that $\bff_0=
\bff |L_0$, i.e., $\bff =\langle\bff_0,<^\bff \rangle$.  In particular,
$F_0=F$.  Put $<^\bff =\prec_0$.  It is natural now to look at the action
of Aut$(\bff_0 )$ on the space of all linear orderings on $F_0$.  Denote
then by $X_\k$ the orbit closure $\overline{G\cdot \prec_0}$ of $\prec_0$
in this action.  It is easy to see that $X_\k$ is the space of all linear
orderings $\prec$ on $F_0$ which have the property that for any finite
substructure $\bfb_0$ of $\bff_0,\ \bfb =\langle\bfb_0,\prec |B_0\rangle
\in\k$.  We call these $\k$-{\it admissible} orderings.  This is clearly a compact
Aut$(\bff_0)$-invariant subset of $2^{F_0\times F_0}$ in the natural
action of Aut$(\bff_0)$ on $2^{F_0\times F_0}$, so $X_\k$ is an Aut$(\bff_0)$-flow.  
If $\k$ has the Ramsey property, it turns out that
any minimal subflow of $X_\k$ is the universal minimal flow of Aut$(\bff_0)$,
and $X_\k$ is itself minimal
precisely when $\k$ additionally satisfies a natural
property called the ordering property, which also plays an important role
in structural Ramsey theory (see Ne\v set\v ril-R\"odl \cite{NR78} and
Ne\v set\v ril \cite{N95}).  We say that $\k$ satisfies the {\it ordering
property} if for every $\bfa_0\in\k_0$, there is $\bfb_0\in\k_0$ such that
for every linear ordering $\prec$ on $A_0$ and every linear ordering $\prec'$
on $B_0$, if $\bfa =\langle\bfa_0,\prec\rangle\in\k$ and $\bfb =\langle
\bfb_0,\prec'\rangle\in\k$, then $\bfa\leq\bfb$.  Now Theorems 
\ref{7.4}, \ref{7.5}, and \ref{10.8} provide a toolbox for computing 
universal minimal flows of automorphism groups, which can be summarized
as follows. 
\medskip

$\blacklozenge$ 
{\it In the above assumptions, let
$X_\k$ be the ${\rm Aut} (\bff_0)$-flow of $\k$-admissible
orderings on $F_0 \;(=F)$.  Then the following are equivalent:

(i) $X_\k$ is a minimal ${\rm Aut} (\bff_0)$-flow,

(ii) $\k$ satisfies the ordering property.

\noindent Moreover the following are equivalent:

(iii) $X_\k$ is the universal minimal ${\rm Aut} (\bff_0)$-flow.

(iv) $X_\k$ satisfies the Ramsey and ordering properties.}
\medskip

Now all the classes $\k$, considered in {\bf (F)} above, satisfy the ordering property.
This is due to Ne\v set\v ril-R\"odl \cite{NR78} for the case of graphs and
hypergraphs, Ne\v set\v ril \cite{Nes04a} for metric spaces, and is easily verified in 
all the other cases.  Therefore,
we have the following computations of universal minimal flows 
(see Section 8):
\medskip

$\blacklozenge$
{\it If $\bff_0$ is one of the following structures, then the 
universal minimal flow $M({\rm Aut} (\bff_0))$ of the group
${\rm Aut} (\bff_0)$ of automorphisms of $\bff_0$ is its action on 
the space of linear orderings on the universe $F_0$ of $\bff_0$:

(a)  The random graph.

(b) The random $\k_n$-free graph, $n=2,3,\dots$

(c)  The structure $\langle\bbN\rangle$ (where ${\rm Aut} (\bff_0)
=S_\infty$).

(d) The random hypergraph of type $L_0$.

(e) The random $\a$-free hypergraph of type $L_0$,
where $\a$ is a class of irreducible finite hypergraphs of type $L_0$.

(f) The rational Urysohn space.

\smallskip
$\blacklozenge$ If $\bff_0$ is the equivalence relation on a countable
set with infinitely many classes, each of which is infinite, then
$M({\rm Aut} (\bff_0))$ is the 
space of all linear orderings on that set for which each equivalence
class is convex.

\smallskip
$\blacklozenge$ 
If $\bff_0=\bfv_F$ is the $\aleph_0$-dimensional vector 
space $\bfv_F$ over a finite field $F$, then $M({\rm Aut} (\bff_0))$ is
the space of all orderings on $V_F$, whose
restrictions to finite-dimensional subspaces are natural.

\smallskip
$\blacklozenge$ 
If $\bff_0=\bfb_\infty$ is the countable
atomless Boolean algebra $\bfb_\infty$, then  $M({\rm Aut} (\bff_0))$
the space of all
linear orderings on $B_\infty$, whose restrictions to finite subalgebras
are natural.}

\medskip
In particular, in all these cases, the universal minimal flow is metrizable.

Of course (i), (c) is the result of Glasner-Weiss \cite{GlW02}.
Very recently, Glasner-Weiss \cite{GlW03} computed the universal minimal
flow of $H(2^\bbN )$, the homeomorphism group of the Cantor space $2^\bbN$,
as the space of maximal chains of compact subsets of $2^\bbN$, which is 
metrizable.  Since the group in (iv) above is, by Stone duality, isomorphic
to $H(2^\bbN )$, we have another proof that the universal minimal flow is
metrizable and a different description of this flow.  Of course these two
flows are isomorphic and in fact an explicit isomorphism can be found.

There is actually quite a bit more that we can say in this context and this has some further interesting connections with Ramsey theory.
Let $\k_0$ be a Fra\"\i ss\'e class in a signature $L_0$, let $L=L_0
\cup\{<\}$ and call any class $\k$ in $L$ with $\k |L_0=\k_0$ an 
{\it expansion} of $\k_0$.  We define in Section 9 a canonical notion
of equivalence of expansions of $\k_0$, called {\it simple bi-definability}.
Intuitively, if $\k' ,\k''$ are simply bi-definable expansions, we view 
$\k' ,\k''$ as ``trivial'' perturbations of each other.  Using dynamical
ideas, e.g., the uniqueness of the universal minimal flow, we prove the
following (see \ref{9.2}, \ref{10.7}):
\medskip


$\blacklozenge$
{\it If $\k$ is a reasonable Fra\"\i ss\'e order class in $L$ which
is an expansion of $\k_0$ and has the Ramsey property, then there
is $\k'\subseteq\k$, a reasonable Fra\"\i ss\'e order class in $L$
which is an expansion of $\k_0$ and has both the Ramsey and ordering
properties.

$\blacklozenge$ If $\k',\k''$ are reasonable Fra\"\i ss\'e order classes in $L$
that are expansions of $\k_0$ and satisfy both the Ramsey and ordering
properties, then $\k',\k''$ are simply bi-definable.}
\medskip

Thus among all the expansions of $\k_0$ that satisfy the Ramsey property
there are canonical ones, those that also satisfy the ordering property.
These are unique, up to simple bi-definability, and are ``extremal'',
i.e., least under inclusion, again up to simple bi-definability.
In several cases, like, e.g., Boolean algebras and vector spaces, we can
also list explicitly all such classes.

If $\k_0$ is as above and $\bfa_0\in\k_0$, let
$t(\bfa_0,\k_0)$ be the smallest
number $t\in\bbN$, if it exists, that satisfies:  For every $\bfa_0
\leq\bfb_0$ in $\k_0, k\geq 2$, there is $\bfc_0\in\k_0$ with $\bfc_0\geq
\bfb_0$ and
\[
\bfc_0\rightarrow (\bfb_0)^{\bfa_0}_{k,t},
\]
where this notation means that for every coloring $c:\begin{pmatrix}
\bfc_0\\ \bfa_0\end{pmatrix}\rightarrow\{1,\dots ,k\}$, there is
$\bfb'_0\in\begin{pmatrix}\bfc_0\\ \bfb_0\end{pmatrix}$ such that
$c$ on $\begin{pmatrix}\bfb'_0\\ \bfa_0\end{pmatrix}$ takes at most
$t$ many values.  (Thus $\bfc_0\rightarrow (\bfb_0)^{\bfa_0}_k
\Leftrightarrow\bfc_0\rightarrow (\bfb_0)^{\bfa_0}_{k,1}$.)  
This variation of the original arrow-notation of Erd\"{o}s-Rado \cite{ER56} is due to Erd\"{o}s-Hajnal-Rado \cite{EHR65} and has already reappeared in several other areas of Ramsey theory. Following Fouch\'e \cite{Fo98}, we call $t(\bfa_0,\k_0)$ the {\it Ramsey
degree} of $\bfa_0$ in $\k_0$.
For every expansion
$\k$ of $\k_0$ let, for $\bfa_0\in\k_0$
\[
X^{\bfa_0}_\k =\{\prec:\prec\text{ is a linear ordering on }A_0\text{ and }
\langle\bfa_0,\prec\rangle\in\k\}.
\]
Then Aut$(\bfa_0)$ acts in the obvious way on $X^{\bfa_0}_\k$ and we let
$t_\k (\bfa_0)$ be the number of orbits.  Clearly,
\[
t_\k (\bfa_0)=\frac{{\rm card}(X^{\bfa_0}_\k )}{{\rm card (Aut}(\bfa_0))}.
\]
For example, if $\k_0=$ the class of finite graphs, $\k =$ the class
of finite ordered graphs, then $t_\k (\bfa_0)=\tfrac{{\rm card}(\bfa_0)!}
{{\rm card (Aut}(\bfa_0))}$.  The following can be proved by standard methods in Ramsey theory (see Section 10).
\medskip

$\blacklozenge$ {\it Let $\k_0$ be a Fra\"\i ss\'e class
in a signature $L_0, L=L_0\cup\{<\}$.  If $\k_0$ admits a reasonable
Fra\"\i ss\'e order class in $L$ which is an expansion of $\k_0$ and
has the Ramsey property, then $t(\bfa_0,\k_0)$ exists for all $\bfa_0\in\k_0$.
Moreover, $t(\bfa_0,\k_0)\leq t_\k (\bfa_0)$, for any reasonable Fra\"\i ss\'e
order class $\k$ in $L$ which is an expansion of $\k_0$ and has the Ramsey
property, and $t(\bfa_0,\k_0)=t_{\k'}(\bfa_0)$ for any such class 
$\k'$ that also has the ordering property.}
\medskip

For example, if $\k_0=$ the class of finite graphs, then the Ramsey-theoretic
results of Ne\v{s}et\v{r}il-R\"odl \cite{NR77}, \cite{NR78}, \cite{NR83} 
and Abramson-Harrington
[78] can be interpreted as saying that $t(\bfa_0,\k_0)=\tfrac{{\rm card}
(A_0)!}{{\rm card(Aut}(\bfa_0))}$.  This seems to have been first pointed
out in print by Fouch\'e \cite{Fo99}, who in a series of other papers 
\cite{Fo97}, \cite{Fo99}, \cite{Fo99a}
has computed Ramsey degrees of several classes of finite structures.

\medskip
{\it Acknowledgments.} A.S.K. was partially supported by NSF Grant 
DMS-9987437, the Fields Institute, Toronto, the Centre de Recerca
Matem\` atica, Bellatera, 
and a Guggenheim Fellowship.  V.G.P. was partially supported by the
Marsden Fund of the Royal Society of New Zealand, by a University of
Ottawa start-up grant, and by an NSERC operating grant and would like to
thank Vitali Milman and Pierre de la Harpe for stimulating discussions and the Caltech Department
of Mathematics for hospitality extended in October 2001 and April 2003.
S.T. would like to acknowledge the support of the CNRS, Paris,
the Centre de Recerca Matem\` atica, Bellatera, 
and the Fields Institute, Toronto.  We would also like to thank G. Debs,
J. Ne\v{s}et\v{r}il and L. Nguyen Van The for helpful comments and suggestions.

\section{Topological dynamics}  {\bf (A)} We will survey here some basic
concepts and results of topological dynamics, which we will need in this
paper.  More detailed treatments can be found in Ellis \cite{E69}, Auslander
\cite{A}, de Vries \cite{dV}, Pestov \cite{P98}, \cite{P99}, 
Glasner \cite{Gl00}, Uspenskij \cite{Usp02}.

Recall that an action $(g,x)\in G\times X\mapsto g\cdot x\in X$ of a
topological group $G$ on a topological space $X$ is {\it continuous} if 
it is continuous as a map from $G\times X$ into $X$.  We will consider
continuous actions of (Hausdorff) topological groups $G$  on (non-$\emptyset$)
compact, Hausdorff spaces $X$.  Actually we are primarily interested in
metrizable topological groups $G$ and in fact only separable metrizable
ones.  So although we will state in this survey several standard results
for general topological groups, we will often give a sketch of an argument
for the metrizable case only.  It should be kept in mind that if $G$ is
metrizable (equivalently is Hausdorff and has a countable nbhd basis
at the identity), then $G$ admits a right-invariant compatible metric
$d_r$, which of course can always be taken to be bounded by 1, by replacing
it, if necessary, by $\tfrac{d_r}{1+d_r}$.  See, e.g., \cite{Be}, p. 28.

Let $G$ be a topological group and $X$ a compact, Hausdorff space. If we equip $H(X)$, the group of homeomorphisms
of $X$, with the compact-open topology, i.e., the topology with subbasis
$\{f\in H(X):f(K)\subseteq V\}$, with $K\subseteq X$ compact, $V\subseteq
X$ open, then $H(X)$ is a topological group, and a continuous action of $G$
on $X$ is simply a continuous homomorphism of $G$ into $H(X)$.  We will
also refer to a continuous action of $G$ on $X$ as a $G$-{\it flow} on $X$.
If the action is understood, we will often simply use $X$ to refer to the 
flow.

Given a $G$-flow on $X$ and a point $x\in X$, the {\it orbit} of $x$ is the
set
\[
G\cdot x=\{g\cdot x:g\in G\}
\]
and the {\it orbit closure} of $x$, the set
\[
\overline{G\cdot x}.
\]
This is a $G$-invariant, compact subset of $X$.  In general, a (non-$\emptyset$)
compact, $G$-invariant subset $Y\subseteq X$ defines a {\it subflow} by
restricting the $G$-action to $Y$.  A $G$-flow on $X$ is {\it minimal} if
it contains no proper subflows, i.e., there is no (non-$\emptyset$) compact
$G$-invariant set other than $X$.  Thus $X$ is minimal iff every orbit is
dense.  A simple application of Zorn's Lemma shows that every $G$-flow
$X$ contains a minimal subflow $Y\subseteq X$.

Among minimal flows of a given group $G$, there is a largest (universal)
one, called the {\it universal minimal flow}.  To define this, we first
need the concept of homomorphism of $G$-flows.  Let $X,Y$ be two $G$-flows.
A {\it homomorphism} of the $G$-flow $X$ to the $G$-flow $Y$ is a continuous
map $\pi :X\rightarrow Y$, which is also a $G$-{\it map}, i.e.,
\[
\pi (g\cdot x)=g\cdot \pi (x),\ \ \ x\in X, g\in G.
\]
Notice that if $Y$ is minimal, then any homomorphism of $X$ into $Y$ is
surjective.  An {\it isomorphism} of $X$ to $Y$ is a bijective homomorphism
$\pi :X\rightarrow Y$ (notice then that $\pi^{-1}$ is also a homomorphism).
We now have the following basic fact in topological dynamics.  (For a 
proof see Auslander \cite{A}, Ch. 8, or Uspenskij \cite{Usp02}, \S 3.)

\begin{theorem}
\label{1.1} 
Given a topological group $G$, there is a
minimal $G$-flow $M(G)$ with the following property:  For any minimal
$G$-flow $X$ there is a homomorphism $\pi :M(G)\rightarrow X$.  Moreover,
$M(G)$ is uniquely determined up to isomorphism by this property. \qed
\end{theorem}

The space $M(G)$ is called the {\it universal minimal flow} of $G$.  In order
to get an intuition about this space, we will discuss a standard way of
looking at it.  For that we first need to discuss the concept of a {\it
pointed} $G$-{\it flow} or {\it ambit}.

Let $G$ be a topological group.  A $G$-{\it ambit} is a $G$-flow $X$ with a
distinguished point $x_0\in X$, whose orbit is dense in $X$.  We often
abbreviate this by $(X,x_0)$.  A {\it homomorphism of} $G$-{\it ambits}
$(X,x_0), (Y,y_0)$ is a homomorphism $\pi :X\rightarrow Y$ of the $G$-flows
such that $\pi (x_0)=y_0$. If such a homomorphism exists, it is clearly
unique.  Similarly we define the concept of isomorphism of $G$-ambits.

It is another basic fact of topological dynamics that there is again
a largest (universal) $G$-ambit.

\begin{theorem}
\label{1.2}
Given a topological group $G$, there is a
$G$-ambit ($S(G), s_0$) with the following property:  For any $G$-ambit
$(X,x_0)$ there is a homomorphism of $(S(G), s_0)$ to $(X,x_0)$.  Moreover,
$(S(G), s_0)$ is uniquely determined up to isomorphism by this property.
\qed
\end{theorem}

The space $(S(G),s_0)$, often simply written as $S(G)$, is called the {\it
greatest G-ambit}.

The uniqueness part of the preceding result is obvious from the definitions.
To establish existence, we will describe a particularly useful way of
constructing the greatest ambit.

Consider the space $RUC^b(G)$ of bounded right-uniformly continuous functions $x:G\rightarrow\bbC$.
Recall that $x:G\rightarrow\bbC$ is {\it right-uniformly continuous} if for
each $\epsilon >0$ there is a nbhd $V$ of the identity $1_G$ of $G$ so that
\[
gh^{-1}\in V\Rightarrow |x(g)-x(h)|<\epsilon .
\]
If $G$ is metrizable, with right-invariant compatible metric $d_r$, then
right-uniformly continuous means uniformly continuous with respect to $d_r$:
\[
\forall\epsilon\exists\delta (d_r(g,h)<\delta\Rightarrow |x(g)-x(h)|<
\epsilon ).
\]
Under pointwise addition, multiplication and conjugation, and with the
sup norm $\Vert x\Vert_\infty =\sup\{|x(g)|:g\in G\}, RUC^b(G)$ is an
abelian $C^*$-algebra which is unital (with multiplicative identity
the constant 1 function).  Denote by $S(G)$ the maximal ideal space of
$RUC^b(G)$, i.e., the space of all continuous homomorphisms $\varphi
:RUC^b(G)\rightarrow\bbC$.  Equipped with the topology generated by
the maps $\hat x:S(G)\rightarrow\bbC,\hat x(\varphi )=\varphi (x)$, for
$x\in RUC^b(G)$ (i.e., the smallest topology in which all $\hat x,
x\in RUC^b(G)$, are continuous), this is a compact, Hausdorff space.
Moreover, by the Gelfand-Naimark theorem, $RUC^b(G)$ can be canonically 
identified with $C(S(G))$, the $C^*$-algebra of all continuous
complex-valued functions on $S(G)$, identifying $x\in RUC^b(G)$
with $\hat x$ (see, e.g., \cite{Ru}, 11.18).

Now $G$ acts continously by $C^*$-automorphisms on $RUC^b(G)$ by 
left-shift
\[
g\cdot x(h) =x(g^{-1}h)
\]
and thus acts canonically on $S(G)$ via
\[
g\cdot\varphi (x)=\varphi (g^{-1}\cdot x).
\]
It is also easy to check that this action is continuous, so $S(G)$ is a
$G$-flow.  We will now identify a canonical element of $S(G)$, that will
turn it into an ambit.

For each $g\in G$, let $\varphi_g\in S(G)$ be defined by $\varphi_g
(x)=x(g),x\in RUC^b(G)$.  Then one can see that $g\mapsto\varphi_g$
is a homeomorphism of $G$ with a dense subset of $S(G)$.  For example,
when $G$ is metrizable with bounded compatible right-invariant metric
$d_r$, and $g_0\neq h_0$, then for $x(g)=d_r(g,h_0),\varphi_{g_0}(x)\neq
\varphi_{h_0}(x)$ so $\varphi_{g_0}\neq\varphi_{h_0}$, i.e., this map is
1-1.  The verification that it is homeomorphism is straightforward.
Finally $\{\varphi_g:g\in G\}$ is dense in $S(G)$, since, otherwise,
there is $f\in C(S(G))$, so that $f=\hat x$ for some $x\in RUC^b(G)$,
with $f\neq 0$ but $f(\varphi_g)=\hat x(\varphi_g)=x(g)=0,\forall g\in G$,
which implies that $x=0$, thus $f=0$, a contradiction.

So from now on we will identify $g$ with $\varphi_g$ and think of $G$ as
a dense subset of $S(G)$.  Moreover $G$ is an invariant subset of 
$S(G)$ and the restriction of the action to $G$ is simply left-translation:
$(g,h)\mapsto gh$.  We now have the following standard fact.

\begin{theorem}
\label{1.3}
The $G$-ambit $(S(G), 1_G)$ is the
greatest $G$-ambit.
\end{theorem}

\begin{proof}
Since the orbit of $1_G$ in $S(G)$ is $G$, which is dense,
clearly $(S(G),1_G)$ is a $G$-ambit.  Consider now an arbitrary $G$-ambit
$(X,x_0)$.  Suppose $f\in C(X)$.  Define then $f^*:G\rightarrow\bbC$ by
\[
f^*(g)=f(g\cdot x_0).
\]
We verify that $f^*\in RUC^b(G)$.  Since the action of $G$ on $X$
is continuous, an easy compactness argument shows that given $\epsilon >0$,
there is a nbhd $V$ of the identity of $G$ such that $g\in V$ implies
$|f(g\cdot x)-f(x)|<\epsilon ,\forall x\in X$.  So if $gh^{-1}\in V$, then
$|f^*(g)-f^*(h)|=|f(g\cdot x_0)-f(h\cdot x_0)|=|f(gh^{-1}\cdot (h\cdot x_0))
-f(h\cdot x_0)|<\epsilon$, and thus $f^*\in RUC^b(G)\ (f^*$ is clearly
bounded).

Identifying, as usual, $RUC^b(G)$ with $C(S(G))$, the map $f\mapsto
f^*$ is a unital $C^*$-algebra monomorphism of $C(X)$ into $C(S(G))$.
Now it is a well known fact that a unital $C^*$-algebra monomorphism
$\pi :C(K)\rightarrow C(L)$, where $K,L$ are (non-$\emptyset$) compact
spaces, is of the form $\pi (f)=f\circ\Pi$ for a {\it uniquely}
determined continuous surjection $\Pi :L\rightarrow K$ (see, e.g.,
\cite{Co}, 2.4.3.6).  From this it also follows that if $K,L$
are actually $G$-flows, and we let $G$ act on $C(K),C(L)$ by shift,
$g\cdot f(x)=f(g^{-1}\cdot x)$, then if $\pi$ is a $G$-map, $\Pi$ is
also a $G$-map.  Applying this to $f\mapsto f^*$, we see that there is a
homomorphism of $G$-flows $\Phi :S(G)\rightarrow X$ with $f^*=f\circ\Pi$.
It only remains to check that $\Pi (1_G)=x_0$.  But for any $f\in C(X),
f^*(1_G)=f(x_0)=f(\Pi (1_G))$, so we must have $\Pi (1_G)=x_0$, and
the proof is complete.\hfill$\dashv$
\end{proof}

Using this it is now immediate to obtain the following description of the
universal minimal flow of $G$.

\begin{corollary}
\label{1.4}
Let $M(G)$ be a minimal subflow of
$S(G)$ (i.e., $M(G)$ is a minimal $G$-invariant compact subset of $S(G)$).
Then $M(G)$ is the universal minimal flow (up to isomorphism).
\end{corollary}

\begin{proof}
Let $X$ by any minimal $G$-flow.  Fix $x_0\in X$.  Then
$(X,x_0)$ is a $G$-ambit, so let $\pi :(S(G), 1_G)\rightarrow (X,x_0)$
be a homomorphism.  Then clearly the restriction of $\pi$ to $M(G)$ is also
a homomorphism, and we are done.
\end{proof}

In particular, it follows from the uniqueness part of Theorem \ref{1.1}, that all
minimal subflows of $S(G)$ are isomorphic.  (This uniqueness part, which
is not proved here, is based on techniques from semigroup theory.)
As we will soon see, the space $M(G)$ can be extremely complicated, e.g.,
non-metrizable, even when the group $G$ is very ``small'', e.g.,
an infinite 
countable discrete $G$.  However, we will verify that when $G$ is
separable, metrizable, $M(G)$ is at least an inverse limit of metrizable
$G$-flows.

Fix a topological group $G$.  An {\it inverse system} of $G$-flows consists
of a directed set $\langle I,\preceq\rangle$, a family $\{X_i\}_{i\in I}$
of $G$-flows and a family of homomorphisms $\pi_{ij}:X_j\rightarrow X_i$,
for each $i\preceq j$, such that $\pi_{ii}=$ the identity of $X_i$, and
$i\preceq j\preceq k\Rightarrow\pi_{ik}=\pi_{ij}\circ\pi_{jk}$.  The
{\it inverse limit} $\lim_{\leftarrow}X_i$ is the $G$-flow defined as
follows:  Consider the product topological space $\prod_{i\in I}X_i$ and
let
\[
\lim_{\leftarrow}X_i=\{\{x_i\}_i\in\prod_iX_i:\forall i\preceq j(\pi_{ij}
(x_j)=x_i)\}.
\]
By a simple application of compactness, $\lim_{\leftarrow}X_i\neq\emptyset$ and is 
clearly a compact subset of $\prod_{i\in I}X_i$.  The group $G$ acts on
$\lim_{\leftarrow}X_i$ coordinatewise:  $g\cdot\{x_i\}=\{g\cdot x_i\}$ and
this is clearly a continuous action.  Define
\[
\pi_i:\lim_{\leftarrow}X_i\rightarrow X_i
\]
by $\pi_i(\{x_i\}_{i\in I})=x_i$.  Then $\pi_i$ is a homomorphism and
if $i\preceq j$, then $\pi_i=\pi_{ij}\circ\pi_j$.  Finally, if $X$ is a
$G$-flow and there are homomorphisms $\varphi_i:X\rightarrow X_i$ with
$i\preceq j\Rightarrow\varphi_i=\varphi_{ij}\circ\varphi_j$, then there is
a unique homomorphism $\varphi :X\rightarrow\lim_{\leftarrow}X_i$ such
that $\varphi_i=\pi_i\circ\varphi$.

Similarly we define inverse systems of $G$-ambits.  We now have the following
folklore fact.

\begin{theorem}
\label{1.5}
Let $G$ be a separable, metrizable group.
Then the greatest ambit $(S(G),1_G)$ is the inverse limit of a system of
metrizable $G$-ambits.  Similarly, the universal minimal flow is the inverse
limit of a system of metrizable minimal $G$-flows.
\end{theorem}

\begin{proof} We will first derive the second assertion from the first.
Suppose $\{(X_i,x^0_i)\}$ is an inverse system of metrizable $G$-ambits,
so that $(X,x_0)=\lim_{\leftarrow}(X_i,x^0_i)$ is the greatest $G$-ambit.
Let $\pi_i:X\rightarrow X_i$ be the corresponding homomorphism.  In
particular, $X=\lim_{\leftarrow}X_i$ as a $G$-flow.  Now fix a minimal
subflow $M\subseteq X$, so that, as we have seen earlier, $M$ is the 
universal minimal flow.  Put $\pi_i(M)=M_i\subseteq X_i$.  Then $M_i$ is a
subflow of $X_i$ and it is clearly minimal and metrizable.  So it is enough
to check that $M=\lim_{\leftarrow}M_i$, which is an easy compactness
argument.

For the first assertion, fix a countable dense set $D\subseteq G$.  Let
$(I,\preceq)$ be the following directed set:  $I$ consists of all separable,
closed, unital $G$-invariant $C^*$-subalgebras of $RUC^b(G)$.  Since a 
closed, unital $C^*$-subalgebra of $RUC^b(G)$ is $G$-invariant iff it
is $D$-invariant, clearly $\bigcup I=RUC^b(G)$.  Let for $A,B\in I$
\[
A\preceq B\Leftrightarrow A\subseteq B.
\]
For $A\in I$ denote by $X_A$ the maximal ideal space of $A$, which is
compact metrizable, since $A$ is separable.  Since $G$ acts continuously
by $C^*$-automorphisms on $A$, it acts continuously on $X_A$.  Moreover,
as before,  we can identify each $g\in G$ with an element of $X_A$, so that
$G$ can be thought as a dense invariant subset of $X_A$ and the $G$ action
on it is by left-translation.  Thus again $(X_A,1_G)$ is a metrizable
$G$-ambit.  We view as usual $A$ as identified with $C(X_A)$ via $x\mapsto
\hat x_A$.  When $A\preceq B$, the identity is an injective unital 
$C^*$-homomorphism from $A$ to $B$, so there is a unique surjection 
$\pi_{AB}:X_B\rightarrow X_A$ such that for $x\in A$, $\hat x_B=\hat x_A
\circ\pi_{AB}$, therefore for $\varphi\in X_B, \hat x_A(\pi_{AB}(\varphi
))=\hat x_B(\varphi )$ or $\pi_{AB}(\varphi )(x)=\varphi (x)$, i.e.,
$\pi_{AB}(\varphi )=\varphi |A$.  Similarly, there is a surjection
$\pi_A:S(G)\rightarrow X_A$ given by $\pi_A(\varphi )=\varphi |A$ for
each $A\in I$, so that $\pi_A=\pi_{AB}\circ\pi_B$ for $A\preceq B$.
Moreover, $\pi_A (1_G)= \pi_{AB} (1_G) = 1_G$.  Thus 
$\varphi\mapsto (\pi_A(\varphi))_{A\in I}$
is a homomorphism from $(S(G),1_G)$ to $\lim_{\leftarrow}(X_A,1_G)$.
Also given $(\varphi_A)\in\lim_{\leftarrow}X_A$, we have, for
$A\preceq B$, that $\varphi_A=\varphi_B|A$, so that there is a unique
$\varphi\in S(G)$ with $\pi_A(\varphi )=\varphi_A$.  Thus $\varphi\mapsto
(\pi_A(\varphi ))_{A\in I}$ is an isomorphism of the $G$-ambit
$(S(G),1_G)$ and $\lim_{\leftarrow} (X_A,1_G)$.
\end{proof}

{\bf (B)}  We will now discuss the case of infinite countable discrete
groups $G$ and see that in this case $M(G)$ is an extremely large space,
in particular it is not metrizable.

For an infinite countable discrete $G$, it is clear that $RUC^b(G)$ is
identical with $\ell^\infty (G)$, the $C^*$-algebra of bounded complex
functions on $G$ with the supremum norm.  It is then easy to see that 
$S(G)$ is identical with $\beta G$, the space of ultrafilters on $G$ with
the topology whose basis consists of the sets $\hat A=\{U\in\beta G:
A\in U\}$, for $\emptyset\neq A\subseteq G$.  This is a non-metrizable,
compact, Hausdorff space.  The action of $G$ on $S(G)$ is given by
\[
A\in g\cdot U\Leftrightarrow g^{-1}A\in U
\]
for $A\subseteq G$.

The copy of $G$ in $S(G)$ consists simply of the principal ultrafilters.
So the distinguished point is the principal ultrafilter on $1_G$.  Finally, $M(G)$ is any minimal subflow
of $\beta G$.

First we point out that $G$ acts {\it freely} on $S(G)$, i.e., $g\cdot
x\neq x,\ \forall g\neq 1_G,\ x\in S(G)$ (this is due to Ellis \cite{E60}).
For that it is of course enough to find a free $G$-flow $X$
(since $S(G)$ can be homomorphically mapped to $X$).  Again it is enough
to find for each $g\in G,\ g\neq 1_G$, a $G$-flow $X_g$, such
that $g\cdot x\neq x,\ \forall x\in X_g$.  Because then $X=\prod_{g\in G
\setminus\{1_G\}}X_g$ with the action $h\cdot (x_g)=(h\cdot x_g)$ works.
Consider the shift action of $G$ on $3^G,\ g\cdot p(h)=p(g^{-1}h)$. 
Suppose we can find $p_g\in 3^G$ such that $(*):\forall h\in G(g\cdot
(h\cdot p_g)(1_G)\neq (h\cdot p_g)(1_G))$.  Then we can take $X_g$ to be
the orbit closure of $p_g$.  Now $(*)$ is equivalent to:  $\forall h\in G
(p_g(h^{-1}g^{-1})\neq p_g(h^{-1}))$ or $\forall h(p_g(h)\neq p_g(hg))$.
Consider a coset $h\langle g\rangle$ of the cyclic subgroup $\langle g
\rangle$.  Define $p_g$ on $h\langle g\rangle$, so that $p_g(hg^k)\neq
p_g(hg^{k+1})$, for any $k\in\bbZ$.  (We used $3^G$ instead of $2^G$
to take care of the case when $g$ has finite order.)  This clearly works.

The space $M(G)$ is quite big, see e.g., the references in \cite{dV}, 
p. 391, {\bf 11}.  Let us verify for instance that it is not metrizable.
The space $M(G)$ is a closed subset of $\beta G$.  If it was metrizable
and infinite, it would have a non-trivial convergent sequence, which is
impossible in the extremally disconnected space $\beta G$ (see \cite{Eng},
Ex. 6.2.G(a) on p. 456).  So if $M(G)$ is metrizable, it has to be
finite, contradicting the fact that $G$ acts freely on $M(G)$.
\medskip

{\bf (C)} More generally, Veech \cite{V} has shown that when $G$ 
is locally compact, then
$G$ acts freely on $S(G)$ and thus on $M(G)$.  For a simpler version of
the proof see Pym \cite{Py}.  See also \cite{AS} for the second countable
case.  In an Appendix to this paper, we also give a new proof of Veech's
Theorem.  Note that Veech's Theorem implies that if $G$ 
is second countable, then $G$
admits a free metrizable $G$-flow.  To see this, notice that it is enough
to find for each $1_G\neq g\in G$ a metrizable $G$-flow $X_g$ with
$g\cdot x\neq x$, if $x\in X_g$.  Indeed, if we can do that, then, by 
compactness, there is an open nbhd $V_g$ of $g$ with $1_G\not\in V_g$ and $h\in
V_g\Rightarrow (h\cdot x\neq x,\ \forall x\in X_g)$.  Find now $g_0,g_1,
\dots \in G\setminus\{1_G\}$ so that $\{V_{g_n}\}_{n\in\bbN}$ is an open
cover of $G\setminus\{1_G\}$.  Then $X=\prod_nX_{g_n}$ with the
coordinatewise action is a free metrizable $G$-flow. (This argument comes
from \cite{AS}.)
So fix $g\in G,\ g\in 1_G$, in order to find $X_g$.  Write $S(G)=
\lim_{\leftarrow}X_i,\ X_i$ a metrizable $G$-flow.  If none of the
$X_i$ can be $X_g,\ \{x_i\in X_i:g\cdot x_i=x_i\}=Y_i$ is non-$\emptyset$,
and if $\pi_{ij}:X_j\rightarrow X_i$, for $i\preceq j$, are the
corresponding homomorphisms, then $\pi_{ij}(Y_j)\subseteq Y_i$, so the
inverse limit of $(Y_i,\pi_{ij}|Y_j)$ is non-$\emptyset$ and thus there 
is $\{y_i\}\in S(G)$ with $g\cdot \{y_i\}=\{y_i\}$, a contradiction, as $G$
acts freely on $S(G)$.

We also show in an Appendix that when $G$ is non-compact, locally compact,
$M(G)$ is non-metrizable.  Stronger results in
special cases were obtained in Turek \cite{Tu}, Lau-Milner-Pym \cite{LMP}.  Of
course when $G$ is compact, $M(G)$ is $G$ itself with the left-translation
action.
\medskip

{\bf (D)}  Rather remarkably, there are nontrivial groups $G$ for which $M(G)$ 
trivializes, i.e., consists of a single point.  Such groups are called 
extremely amenable.  Thus a topological group is {\it extremely
amenable} if any $G$-flow $X$ has a fixed point, i.e., there is $x\in X$
with $g\cdot x = x,\ \forall g\in G$.  (For this reason, sometimes
extremely amenable groups are described as groups having the {\it fixed
point on compacta property}.)  By Veech's Theorem nontrivial such 
groups cannot be
locally compact.  As it turned out, a number of important, non-locally
compact Polish groups are extremely amenable.  Among them are: the unitary
group $U(H)$ of the infinite dimensional separable Hilbert space $H$
(Gromov-Milman \cite{GroM}); $L(I,\mathbb T )$, the group of measurable maps from $I=
[0,1]$ to $\mathbb T$, with pointwise multiplication, and the topology of
convergence in measure (Furstenberg-Weiss, Glasner \cite{Gl98}); $H_+(I)$ and
$H_+(\bbR )$, the groups of orientation preserving homeomorphisms of
$I$ and $\bbR$, with the compact-open topology (Pestov \cite{P98a}); Aut$(I,
\lambda )$ (resp., Aut$^*(I,\lambda ))$, the groups of measure preserving
(resp., measure-class preserving) automorphisms of Lebesgue measure
$\lambda$ on $I$, with the weak topology (Giordano-Pestov \cite{GiP});
Iso($\bf U$), the group of isometries of the Urysohn space, with the
pointwise convergence topology (Pestov \cite{P02})), and Aut$(\langle\bbQ ,
<\rangle )$, the group of automorphisms of the rationals with the usual
ordering, with the pointwise convergence topology (Pestov \cite{P98a}).  For
more about extreme amenability, see also Pestov \cite{P99},
\cite{P02a}, \cite{P02b} and Uspenskij \cite{Usp02}.

In case the group $G$ is separable metrizable, we can restate the definition
of extreme amenability in terms of metrizable flows only.  In
other words, a separable metrizable group $G$ is extremely amenable iff
every metrizable $G$-flow has a fixed point.  Indeed, if every
metrizable $G$-flow has a fixed point, every minimal metrizable $G$-flow
is a singleton, and thus so is $M(G)$, being an inverse limit of such $G$-flows.
\medskip

{\bf (E)} Except for the case of compact metrizable or extremely amenable groups, there were very
few cases where the universal minimal flow $M(G)$ was known to be metrizable.
Pestov \cite{P98a} first computed that for the group $H_+(\mathbb T )$ of
orientation-preserving homeomorphisms of the circle $\mathbb T$, the canonical
evaluation action on $\mathbb T$ is the universal minimal flow.  Then 
Glasner-Weiss \cite{GlW02} computed the universal minimal flow of $S_\infty$,
the group of permutations of $\bbN$ with the pointwise convergence topology.
It turns out to be the canonical action on the compact, metrizable space of
linear orderings on $\bbN$.  Finally Glasner-Weiss \cite{GlW03} computed the universal minimal flow of the
group $H(2^\bbN )$ of the homeomorphisms of the Cantor space $2^\bbN$.
It is the action of $H(2^\bbN )$ on the space of maximal chains of compact
subsets of $2^\bbN$, invented by Uspenskij \cite{Usp00}.

Let us point out here that only one of the following is possible:

(i) $M(G)$ is finite.

(ii) $M(G)$ is perfect (i.e., has no isolated points), and thus card$(M(G))=
2^{\aleph_0}$, if $M(G)$ is metrizable.

To see this, note that if $M(G)$ has an isolated point $x_0$, then, as
every orbit of $M(G)$ is dense, $x_0$ belongs to every orbit, thus there is
only one orbit in $M(G)$, and so every point of $M(G)$ is isolated.
Since $M(G)$ is compact, $M(G)$ is finite.

We have already discussed examples of metrizable $M(G)$ which
consist of exactly one point or are perfect.  It is easy to see that also
any finite cardinality for $M(G)$ is possible.  Because if $H$ is extremely
amenable, and $G=H\times\bbZ_n$, then the obvious action of $G$ on $\bbZ_n$
is the universal minimal flow of $G$, thus card$(M(G))=n$.

Our goal in this paper is to study extreme amenability and universal
minimal flows of closed subgroups of $S_\infty$, i.e., automorphism
groups of countable structures.  In particular, we find new examples of
extremely amenable groups and also new cases where the universal minimal
flow is metrizable and can be computed.

\section{Fra\"{\i}ss\'{e} Theory}

We will review here some basic ideas of model theory concerning the Fra\"{\i}ss\'{e}
construction and ultrahomogeneous countable structures.  Our main reference
here is Hodges \cite{H}, Ch. 7.  See also \cite{Ch} and \cite{Ca}.

A (countable) {\it signature} is a countable collection $L=\{R_i\}_{i\in I}
\cup\{f_j\}_{j\in J}$ of (distinct) {\it relation} and {\it function symbols}
each of which has an associated number, called its {\it arity}.  The arity
$n(i)$ of each relation symbol $R_i$ is a positive integer and the arity
$m(j)$ of each function symbol $f_j$ is a non-negative integer.  A structure
for $L$ is an object of the form
\[
\bfa=\langle A,\{R^\bfa_i\}_{i\in I},\{f^\bfa_j\}_{j\in J}\rangle ,
\]
where $A$ is a {\it non-empty} set, called the {\it universe} of $\bfa$,
$R^\bfa_i\subseteq A^{n(i)}$, i.e., $R^\bfa_i$ is a $n(i)$-ary
relation on $A$, and $f_j:A^{m(j)}\rightarrow A$, i.e., $f_j^{\bfa}$ is an
$m_j$-ary function on $A$.  (When $m(j)=0,\ f_j^{\bfa}$ is a distinguished
element of $A$.)

Given two structures $\bfa ,\bfb$ of the same signature $L$, a {\it
homomorphism} of $\bfa$ to $\bfb$ is a map $\pi :A\rightarrow B$
such that
\[
R^\bfa_i(a_1,\dots ,a_{n(i)})\Leftrightarrow R^\bfb_i(\pi (a_1),
\dots ,\pi (a_{n(i)}))
\]
and
\[
\pi (f^\bfa_j(a_1,\dots ,a_{m(j)}))=f^\bfb_j(\pi (a_1),\dots ,\pi (a_{m(j)}
)).
\]
We write also in this case $\pi :\bfa\rightarrow\bfb$. ({\it Caution.} Sometimes in the definition of homomorphism, one only requires the left-to-right implication concerning $R_i^{\bfa}, R_i^{\bfb}$. We will use here only the stronger version above.) If $\pi$ is also
1-1, it is called a {\it monomorphism} or {\it embedding}.  Finally,
if $\pi$ is 1-1 and onto it is called an {\it isomorphism}.  If there
is an isomorphism from $\bfa$ to $\bfb$, we say that $\bfa ,\bfb$ are
{\it isomorphic}, in symbols $\bfa\cong\bfb$.  An {\it automorphism}
of $\bfa$ is an isomorphism of $\bfa$ to itself.  We denote by Aut$(\bfa )$
the group of automorphisms of $\bfa$.

We will be primarily interested in countable structures $\bfa$ (i.e., the
universe $A$ is countable).  For countable $\bfa$, the group Aut$(\bfa )$,
with the pointwise convergence topology, is Polish, in fact it is a
closed subgroup of $S_A$, the Polish group of permutations of $A$ with the
pointwise convergence topology. Conversely, given a closed subgroup $G
\subseteq S_A$, there is a signature $L$ and a structure $\bfa_G$
with universe $A$, so that Aut$(\bfa_G )=G$.  To see this, let for each
$n\geq 1,\o_1^n,\o_2^n,\dots$ be the orbits of $G$ on $A^n$, the action of
$G$ on $A^n$ being defined by $g\cdot (a_1,\dots ,a_n)=(g(a_1),\dots ,
g(a_n))$.  Let $L=\{R_{n,i}\}_{n\geq 1}$, where each $R_{n,i}$ is an $n$-ary
relation symbol.  Define $\bfa_G$ by letting $R^{\bfa_G}_{n,i}=\o_i^n\subseteq
A^n$.  Then it is easy to check that Aut$(\bfa_G)=G$.  We call $\bfa_G$
the {\it induced structure} associated to $G$ (see Hodges \cite{H}, 4.1.4, 
where this is called the canonical structure for $G$ - we use however the
term canonical for other purposes in this paper).

A substructure $\bfb$ of $\bfa$ has as universe a (non-empty) subset
$B\subseteq A$ closed under each $f^\bfa_j$, and $R^\bfb_i=R^\bfa_i
\cap B^{n(i)},f^\bfb_j=f^\bfa_j|B^{m(j)}$. We write $\bfb\subseteq\bfa$ in this case. For each $X\subseteq A$, there 
is a smallest substructure containing $X$, called the {\it substructure
generated} by $X$.  A substructure is {\it finitely generated} it is
generated by a finite set.  A structure is {\it locally finite} if all its
finitely generated substructures are finite.  For example, if $L$ is
{\it relational}, i.e., $J=\emptyset$, the substructure generated by $X$
has universe $X$ and so every finitely generated substructure is finite.
This is also true if $J$ is finite and each $f_j$ has arity 0.

A structure $\bfa$ is called {\it ultrahomogeneous} if every isomorphism
between finitely generated substructures $\bfb ,\bfc$ of $\bfa$ can be 
extended to an automorphism of $\bfa$.  For example, $\langle\bbQ ,<
\rangle$, the rationals with the usual order, is an ultrahomogeneous
structure.  Fraiss\'e's theory provides a general analysis of
ultrahomogeneous countable structures.

Let $\bfa$ be a structure for $L$.  The {\it age} of $A$, Age$(\bfa )$
is the collection of all finitely generated structures in $L$ that can
be embedded in $\bfa$, i.e., the closure under isomorphism of the 
collection of finitely generated substructures of $\bfa$.  Clearly the class $\k =$ Age$(\a )$ is non-empty, and satisfies the following
two properties:

(i) {\it Hereditary property (HP)}:  If $\bfb\in\k$ and $\bfc$ is a finitely
generated structure that can be embedded in $\bfb$, then $\bfc\in\k$.

(ii) {\it Joint embedding property (JEP)}:  If $\bfb ,\bfc\in\k$, there
is $\bfd\in\k$ such that $\bfb ,\bfc$ can be embedded in $\bfd$.

When $\bfa$ is moreover ultrahomogeneous, it is easy to see that $\k =$
Age$(\bfa )$ satisfies also the following crucial property:

(iii) {\it Amalgamation property (AP)}:  If $\bfb ,\bfc ,\bfd\in\k$
and $f:\bfb\rightarrow\bfc ,g:\bfb\rightarrow\bfd$ are embeddings, then
there is $\bfe\in\k$ and embeddings $r:\bfc\rightarrow\bfe ,s:\bfd
\rightarrow\bfe$ such that $r\circ f=s\circ g$.

We summarize:

\begin{proposition}
\label{2.1}
Let $\bfa$ be an ultrahomogeneous
structure.  Then $\k =\;${\rm Age}$(\bfa )$ is non-empty, 
and satisfies HP, JEP and AP. 
\qed
\end{proposition}

If $\bfa$ is countable, then clearly {\rm Age}$(\bfa )$ contains only
countably many isomorphism types.  Abusing language, we say that a class
$\k$ of structures is {\it countable} if it contains only countably many
isomorphic types.

We now have the following main result of Fra\"{\i}ss\'{e} \cite{Fr}.

\begin{theorem}[Fra\"{\i}ss\'{e}]
\label{2.2}
Let $L$ be a signature and
$\k$ a class of finitely generated structures for $L$, which is non-empty, 
countable, and satisfies HP, JEP and AP.
Then there is a unique, up to isomorphism, countable structure $\bfa$
such that $\bfa$ is ultrahomogeneous and $\k =$ {\rm Age}$(\bfa )$.
\qed
\end{theorem}

We call this structure the {\it Fra\"{\i}ss\'{e} limit} of $\k$,
\[
\bfa ={\rm Flim} (\k ).
\]
Thus a countable ultrahomogeneous structure is the Fra\"{\i}ss\'{e} limit of
its age.  For example, Age$(\langle\bbQ ,<\rangle )=$ the class of
finite linear orderings $=\l O$ and so $\langle\bbQ ,<\rangle =
{\rm Flim} (\l O)$.

We note here the following alternative characterization of ultrahomogeneity.

\begin{proposition}
\label{2.3}
Let $\bfa$ be a countable structure.
Then $\bfa$ is ultrahomogeneous iff it satisfies the following {\it
extension property}:
\par
If $\bfb ,\bfc$ are finitely generated and can be embedded in $\bfa$, $f:\bfb\rightarrow\bfa ,g:\bfb
\rightarrow\bfc$ are embeddings, then there is 
an embedding $h:\bfc\rightarrow\bfa$
such that $h\circ g=f$.
\qed
\end{proposition}

It follows that if $\bfa$ is countable, $\bfb$ is countable ultrahomogeneous
with Age$(\bfa )\subseteq {\rm Age} (\bfb )$ and $\bfc\subseteq\bfa$ is
finitely generated, then every embedding $f:\bfc\rightarrow\bfb$ can be 
extended to an embedding $g:{\bfa}\rightarrow {\bfb}$.

In particular, $\bfb$ is {\it universal} for the class of all countable
structures $\bfa$ whose age is contained in that of $\bfb$, i.e., every
such $\bfa$ can be embedded in $\bfb$.  For example, any countable linear
ordering can be embedded in $\langle\bbQ ,<\rangle$.

We will be primarily interested in the case when the classes $\k$ in \ref{2.2}
actually consist of {\it finite} structures and the Fra\"{\i}ss\'{e}
limit of $\k$ is countably infinite.  It will be convenient then to 
introduce, for later use, the following terminology, where the {\it
cardinality} of a structure $\bfa =\langle A,\dots\rangle$ is the
cardinality of its universe $A$.

\begin{definition}
\label{2.4}
Given a signature $L$, a {\rm Fra\"{\i}ss\'{e}
class} in $L$ is a class of {\it finite} structures in $L$, which 
contains structures of arbitrary
large (finite) cardinality, is countable, and satisfies HP, JEP and AP.  A
{\rm Fra\"{\i}ss\'{e} structure} in $L$ is a countably {\it infinite} structure
which is locally finite and ultrahomogeneous.
\end{definition}

Thus the map $\k\mapsto{\rm Flim} (\k )$ is a bijection between Fra\"{\i}ss\'{e}
classes and Fra\"{\i}ss\'{e} structures (up to isomorphism) with inverse the
map $\bfa\mapsto {\rm Age} (\bfa )$.

We would like to point out here that for $G$ a closed subgroup of $S_A$,
the induced structure $\bfa_G$ is ultrahomogeneous and, since the
associated signature is relational, it is locally finite, so it is a Fra\"{\i}ss\'{e}
structure, provided $A$ is infinite.

Finally, for further reference, we recall the following definition. A
class $\k$ of structures for $L$ satisfies the {\it strong amalgamation
property} (SAP) if for any $\bfa ,\bfb ,\bfc \in\k$ and embeddings $f:
\bfa\rightarrow\bfb ,g:\bfa\rightarrow\bfc$, there is $\bfd\in\k$ and
embeddings $r:\bfb\rightarrow\bfd ,s:\bfc\rightarrow\bfd$ with $r\circ
f=s\circ g$, such that moreover $r(B)\cap s(C)=r(f(A))\; (=s(g(A))$.

Similarly, we say that $\k$ satisfies the {\it strong joint embedding
property} (SJEP) if for any $\bfa ,\bfb\in\k$ there is $\bfc\in\k$
and embeddings $f:\bfa\rightarrow\bfc ,g:\bfb\rightarrow\bfc$ such that 
$f(A)\cap g(B)=\emptyset$.

\begin{remark}
In retrospect, one can say that Fra\"{\i}ss\'e's construction
was anticipated by Urysohn \cite{U}, who considered the special case of
the class of finite metric spaces with rational distances. 
He constructed a countable metric space $\bfu_0$ (see Section 6, {\bf (E)} below)
whose completion $\bfu$, known as the
{\it Urysohn space}, is the unique universal ultrahomogeneous (with
respect to isometries) complete separable metric space.  Note that we
can view metric spaces $(X,d)$ as structures in a countable signature
$L=\{R_q\}_{q\in\bbQ},\ R_q$ binary, identifying $(X,d)$ with $ \bfx
=(X,R^\bfx_q)$, where $(x,y)\in R^\bfx_q\Leftrightarrow d(x,y)<q$.
\end{remark}

\section{Structural Ramsey Theory}

We will now recall some concepts and results from Ramsey theory, for which
we refer the reader to Ne\v{s}et\v{r}il \cite{N95}, Ne\v{s}et\v{r}il-R\"odl 
\cite{NR90} and Graham-Rothschild-Spencer \cite{GRS}.

Let $\bfa ,\bfb$ be structures in a signature $L$.  We write
\[
\bfa\leq\bfb
\]
if $\bfa$ can be embedded into $\bfb$.  If $\bfa\leq\bfb$, we let
\[
\begin{pmatrix}\bfb\\ \bfa\end{pmatrix}=\{\bfa_0:\bfa_0\text{ is a
substructure of }\bfb\text{ isomorphic to }\bfa\}.
\]
For $\bfa\leq\bfb\leq\bfc ,\ k=2,3,\dots$, we write, using the Erd\"{o}s-Rado \cite{ER56} arrow-notation,
\[
\bfc\rightarrow (\bfb )^\bfa_k
\]
if for any coloring $c:\begin{pmatrix}\bfc\\ \bfa\end{pmatrix}
\rightarrow\{1,\dots ,k\}$ with $k$ colors, there is $\bfb_0\in
\begin{pmatrix}\bfc \\ \bfb\end{pmatrix}$ which is {\it homogeneous}, in the
sense that for some $1\leq i\leq k$, and all $\bfa_0\in\begin{pmatrix}
\bfb_0\\ \bfa\end{pmatrix},\ c(\bfa_0)=i$, i.e., $\begin{pmatrix}\bfb_0\\
\bfa\end{pmatrix}$ is monochromatic.

Let now $\k$ be a class of finite structures in a signature $L$.  We
say that $\k$ satisfies the {\it Ramsey property} if $\k$ is hereditary (i.e., satisfies HP) and for any $\bfa
\leq\bfb$ in $\k ,\ k=2,3,\dots$, there is $\bfc\in\k$ with $\bfb\leq\bfc$
such that
\[
\bfc\rightarrow (\bfb )^\bfa_k.
\]
Note that by a simple induction, we can restrict this condition to $k=2$.

Let us now mention some examples of classes with the Ramsey property:

(i) Let $L=\{<\}$ and let $\l\o$ be the class of finite linear orderings.
Then $\l\o$ has the Ramsey property, by the classical Ramsey theorem.

(ii) Let $L=\{<,E\},\ <,E$ binary relation symbols and let $\o\g$ be the
class of all finite ordered graphs $\bfa =\langle A,<^\bfa ,E^\bfa\rangle$
(i.e., $<^\bfa$ is a linear ordering of $A$ and $E^\bfa$ is a symmetric,
irreflexive relation).  Then Ne\v{s}et\v{r}il-R\"odl \cite{NR77}, 
\cite{NR83} showed that $\o\g$ has the Ramsey property.

(iii) Let $F$ be a finite field and let $L=\{+\}\cup\{f_\alpha
\}_{\alpha\in F}$ (all function symbols), where $+$ has arity 2 and each
$f_\alpha$ is unary. Any vector space over $F$ can be viewed as a structure
in this language with $+$ representing addition and $f_\alpha$ scalar
multiplication by $\alpha\in F$.  A substructure of a vector space is
clearly a subspace.  Let $\calv_F$ be the class of all finite-dimensional vector spaces
over $F$.  Clearly $\bfa\leq\bfb\Leftrightarrow\dim\bfa\leq\dim\bfb$.  It 
was shown in Graham-Leeb-Rothschild \cite{GLR}, that $\calv_F$ has the Ramsey
property.

(iv) Let now $L=\{0, 1, - ,\wedge ,\vee\}$ (all function
symbols), where $0, 1$ have arity 0 and $-,\wedge ,\vee$ have
arities 1, 2, 2 resp.  Any Boolean algebra is a structure in this language
(with $-$ representing Boolean complementation).  Substructures are again
subalgebras.  Let $\b\a$ be the class of all finite Boolean algebras.
Then the so-called Dual Ramsey Theorem of Graham-Rothschild \cite{GR} can be
equivalently reformulated by saying that $\b\a$ has the Ramsey property
(see Ne\v{s}et\v{r}il \cite{N95}, 4.13).

Finally, let us point out the following connection between the Ramsey
property and the amalgamation property, discussed in the previous section
(see Ne\v{s}et\v{r}il-R\"odl \cite{NR77}, p. 294, Lemma 1).

Let $\k$ be a class of finite {\it rigid} (i.e., having no non-trivial
automorphisms) structures in a signature $L$, which is hereditary.  If $\k$ has the JEP and the Ramsey
property, then $\k$ has the AP.  To see this, fix $\bfa ,\bfb ,\bfc\in\k$
and embeddings $f:\bfa\rightarrow\bfb ,\ g:\bfa\rightarrow \bfc$.  By the
JEP, find $\bfe\in\k$ in which both $\bfb ,\bfc$ can be embedded.
Then find $\bfd\in\k$ such that $\bfd\rightarrow (\bfe )^\bfa_4$ and consider
the coloring $c:\begin{pmatrix}\bfd\\ \bfa\end{pmatrix}\rightarrow\{x:
x\subseteq\{\bfb ,\bfc\}\}$ defined as follows:  Given $\bfa_0\in\begin{pmatrix}
\bfd \\ \bfa\end{pmatrix},\bfb\in c(\bfa_0)\Leftrightarrow$ there is an
embedding $r:\bfb\rightarrow\bfd$ with $r\circ f(\bfa )=\bfa_0$, and 
similarly for $\bfc$.  Let $\bfe_0\in\begin{pmatrix}\bfd \\ \bfe
\end{pmatrix}$ be a homogeneous set.  Then $c(\bfa_0)=\{\bfb ,\bfc\}$, for
all $\bfa_0\in\begin{pmatrix}\bfe_0\\ \bfa\end{pmatrix}$.  For such $\bfa_0$,
there is $r:\bfb\rightarrow\bfd$ with $f\circ r(\bfa )=\bfa_0$ and $s:
\bfc\rightarrow\bfd$ with $g\circ s(\bfa )=\bfa_0$.  So $r\circ f,\ g\circ s$
are isomorphisms of $\bfa$ with $\bfa_0$.  Since $\bfa ,\bfa_0$ are rigid,
it follows that $r\circ f=s\circ g$, so $\bfd ,r,s$ verify the amalgamation
property for $f:\bfa\rightarrow\bfb ,g:\bfa\rightarrow \bfc$.

In particular, if $\k$ is a non-empty class of rigid finite structures which is countable, contains structures of arbitrarily
large cardinality, and satisfies HP and JEP and the Ramsey property, then
$\k$ is a Fra\"{\i}ss\'{e} class, i.e., the age of a countably infinite 
ultrahomogeneous structure.
\medskip

\section{Characterizing extremely amenable automorphism groups by a
Ramsey property}

We will first reformulate the condition that a $G$-flow has a fixed point
in the following manner.

\begin{lemma}
\label{4.1}
Let $G$ be a topological group and $X$
a $G$-flow.  Then the following are equivalent:
\par
(i) The $G$-flow $X$ has a fixed point.
\par
(ii) For every $n=1,2,\dots$, and continuous $f:X\rightarrow\bbR^n,\ 
\epsilon >0,\ F\subseteq G$ finite, there is $x\in X$, such that $|f(x)
-f(g\cdot x)|\leq\epsilon ,\ \forall g\in F$, where $|\ |$ refers to
Euclidean norm.
\end{lemma}

\begin{proof}
(i) $\Rightarrow$ (ii).  This is obvious.

(ii) $\Rightarrow$ (i).  We use a compactness argument.  For $f:X\rightarrow
\bbR^n$ continuous, $\epsilon >0,\ F\subseteq G$ finite, put
\[
A_{f,\epsilon ,F}=\{x\in X:\forall g\in F(|f(x)-f(g\cdot x)|\leq\epsilon
)\}.
\]

{\bf Claim}.  $\bigcap_{f,\epsilon ,F}A_{f,\epsilon ,F}\neq\emptyset$.

\medskip
Granting this, fix $x\in\bigcap_{f,\epsilon ,F}A_{f,\epsilon ,F}$.  Then 
$x$ is a fixed point, since otherwise there is $g\in G$ with $g\cdot x\neq
x$, so there is a continuous $f:X\rightarrow\bbR$ with $f(x)=0, f(g\cdot
x)=1$, thus $x\not\in A_{f,1,\{g\}}$.
\medskip

{\it Proof of the claim}.  Notice that $A_{f,\epsilon , F}$ is closed, so,
by compactness, it is enough to show that for any finite collection
$(f_j,\epsilon_j,F_j),\ j=1,\dots ,m$, we have $\bigcap^m_{j=1}
A_{f_j,\epsilon_j,F_j}\neq\emptyset$.  Put
\[
\bar F=F_1\cup\dots\cup F_m,\ \bar\epsilon =\min\{\epsilon_1,\dots ,
\epsilon_m\}
\]
\[
\bar f=(f_1,\dots ,f_m):X\rightarrow\bbR^{n_1+\dots +n_m},
\]
where $f_i:X\rightarrow\bbR^{n_i}$.  Then $A_{\bar f,\bar\epsilon ,\bar F}
\subseteq\bigcap^m_{j=1} A_{f_j,\epsilon_j,F_j}$.  But $A_{\bar f,\bar
\epsilon ,\bar F}$ is non-empty by (ii).
\end{proof}

We use this to prove the following preliminary characterization.

\begin{proposition}
\label{4.2}
Let $S_\infty$ be the Polish group
of permutations of $\bbN$ with the pointwise convergence topology.
If $G\leq S_\infty$ is a closed subgroup, then the following are
equivalent:
\par
(i) $G$ is extremely amenable.
\par
(ii) For any open subgroup $V$ of $G$, every coloring $c:G/V\rightarrow
\{1,\dots ,k\}$, of the set of left-cosets $hV$ of $V$, and every finite
$A\subseteq G/V$, there is $g\in G$ and $1\leq i\leq k$, such that
$c(g\cdot a)=i,\ \forall a\in A$, where $G$ acts on $G/V$ in the usual
way $g\cdot hV=ghV$.
\end{proposition}

\begin{proof}
(i) $\Rightarrow$ (ii):  Fix $V,k,c$ as in (ii), and consider
the shift action of $G$ on $Y=\{1,\dots ,k\}^{G/V},\ g\cdot p(x)=p
(g^{-1}\cdot x)$, for $p\in Y,\ x\in G/V$.  This is a $G$-flow and $c\in Y$.
Let $X=\overline{G\cdot c}$.  By (i) find a fixed point $\gamma\in X$.
Since $G$ acts transitively on $G/V$, clearly $\gamma :G/V\rightarrow\{
1,\dots ,k\}$ is a constant function, say $\gamma (a)=i,\ \forall a\in G/V$.
Fix now finite $A\subseteq G/V$.  Since $\gamma\in\overline{G\cdot c}$,
there is $g\in G$ such that $g^{-1}\cdot c|A=\gamma |A$, so $c(g\cdot a)=
\gamma (a)=i,\ \forall a\in A$.

(ii) $\Rightarrow$ (i):  Clearly (ii) is equivalent to the corresponding
statement about the space $V/G$ of right-cosets $Vh$ of $V$ on which $G$
acts as usual by $g\cdot Vh=Vhg^{-1}$, and this is what we will use below.
Using \ref{4.1}, it suffices to show that if $X$ is a $G$-flow and $f:X\rightarrow
\bbR^n$ is continuous, $\epsilon >0,\ F\subseteq G$ is finite, then there
is $x\in X$, with $|f(x)-f(h\cdot x)|\leq\epsilon ,\ \forall h\in F$.

As in the proof of \ref{1.3}, there is an open nbhd of $1_G$, say $V$, such
that $\forall h\in V\forall x\in X |f(x)-f(h\cdot x)|\leq\epsilon /3$.
But, since $G$ is a closed subgroup of $S_\infty$, we can assume that $V$
is an open subgroup of $G$ (see \cite{BK}, 1.5).  Partition now
the compact set $f(X)\subseteq\bbR^n$ into sets $A_1,\dots ,A_k$ of
diameter $\leq\epsilon /3$.

Fix $x_0\in X$ and let
\[
U_i=\{g\in G: f(g\cdot x_0)\in A_i\}.
\]
Put $VU_i=V_i$, so that $V_i$ is a union of right-cosets of $V$ and thus
can be viewed as a subset of $V/G$.  Since $\bigcup^k_{i=1}V_i=G/V$, we
can find $c:G/V\rightarrow\{1,\dots ,k\}$ such that $c^{-1}(\{i\})
\subseteq V_i$.  So by (ii) there is $1\leq i\leq k$ and $g\in G$ with
$(F\cup\{1_G\})g\subseteq V_i=VU_i$.  We will now show that $x=g\cdot
x_0$ works.

Indeed, fix $h\in F$.  Let $v\in V$ be such that $vhg\in U_i$, so that
$f(vhg\cdot x_0)=f(vh\cdot x)\in A_i$.  Since $|f(vh\cdot x)-f(h\cdot x)|
\leq\tfrac{\epsilon}{3}$, it follows that $f(h\cdot x)$ is in the
$\tfrac{\epsilon}{3}$-nbhd of $A_i$.  Since $1_G\in F$, it follows that
$|f(x)-f(h\cdot x)|\leq\epsilon$.
\end{proof}

Clearly in \ref{4.2} we can restrict $V$ to any local basis at $1_G$ consisting
of open subgroups.  In particular, if for each non-empty finite $F\subseteq
\bbN$ we let
\[
G_{(F)}=\{g\in G:\forall i\in F(g(i)=i)\}
\]
be the {\it pointwise stabilizer} of $F$, then as $\{G_{(F)}:\emptyset\neq
F\subseteq\bbN ,\ F$ finite\} is a local basis of $1_G$, we can restrict $V$
in \ref{4.2} to be of the form $G_{(F)}$, and moreover it is enough to consider
only $F$ in any cofinal under inclusion collection of finite subsets of
$\bbN$.
\smallskip

\noindent
{\it Remark.}
By the proof of \ref{4.2}, to test extreme amenability of a closed
subgroup $G\leq S_\infty$ it is enough to find fixed points in compact
invariant subsets of the $G$-flow $\{1,\dots ,k\}^{G/V},\ V\leq G$ open.
\medskip

For the next result we will need some further notation and terminology. 
Each $G\leq S_\infty$ acts on the finite subsets of $\bbN$ in the obvious
way
\[
g\cdot F=\{g(i):i\in F\}.
\]
For each finite $\emptyset\neq F\subseteq\bbN$, we let then
\[
G_F=\{g\in G:g\cdot F=F\}
\]
be the stabilizer of $F$ in this action, i.e., the {\it setwise stabilizer}
of $F$.  Clearly $G_{(F)}\leq G_F$ and $[G_F:G_{(F)}]<\infty$.

The $G$-{\it type} of $\emptyset\neq F\subseteq\bbN ,\ F$ finite, is the
orbit $G\cdot F$ of $F$.  A $G$-{\it type} $\sigma$ is the $G$-type of
some finite nonempty $F,\ \sigma =G\cdot F$.  If $\rho ,\sigma$ are $G$-types,
we write
\[
\begin{array}{ll}
\rho\leq\sigma &\Leftrightarrow\exists F\in\sigma\exists F'\in\rho (F'
\subseteq F).\\
&\Leftrightarrow\forall F\in\sigma\exists F'\in\rho (F'\subseteq F)\\
&\Leftrightarrow\forall F'\in\rho\exists F\in\sigma (F'\subseteq F).
\end{array}
\]

Finally, given a signature $L=\{R_i\}_{i\in I}\cup\{f_j\}_{j\in J}$,
we denote by $X_L$ the space of all structures for $L$ with universe
$\bbN$.  Thus
\[
X_L=\prod_i2^{(\bbN^{n_i})}\times\prod_j\bbN^{(\bbN^{m_j})}
\]
If $L$ is relational, i.e., $J=\emptyset ,\ X_L$ is compact (homeomorphic
to $2^\bbN$).  The group $S_\infty$ acts canonically on $X_L$ as follows:
Given $\bfa =\langle\bbN ,\{R^\bfa_i\},\ \{f^\bfa_j\}\rangle$ we let
$g\cdot\bfa =\bfb =\langle\bbN ,\{R^\bfb_i\},\ \{f^\bfb_j\}\rangle$,
where
\[
R^\bfb_i(a_1,\dots ,a_{n(i)})\Leftrightarrow R^\bfa_i(g^{-1}(a_1),\dots
,g^{-1}(a_{n(i)}))
\]
\[
f^\bfb_j(a_1,\dots ,a_{m(j)})=f^\bfa_j(g^{-1}(a_1),\dots ,g^{-1}(a_{m(i)}
)),
\]
so that $g$ is an isomorphism from $\bfa$ to $\bfb$.  This action, called
the {\it logic action}, is clearly continuous.  In particular, if $L$ is
relational, the action of any $G\leq S_\infty$ on $X_L$ is a $G$-flow.
Consider the language $L=\{<\},\ <$ a binary relation symbol, and denote
by LO the compact $S_\infty$-invariant subset of $X_L$ consisting of
all linear orderings $\bfa =\langle\bbN ,<^\bfa\rangle$ on $\bbN$.  Clearly
for any $G\leq S_\infty$, LO is a subflow of the $G$-flow $X_L$.
We say that $G$ {\it preserves an ordering} if this subflow has a fixed
point, i.e., there is an ordering $\prec$ on $\bbN$ such that for every
$g\in G,\ a\prec b\Leftrightarrow g(a)\prec g(b)$.

We now have:

\begin{proposition}
\label{4.3}
Let $G\leq S_\infty$ be a closed
subgroup.  Then the following are equivalent:
\par
(i) $G$ is extremely amenable.
\par
(ii) (a) For any finite $\emptyset\neq F\subseteq\bbN ,\ G_{(F)}=G_F$ and
(b) For any two $G$-types $\rho ,\sigma$ with $\rho\leq\sigma$, and
every finite coloring $c:\rho\rightarrow\{1,\dots ,k\}$, there is
$1\leq i\leq k$ and $F\in\sigma$ such that $c(F')=i,\ \forall F'\subseteq
F,\ F'\in\rho$.
\par
(iii) (a)' $G$ preserves an ordering and (b) as in (ii) above.
\end{proposition}

\begin{proof}
(i) $\Rightarrow$ (iii):  Consider (a)' first.  Since $G$
is extremely amenable and LO is a $G$-flow, there is a fixed point,
i.e., $G$ preserves an ordering.

We next prove (b).  Fix $\rho\leq\sigma ,\ c:\rho\rightarrow\{1,\dots ,k\}$.
Say $G\cdot F'=\rho$.  Then if $V=G_{(F')}$, we note that $V=G_{(F')}=
G_{F'}$ by (a)' and we can identify $G/V$ with $G\cdot F'=\rho$.  Applying
then (ii) of \ref{4.2}, to $V,c,A=\{F'_0\subseteq F_0:F'_0\in\rho\}$, where
$F_0\in\sigma$, we find $1\leq i\leq k$ and $g\in G$ with $c(g\cdot F'_0)
=i,\ \forall F'_0\in A$.  Let $F=g\cdot F_0\in\sigma$.  If $F'\subseteq
F,\ F'\in\rho$, then $g^{-1}\cdot F'=F'_0\subseteq F_0$, and $F'_0\in\rho$,
so $c(g\cdot F'_0)=c(F')=i$.

(iii) $\Rightarrow$ (ii):  Clearly, (a)' $\Rightarrow$
(a).

(ii) $\Rightarrow$ (i):  We verify (ii) of \ref{4.2} for $V$ of the form
$G_{(F)}=G_F,\ \emptyset\not= F\subseteq\bbN$ finite.  If $V=G_F$, then $G/V$
can be identified with $\rho =G\cdot F$.  So fix $c:\rho\rightarrow\{
1,\dots ,k\}$ and $A\subseteq\rho$ finite.  Put $\bigcup A=F_0,\ \sigma
=G\cdot F_0$.  Clearly $\rho\leq\sigma$, so there is $1\leq i\leq k$ and
$g\in G$ such that for all $F'\subseteq g\cdot F_0$ with $F'\in\rho$,
we have $c(F')=i$.  Thus $c(g\cdot F)=i,\ \forall F\in A$.
\end{proof}

We will now use a compactness argument to put this characterization in a 
final form.  It will be convenient first to introduce the following
notation.

\begin{definition}
\label{4.4}
Let $G\leq S_\infty$.  Let $\rho\leq
\sigma$ be $G$-types.  If $F\in\sigma$, we put
\[
\begin{pmatrix}F\\ \rho\end{pmatrix}=\{F'\subseteq F:F'\in\rho\}.
\]
If $\rho\leq\sigma\leq\tau$ are $G$-types, we put
\[
\tau\rightarrow (\sigma )^\rho_k,
\]
where $k=2,3,\dots$, if for every $F\in\tau$ and coloring $c:\begin{pmatrix}
F\\ \rho\end{pmatrix}\rightarrow\{1,\dots ,k\}$, there is $F_0\in
\begin{pmatrix}F\\ \sigma\end{pmatrix}$, which is homogeneous, i.e., $c$
is monochromatic on $\begin{pmatrix}F_0\\ \rho\end{pmatrix}$: for some
$1\leq i\leq k$, and every $F'\in\begin{pmatrix}F_0\\ \rho\end{pmatrix},
\ c(F')=i$.  (Note that this is equivalent to asserting that this is true 
for {\it some} $F\in\tau$.)

We say that $G$ has the {\rm Ramsey property} if for every $G$-types
$\rho\leq\sigma$ and every $k=2,3,\dots$, there is a $G$-type $\tau\geq
\sigma$ with $\tau\rightarrow (\sigma )^\rho_k$.
\end{definition}

We now have

\begin{theorem}
\label{4.5}
Let $G\leq S_\infty$ be a closed subgroup.
Then the following are equivalent:
\par
(i) $G$ is extremely amenable.
\par
(ii) (a) $G$ preserves an ordering and (b) $G$ has the Ramsey property.
\end{theorem}

\begin{proof}
(i) $\Rightarrow$ (ii).  We have already seen (a).  To prove
(b), first note that, by a simple induction, it is enough to restrict
ourselves in \ref{4.4}. to the case $k=2$.  So assume, towards a contradiction,
that for some $G$-types $\rho\leq\sigma$, there is no $\tau\geq\sigma$
with $\tau\rightarrow (\sigma )^\rho_2$.  Fix $F_0\in\sigma$.  Then for
every finite set $E\supseteq F_0$ there is a coloring $c_E:\begin{pmatrix}
E\\ \rho\end{pmatrix}\rightarrow\{1,2\}$ which does not have a homogeneous
set $F\in\begin{pmatrix}E\\ \sigma\end{pmatrix}$.  Pick an ultrafilter
$U$ on the set $I$ of finite non-empty subsets of $\bbN$ such that for
every finite $F\subseteq\bbN ,\ \{E:F\subseteq E\}\in U$.  Then for each
$D\in\rho ,\ \{E\supseteq D\cup F_0:c_E(D)=1\}\in U$ or $\{E\supseteq
D\cup F_0:c_E(D)=2\}\in U$, so put $c(D)=i$ iff $\{E\supseteq D\cup
F_0:c_E(D)=i\}\in U$.  This gives a coloring $c:\rho\rightarrow\{1,2\}$.
Then by \ref{4.3}, (ii) (b), there is $F\in\sigma$ such that $c$ is monochromatic
on $\begin{pmatrix} F\\ \rho\end{pmatrix}$, say with value $i$.  If $D\in
\begin{pmatrix}F\\ \rho\end{pmatrix},\ A_D=\{E\supseteq F\cup F_0:
c_E(D)=c(D)=i\}\in U$, so pick $E\in\bigcap_{D\in\begin{pmatrix}F\\ \rho
\end{pmatrix}}A_D$.  Then $E\supseteq F_0$ and for each $D\in\begin{pmatrix}
F\\ \rho\end{pmatrix},\ c_E(D)=i$, so $F\in\begin{pmatrix}E\\ \sigma
\end{pmatrix}$ is homogeneous for $c_E$, a contradiction.
\medskip

(ii) $\Rightarrow$ (i):  It is of course enough to verify (b) of \ref{4.3} (ii),
which follows trivially from the assumption that $G$ has the Ramsey property.
\end{proof}

\begin{remark}
\label{4.6}
Let $G\leq S_\infty$.  We call a set $T$ of $G$-types
{\it cofinal} if for every $G$-type $\rho$ there is $\sigma\in T$ with 
$\rho\leq\sigma$.  Then it is not hard to see that Theorem 4.5 still holds
if in the definition \ref{4.4} of $G$ having the Ramsey property, we restrict
the $G$-types to be in any given cofinal set of $G$-types.
\end{remark}

We will finally tie-up extreme amenability of automorphism groups with
the structural Ramsey theory of \S 3.

Let $L$ be a signature with a distinguished binary relation symbol $<$
(and perhaps other symbols).  An {\it order structure} for $L$ is a
structure $\bfa$ of $L$ in which $<^\bfa$ is a linear ordering.  If $\k$
is a class of structures of $L$, we say that $\k$ is an {\it order
class} if every $\bfa\in\k$ is an order structure.

We also recall that up to (topological group) isomorphism the closed
subgroups of $S_\infty$ are exactly the same as the automorphism groups
of countable structures and also the same as the Polish groups which admit
a countable nbhd basis at the identity consisting of open subgroups
(see \cite{BK}, 1.5).  So the next result provides a
characterization of the groups in this last class that are extremely
amenable.

\begin{theorem}
\label{4.7}
Let $G\leq S_\infty$ be a closed subgroup.
Then the following are equivalent:
\par
(i) $G$ is extremely amenable.
\par
(ii) $G={\rm Aut} (\bfa )$, where $\bfa$ is the Fra\"{\i}ss\'{e} limit of a
Fra\"{\i}ss\'{e} order class with the Ramsey property.
\end{theorem}

{\bf Proof}. (i) $\Rightarrow$ (ii):  Let $\bfa_G=\langle\bbN ,\dots
\rangle$ be the induced structure for $G$.  As we pointed out in \S
2 (third paragraph before the last remark) $\bfa_G$ is ultrahomogeneous.
Also, since $G$ is extremely amenable, $G$ preserves a linear order 
$\prec$ on $\bbN$.  Let $L$ be the signature obtained from the signature
of $\bfa_G$ by adding a new binary relation symbol $<$.  Let $\bfa$
be the expansion of the structure $\bfa_G$ in which $<^{\bfa}  \;=\;\prec$.
Clearly we have Aut$(\bfa )=G$, in particular $\bfa$ is still 
ultrahomogeneous.  Note also that the signature of $\bfa_G$ and thus of 
$\bfa$ is relational, so $\bfa$ is locally finite.  Thus $\k=$ Age$(\bfa )$
is a Fra\"{\i}ss\'{e} order class.  Noting now that, by ultrahomogeneity, a 
$G$-type is exactly the collection of all substructures of $\bfa$
isomorphic to a given $\bfa_0\in$ Age$(\bfa )$, we see that the $G$
having the Ramsey property is equivalent to Age$(\bfa )$ having the
Ramsey property, so we are done by \ref{4.5}.

(ii) $\Rightarrow$ (i):  Since $\bfa$ is the Fra\"{\i}ss\'{e} limit of a
Fra\"{\i}ss\'{e} order class it is a locally finite order structure.  This
implies that $G$ preserves an ordering and also that the $G$-types of
finite substructures are cofinal in all the $G$-types.  As noted earlier,
the $G$-type of a finite substructure $\bfa_0$ is the collection of all
substructures of $\bfa$ isomorphic to $\bfa_0$, so, by \ref{4.5} and 4.6, it is clear
that $G$ has the Ramsey property, so it is extremely amenable.\hfill$\dashv$
\medskip

We make explicit the following fact observed in
the preceding proof.

\begin{theorem}
\label{4.8}
Let $\k$ be a Fra\"{\i}ss\'{e} order class
and $\bfa={\rm Flim} (\k )$.  Then the following are equivalent:
\par
(i) {\rm Aut}$(\bfa )$ is extremely amenable.
\par
(ii) $\k$ has the Ramsey property.
\qed
\end{theorem}

\section{Reducts}

Let $L=\{R_i\}_{i\in I}\cup \{f_i\}_{j\in J}$ be a signature and
$\bfa =\langle A,\{R^\bfa_i\},\{f^\bfa_j\}\rangle$ a structure for $L$.
If $L_0=\{R_i\}_{i\in I_0}\cup\{f_j\}_{j\in J_0}$, with $I_0\subseteq
I,\ J_0\subseteq J$, so that $L_0\subseteq L$, we let $\bfa_0=\bfa |L_0=
\langle A, \{R^\bfa_i\}_{i\in I_0},\ \{f^\bfa_j\}_{j\in J_0}\rangle$ be the
{\it reduct} of $\bfa$ to $L_0$.  We also call $\bfa$ an {\it expansion} of
$\bfa_0$.  If $\k$ is a class of structures in $L$, we denote by
\[
\k |L_0=\{\bfa |L_0:A\in\k\},
\]
the class of reducts of elements of $\k$, called also the {\it reduct}
of $\k$ to $L_0$.

We have seen that order classes of structures play a crucial role in
extreme amenability of automorphism groups.  We will now examine what
happens to reducts of such classes, when the ordering is dropped.

Let $L$ be a signature with a distinguished binary relation symbol $<$
and let $L_0=L\setminus\{<\}$.  For $\bfa$ a structure for $L$, we denote
by $\bfa_0$ the reduct of $\bfa$ to $L_0$ and for any class $\k$ of
structures for $L$, we denote by $\k_0$ the reduct of $\k$ to $L_0$.
Conversely if $\bfa_0=\langle A_0,\dots\rangle$ is a structure for $L_0$
and $\prec$ a binary relation on $A_0$, we denote by $\langle\bfa_0,
\prec\rangle =\bfa$ the structure for $L$ whose reduct to $L_0$ is
$\bfa_0$ and $\prec\; =\; <^\bfa$ (thus also $A=A_0$).

Consider a Fra\"{\i}ss\'{e} order class $\k$ in a signature $L\supseteq\{<\}$
with Fra\"{\i}ss\'{e} limit $\bff$. We characterize when $\k_0$ is also a Fra\"{\i}ss\'{e}
class with limit $\bff_0=\bff |L_0$.

\begin{definition}
\label{5.1}
Let $L$ be a signature with $L\supseteq
\{<\}$, and put $L_0=L\setminus\{<\}$.  Let $\k$ be a 
class of structures in $L$ and put $\k_0=\k |L_0$.  We say that $\k$ is {\rm reasonable}
if for every $\bfa_0\in\k_0 ,\ \bfb_0\in\k_0$, embedding $\pi: \bfa_0
\rightarrow\bfb_0$, and linear ordering $\prec$ on $A_0$ with $\bfa =
\langle\bfa_0,\prec\rangle\in\k$, there is a linear ordering $\prec'$
on $B_0$, so that $\bfb =\langle\bfb_0,\prec'\rangle\in\k$ and $\pi :\bfa
\rightarrow\bfb$ is also an embedding (i.e., $a\prec b\Leftrightarrow\pi
(a)\prec'\pi (b)$).
\end{definition}

We now have

\begin{proposition}
\label{5.2}
Let $L\supseteq\{<\}$ be a signature
and $\k$ a Fra\"{\i}ss\'{e} order class in $L$.  Let $L_0=L\setminus\{<\},\ 
\k_0=\k |L_0,\ \bff ={\rm Flim}(\k),\ \bff_0=\bff |L_0$.  Then the
following are equivalent:
\par
(i) $\k_0$ is a Fra\"{\i}ss\'{e} class and $\bff_0={\rm Flim}(\k_0)$.
\par
(ii) $\k$ is reasonable.
\end{proposition}

\begin{proof}
(i) $\Rightarrow$ (ii).  Let $\bfa_0\in\k_0,\ \bfb_0\in\k_0,\ 
\pi :\bfa_0\rightarrow\bfb_0$ an embedding, and fix a linear
ordering $\prec$ of $A_0$ with $\bfa =\langle\bfa_0,\prec\rangle\in\k$. Then
there is an embedding $\varphi :\bfa\rightarrow\bff$, and $\varphi$ is
of course also an embedding $\varphi :\bfa_0\rightarrow\bff_0$.  By the
extension property \ref{2.3} (since $\bff_0$ is ultrahomogeneous), there is
an embedding $\psi :\bfb_0\rightarrow\bff_0$ with $\psi\circ\pi =\varphi$.
Let $\prec'=\psi^{-1}(<^\bff |\psi (B_0))$.  Then $\bfb =\langle\bfb_0,
\prec'\rangle\in\k$, as $\bfb$ is isomorphic to a substructure of $\bff$,
and moreover $\pi :\bfa\rightarrow\bfb$ is also an embedding.

(ii) $\Rightarrow$ (i):  Clearly $\k_0=$ Age$(\bff_0)$, so to verify
that $\k_0$ is a Fra\"{\i}ss\'{e} class, we only need to check that $\k_0$
satisfies the AP.  Fix $\bfa_0,\bfb_0,\bfc_0\in\k_0$ and embeddings
$f:\bfa_0\rightarrow\bfb_0,\ g:\bfa_0\rightarrow\bfc_0$.  Let then $\prec$
be a linear order on $A_0$ with $\bfa =\langle\bfa_0,\prec\rangle\in\k$.
Since $\k$ is reasonable, we can find linear orders $\prec',\prec''$ on
$B_0,C_0$ resp., so that the structure $\bfb =\langle B_0,\prec'\rangle\in\k ,\ \bfc
=\langle\bfc_0,\prec''\rangle\in\k$ and $f:\bfa\rightarrow \bfb, g:\bfa
\rightarrow\bfc$ are still embeddings.  By AP for $\k$ find $\bfd\in\k$
and embeddings $r:\bfb\rightarrow\bfd ,\ s:\bfc\rightarrow\bfd$ with $r
\circ f=s\circ g$.  Let $\bfd_0=\bfd |L_0$.  Clearly $r:\bfb_0\rightarrow
\bfd_0, s:\bfc_0\rightarrow\bfd_0$ and we are done.

Finally, we check that Flim($\k_0)=\bff_0$, for which it is enough to
verify that $\bff_0$ has the extension property \ref{2.3}.  So fix $\bfa_0,
\bfb_0\in\k_0,\pi :\bfa_0\rightarrow\bfb_0$ an embedding, and $\varphi
:\bfa_0\rightarrow\bff_0$ an embedding.  Then let $$\prec \; =\;\varphi^{-1}
(<^\bff |\varphi (\bfa_0)),$$ so that $\bfa =\langle\bfa_0,\prec\rangle
\in\k$ and $\varphi :\bfa\rightarrow\bff$ is an embedding.  Since $\k$
is reasonable, let $\prec'$ be a linear ordering on $B_0$ with $\langle
\bfb_0,\prec'\rangle\in\k$ and $\pi :\bfa\rightarrow\bfb$ still an
embedding.  Since $\bff$ satisfies the extension property, there is an
embedding $\psi :\bfb\rightarrow\bff$ with $\psi\circ\pi =\varphi$
and, since clearly also $\psi :\bfb_0\rightarrow\bff_0$, we are done.
\end{proof}

A common way to construct an order class $\k$ in $L\supseteq\{<\}$ is to
start with a class $\k_0$ in $L_0=L\setminus\{<\}$ and take
\begin{center}
$\k=\k_0*\l\o=\{\langle\bfa_0,\prec\rangle :\bfa_0\in\k_0$ and $\prec$
is a linear ordering on $A_0$\}.
\end{center}
For example, if $\k_0$ is the class of finite graphs, $\k$ is the class
of all finite ordered graphs.  We now have

\begin{proposition}
\label{5.3}
Let $L\supseteq\{<\}$ be a signature
and let $L_0=L\setminus\{<\}$.  Let $\k_0$ be a class of structures in $L_0$ and put
$\k =\k_0*\l \o$.  Then the following are equivalent:
\par
(i) $\k$ satisfies the amalgamation property.
\par
(ii) $\k$ satisfies the strong amalgamation property.
\par
(iii) $\k_0$ satisfies the strong amalgamation property.
\end{proposition}

\begin{proof}
Suppose that $\k$ satisfies AP in order to show that $\k_0$
satisfies the strong amalgamation property.  Fix $\bfa_0,\bfb_0,\bfc_0\in
\k_0$ and embeddings $f:\bfa_0\rightarrow\bfb_0,\ g:\bfa_0\rightarrow
\bfc_0$.  There are clearly linear orderings $\prec$ on $A_0,\prec'$
on $B_0$, and $\prec''$ on $C_0$ such that if $\bfa =\langle\bfa_0,
\prec\rangle ,\ \bfb =\langle\bfb_0,\prec'\rangle ,\ \bfc =\langle\bfc_0,
\prec''\rangle$ (which are all in $\k$), then $f:\bfa\rightarrow\bfb ,\ 
g:\bfa\rightarrow\bfc$ are still embeddings, and $f(A_0)\prec'(B_0\setminus
f(A_0))$ and $(C_0\setminus g(A_0))\prec'' f(A_0)$ (where if $\prec\prec$
is a linear order on a set $X$ and $Y,Z\subseteq X$, then $Y\prec\prec Z
\Leftrightarrow\forall y\in Y\forall z\in Z(y\prec\prec z))$.  By the AP
for $\k$, let $r:\bfb\rightarrow\bfd ,\ s:\bfc\rightarrow\bfd ,\bfd\in\k$,
be such that $r\circ f=s\circ g$.  If, towards a contradiction, there is
$d\in r(B)\cap s(C),\ d\not\in r(f(A))$ (where $A=A_0,\ B=B_0,\ C=C_0)$,
then $r(f(A))<^\bfd d$, since $d\in r(B\setminus f(A))$ and $d<^\bfd
r(f(A))=s(g(A))$, since $d\in s(C\setminus g(A))$, which is absurd.
So $r(B)\cap s(C)=r(f(A))$.

Now assume that $\k_0$ satisfies strong amalgamation.  We verify that $\k$
satisfies strong amalgamation.  Let $\bfa ,\bfb ,\bfc\in\k$ and let $f:
\bfa\rightarrow\bfb ,\ g:\bfa\rightarrow\bfc$ be embeddings.  Then also
$f:\bfa_0\rightarrow\bfb_0,\ g:\bfb_0\rightarrow\bfc_0$ are embeddings, so,
by strong amalgamation for $\k_0$, find $r:\bfb_0\rightarrow\bfd_0,\ 
s:\bfc_0\rightarrow\bfd_0,\ \bfd_0\in\k_0$, so that $r\circ f=s\circ g$
and $r(B)\cap s(C)=r(f(A))$.  Then $r(B)\setminus r(f(A)),\ s(C)\setminus
s(g(A)),\ r(f(A))\ (=s(g(A))$ are pairwise disjoint, so clearly there is
an order $\prec$ on $D_0$ such that if $\bfd =\langle\bfd_0,\prec\rangle$,
which is in $\k$, then $r:\bfb\rightarrow\bfd ,\ s:\bfc\rightarrow\bfd$.
\end{proof}

The following is quite obvious:

\begin{proposition}
\label{5.4}
Let $L,L_0,\k ,\k_0$ be as in \ref{5.3}.
Then $\k_0$ satisfies the strong joint embedding property iff $\k$
satisfies the strong joint embedding property.
\qed
\end{proposition}

Clearly it is not true that if $\k$ satisfies the joint embedding property,
then $\k_0$ satisfies the strong joint embedding property.  Consider, e.g.,
$L=\{c\}\cup\{<\}$, where $c$ is a $0$-ary function symbol, $\k =$
all finite structures in $L$.

Finally, we note a condition that implies a connection between the Ramsey
property for $\k$ and $\k_0$.

\begin{definition}
\label{5.5}
Let $L\supseteq\{<\}$ be a signature,
$L_0=L\setminus\{<\},\ \k$ a class of structures in $L$ and $\k_0=\k |L_0$.
We say that $\k$ is {\rm order forgetful} if for every $\bfa ,\bfb\in\k$,
letting $\bfa_0=\bfa |L_0,\ \bfb_0=\bfb |L_0$, we have
\[
\bfa\cong\bfb\Leftrightarrow\bfa_0\cong\bfb_0.
\]
(Notice that this does not say that any isomorphism of $\bfa_0$ with $\bfb_0$
is also an isomorphism of $\bfa$ with $\bfb$.)
\end{definition}

An example of an order forgetful class is the class of 
finite-dimensional vector
spaces $\bfv$ over a finite field $F$ with antilexicographical ordering
induced by an ordering of a basis of $\bfv$.  This was considered in Thomas 
\cite{Th}.

We now have the following fact, which is easy to prove:

\begin{proposition}
\label{5.6}
Let $L\supseteq\{<\}$ be a 
signature, $L_0=L\setminus\{<\}$ and let $\k$ be a class of finite
structures in $L$ which is hereditary. Put
$\k_0=\k |L_0$.  If $\k$ is order forgetful, 
then the following are equivalent:
\par
(i) $\k$ satisfies the Ramsey property.
\par
(ii) $\k_0$ satisfies the Ramsey property.
\qed
\end{proposition}

\section{Extremely amenable automorphism groups}

We will now apply the preceding general results to find many new examples
of extremely amenable automorphism groups.  We will use the following
immediate consequence of earlier results.

\begin{theorem}
\label{6.1}
Let $L$ be a signature with $L\supseteq
\{<\}$ and let $\k$ be a Fra\"{\i}ss\'{e} order class in $L$.  Let $\bff =
{\rm Flim} (\k )$ be the Fra\"{\i}ss\'{e} limit of $\k$, so that $\bff$ is an 
order structure.  Then the following are equivalent:
\par
(i) $G={\rm Aut} (\bff )$ is extremely amenable.
\par
(ii) $\k$ has the Ramsey property.
\end{theorem}

\begin{proof}
(ii) $\Rightarrow$ (i) follows from \ref{4.7}, (ii) $\Rightarrow$
(i).  (i) $\Rightarrow$ (ii) is as in the proof of (i) $\Rightarrow$ (ii)
of \ref{4.7}.
\end{proof}

Ramsey theory provides now many examples of $\k$ satisfying the
Ramsey property and we use them to produce new examples of extremely
amenable groups.

\medskip
{\bf (A) Graphs}

Let $L_0=\{E\}$ be the signature with one binary relation symbol $E$.
Let also $L=\{E,<\}$.  A structure $\bfa_0=\langle A_0,E^{\bfa_0})$ is a
{\it graph} if $E^{\bfa_0}$ irreflexive and symmetric.  An {\it ordered
graph} is a structure $\bfa =\langle\bfa_0,<^\bfa\rangle$ for $L$ in which
$\bfa_0$ is a graph and $<^{\bfa}$ a linear ordering.

Lachlan-Woodrow \cite{LW} classified all Fra\"{\i}ss\'{e} classes $\k_0$ of finite
graphs.  They are exactly the following:

(i) $\g\r =$ all finite graphs.

(ii) For $n=3,4,\dots$, ${\mathcal Forb} (K_n)=$ the class of all finite graphs
omitting $K_n$, the complete graph on $n$ vertices (i.e., the class of finite 
graphs that do not contain $K_n$ as a substructure).

(iii) $\e\q =$ the class of finite equivalence relations.

(iv) For $n=1,2,\dots ,\ \e\q_n=$ the class of finite equivalence relations
with at most $n$ classes.

(v) For $n=1,2,\dots ,\ \e\q^*_n=$ the class of finite equivalence
relations, all of whose classes have at most $n$ elements. 

(vi) The complement $\overline{\k_0}$ of one of the classes $\k_0$ listed
in (ii)-(iv) above, where for any graph $\bfa_0=\langle A_0,E^{\bfa_0}\rangle$
its {\it complement} is $\overline{\bfa_0}=\langle A_0,\overline{E^{
\bfa_0}}\rangle$, where $(x,y)\in \overline{E^{\bfa_0}}\Leftrightarrow
x\neq y$ and $(x,y)\not\in E^{\bfa_0}$, and $\overline{\k_0}
=\{\overline{\bfa_0}:\bfa_0\in\k_0\}$.
\medskip

\noindent{\it Remark}.  Strictly speaking
an equivalence relation is not a graph,
because it is reflexive.  So when we think of an equivalence relation $\bfx
=\langle X,R\rangle$ as a graph, we identify it with $\langle X,R\setminus
\{(x,x):x\in X\}\rangle$.

\medskip
Since the automorphism group of the complement of a given graph $\bfa_0$
is the same as the automorphism group of $\bfa_0$, we do not need to
consider the classes of type (vi).  For any one of the classes $\k_0$
of type (i)-(iv), let $\o\k_0=\k_0*\l O=$ the class of finite ordered graphs
$\bfa =\langle\bfa_0,<^\bfa\rangle$ with $\bfa_0\in\k_0$.  Now for $\k_0$
of type (v) with $n\geq 2$, it is easy to check that $\k_0$ does not have
the strong amalgamation property, so $\o\k_0$ is not a Fra\"{\i}ss\'{e} order
class, by \ref{5.3}.  For $\k_0$ of type (iv) with $n\geq 2,\ \o\k_0$ is a
Fra\"{\i}ss\'{e} order class, whose Fra\"{\i}ss\'{e} limit $\bff =\langle F,E^\bff ,
<^\bff\rangle$, consists of an equivalence relation $E^\bff$ on $F$ with
exactly $n$ classes, which are all infinite, and in which $<^\bff$ is an ordering isomorphic
to the rationals such that every equivalence class is dense.  But then it is
easy to check that Aut$(\bff )$ is not extremely amenable, since it acts
(continously) on the finite (discrete) space $X=$ the set of $E^\bff$
classes, without fixed point.  Finally, notice that $\overline{\e\q_1}=\e
\q^*_1$, so we only need to consider $\k_0 =\g\r, {\mathcal Forb}
 (K_n),n=3,4,\dots
,\ \e\q$ and $\e\q_1$.

Each one of the classes $\k_0=\g\r ,{\mathcal Forb}
 (K_n),n=3,4,\dots ,\ \e\q$
and $\e\q_1$, clearly satisfies the strong amalgamation and strong joint
embedding properties, so $\k =\o\k_0$ is a Fra\"{\i}ss\'{e} order class.  $\k$ is
reasonable in each case.  Finally, Ne\v set\v ril and R\"odl \cite{NR77},
\cite{NR83}, Ne\v set\v ril \cite{N89}, have shown that each one of these classes $\k$, satisfies the Ramsey property (the case $\k_0=\e\q_1$ is of course the
classical Ramsey Theorem), except for $\o\e\q$, for which, despite the claim in \cite{N89}, it fails. Thus if $\bff$ is the Fra\"{\i}ss\'{e} limit of $\k$,
then Aut$(\bff )$ is extremely amenable.  We discuss now each case in
some more detail:

(i) $\k_0=\g\r$:  Then $\bfr ={\rm Flim} (\g\r )$ is called the {\it random
graph}.  It is natural to call $\bfo\bfr ={\rm Flim}(\o\g\r )$ the {\it
random ordered graph}.  It is of the form $\bfo\bfr =\langle\bfr ,
<^{\bfo\bfr}\rangle$, where $<^{\bfo\bfr}$ is an appropriate linear order
of the random graph, isomorphic to the rationals.  Thus we have:

\begin{theorem}
\label{6.2}
The automorphism group of the random ordered
graph is extremely amenable.
\qed
\end{theorem}

Of course the automorphism group of the random graph itself is not extremely
amenable, since it does not preserve an ordering.

(ii) $\k_0={\mathcal Forb}
 (K_n)$ is called the $K_n$-{\it free random graph}
and so we call $\bfo\bfr^n={\rm Flim}(\o {\mathcal Forb}
 (K_n))$ the {\it
random $K_n$-free ordered graph}.  It is of the form $\bfo\bfr^n=\langle
\bfr^n_0,<^{\bfo\bfr^n}\rangle$, with $<^{\bfo\bfr^n}$ a linear ordering
of the $K_n$-free random graph, isomorphic to the rationals.  Thus we have:

\begin{theorem}
\label{6.3}
The automorphism group of the random
$K_n$-free ordered graph is extremely amenable.
\qed
\end{theorem}

(iii) $\k_0=\e\q$:  Then $\bff_0={\rm Flim} (\e\q )$ is the equivalence
relation with infinitely many classes each of which is infinite.  So 
$\bff ={\rm Flim}(\o\e\q )\cong\langle\bbQ ,E,<\rangle$, where $<$ is the
usual ordering on $\bbQ$ and $E$ is an equivalence relation on $\bbQ$ with
infinitely many classes each of which is dense.  Here we have:

\begin{theorem}
\label{6.4}
The class $\o\e\q$ does not have the Ramsey property and so the automorphism group of the rationals
with the usual order and an equivalence relation with infinitely many
classes, all of which are dense, is not extremely amenable.
\end{theorem}

\begin{proof}
Let $\bfa = \langle\{a,b\}, =, \prec\rangle$, where $a\prec b$, $\bfb = \langle \{ a,b,c\}, E, \prec'\rangle$, where $a\prec' b\prec' c$,and the $E$-equivalence classes are $\{a,c\}$ and $\{ b \}$. Then there is no $\bfc$
verifying the Ramsey property for $\bfa, \bfb$. To see this, order the equivalence classes of any potential such $\bfc$ acccording to the order of their least elements. Then color a copy $\{a',b'\}$ of $\bfa$ in $\bfc$ red if $a'$ is in a lower class than $b'$ and green otherwise. Then any copy of $\bfb$ in $\bfc$ realizes both colors.

One can also directly see that the automorphism group is not extremely amenable, as it acts without fixed points on the space of linear orderings of the set of equivalence classes.
\end{proof}

(iv) $\k_0=\e\q_1$:  Then $\bff_0={\rm Flim}(\e\q_1)$ is clearly the
complete graph on a countable infinite set and so, up to isomorphism, 
Flim$(\o\e\q_1)\cong\langle\bbQ ,E,<\rangle$, where $E$ is the complete graph
on $\bbQ$.  But the automorphism group of this structure is exactly that
of $\langle\bbQ ,<\rangle$, so we have

\begin{theorem}[Pestov \cite{P98a}]
\label{6.5}
The automorphism group of
the rationals with the usual order is extremely amenable.
\qed
\end{theorem}

In the preceding example the automorphism group of $\bff_0$ is of course
exactly $S_\infty$, which is not extremely amenable as proved in Pestov
\cite{P98a}.  This is clear, as $S_\infty$ cannot preserve an ordering.

Finally, we discuss another order class which can be obtained from $\e\q$.
This will also play a role in \S 8.

Let $L_0=\{E\},\ L=\{E,<\}$, where $E,<$ are binary relation symbols.  Let $\k$
be the class of structures $\bfa =\langle A,E^\bfa ,<^\bfa\rangle$ in $L$
such that $E^\bfa$ is an equivalence relation on $A,\ <^\bfa$ a linear 
order on $A$, and for every $a<^\bfa b<^\bfa c$, if $(a,c)\in E^\bfa$,
then $(a,b)\in E^\bfa$, i.e., every $E^\bfa$ class is convex in $<^\bfa$.
We call a structure in $\k$ a {\it convexly ordered finite equivalence
relation}.  Clearly $\k |L_0=\e\q$.  Now it is easy to check that $\k$
is a reasonable, Fra\"{\i}ss\'{e} order class.  The Fra\"{\i}ss\'{e} limit of $\k$ is
of the form $\bff =\langle\bff_0,<^\bff\rangle$, where $\bff_0$ is the
equivalence relation with infinitely many classes each of which is infinite
and $<^\bff$ is an ordering such that each equivalence class is convex
and isomorphic to the rationals, and the equivalence classes are also
ordered like the rationals.  So $\bff$ up to isomorphism, is the same as
$\langle\bbQ^2,E,<_\ell\rangle$, where $<_\ell$ is the lexicographical
ordering on $\bbQ^2$ and $E$ is the equivalence relation on $\bbQ^2$ given 
by $(r,s,)E(r',s')\Leftrightarrow r=r'$.

\begin{theorem}
\label{6.6}
The automorphism group of $\bbQ^2$ with the
lexicographical ordering and the equivalence relation $(r,s)E(r',s')
\Leftrightarrow r=r'$ is extremely amenable.
\qed
\end{theorem}

\begin{proof}
Put $\bff =\langle\bbQ^2,E,<_\ell\rangle$.  We will show that
Aut$(\bff )$ is extremely amenable.

Let $I_r=\{r\}\times\bbQ$, so that $r<s\Leftrightarrow I_r<_{\ell}
I_s$.  Let now $\pi\in{\rm Aut}(\bff )$.  Let $f_\pi :\bbQ\rightarrow\bbQ$
be defined by $f_\pi (r)=s\Leftrightarrow\pi (I_r)=I_s$.  Then clearly
$f_\pi \in{\rm Aut} (\langle\bbQ ,<\rangle )$.  Also let for each $r\in\bbQ ,\ 
(g_\pi )_r:\bbQ\rightarrow\bbQ$ be defined by $(g_\pi )_r(s)=t\Leftrightarrow
\pi (r,s))=(f_\pi (r),t)$.  Thus again $(g_\pi )_r\in{\rm Aut}(\langle\bbQ ,
<\rangle )$.  Put $\Theta (\pi )=(f_\pi ,g_\pi )$, where $g_\pi\in{\rm Aut}
(\bbQ )^\bbQ ,\ g_\pi =\{(g_\pi )_r\}$.  Consider the semi-direct
product Aut$(\langle\bbQ ,<\rangle )\ltimes{\rm Aut}(\langle\bbQ ,<
\rangle )^\bbQ$, where Aut$(\langle\bbQ ,<\rangle )$ acts on Aut$(\langle
\bbQ ,<\rangle )^\bbQ$ by shift:  $g\cdot x(r)=x(g^{-1}(r))$.  Then it is
easy to check that $\Theta$ is a (topological group) isomorphism of Aut$(\bff )$ with the group which is the semidirect product Aut$(\langle\bbQ ,<\rangle )\ltimes{\rm Aut} (\langle\bbQ
,<\rangle )^\bbQ$, so it is enough to check that the latter group is 
extremely amenable.  This follows from the following standard closure
properties of extreme amenability and \ref{6.5}.
\end{proof}

\begin{lemma}
\label{6.7}
i) Let $G$ be a topological group and
$\h$ an upward-directed, under inclusion, family of extremely amenable
subgroups of $G$ such that $\bigcup\h$ is dense in $G$.  Then $G$ is
extremely amenable.
\par
ii) Let $G$ be a topological group, $N\trianglelefteq G$ a closed normal
subgroup.  If $N,G/N$ are extremely amenable, so is $G$.
\par
iii) The product of extremely amenable groups is extremely amenable.
\end{lemma}

\begin{proof}
i) Let $X$ be a $G$-flow, so also an $H$-flow,
for any $H\in\h$.  Put $X_H=\{x\in X:\forall h\in H(h\cdot x=x)\}$.
Then $X_H$ is compact, non-empty and $\{X_H:H\in\h\}$ has the finite
intersection property.  So $\bigcap_{H\in\h}X_H\neq\emptyset$, and any 
$x\in\bigcap_{H\in\h}X_H$ is fixed by an $g\in\bigcup\h$, so by $G$.

ii) Let $X$ be a $G$-flow, so also a $N$-flow.  Then $X_N=\{x\in X:
\forall g\in N(g\cdot x=x)\}$ is a compact, non-empty subset of $X$.
It follows easily by the normality of $N$ that $X_N$ is $G$-invariant:
Let $x\in X_N,g\in G$.  If $h\in N$, then $h\cdot (g\cdot x)=hg\cdot x=
g\cdot (g^{-1}hg)\cdot x=g\cdot x$ (as $g^{-1}hg\in N)$, so $g\cdot x\in
X_N$.

Define now an action of $G/N$ on $X_N$ as follows:  $gN\cdot x=g\cdot x$
(this is clearly well-defined).  It is easy to check that this is a
continuous action, so there is a fixed point $x\in X_N$ which is clearly
a fixed point of the $G$-flow on $X$.

iii) Suppose that each $G_i,i\in I$, is extremely amenable.  Then, by ii),
the products $\prod_{i\in I_0}G_i,\ I_0\subseteq I$ finite, are extremely
amenable, and identifying $\prod_{i\in I_0}G_i$ with the subgroup $\prod_{i\in 
I}G'_i$ of $\prod_{i\in I}G_i$, where $G'_i=G_i$, if $i\in I_0, G'_i=\{
1_{G_i}\}$, if $i\not\in I_0$, the family $\{\prod_{i\in I_0}G_i:I_0\subseteq
I$ finite\} is upwards-directed under inclusion and its union is dense in
$\prod_{i\in I}G_i$, so $\prod_{i\in I}G_i$ is extremely amenable.
\end{proof}

\begin{corollary}
\label{6.8}
The class of convexly ordered 
finite equivalence relations satisfies the Ramsey property.
\end{corollary}

\begin{proof} By \ref{6.7} and \ref{6.1}.
\end{proof}

This Ramsey result can be also proved directly (and in fact its unordered
version has already been considered in the literature; see Rado \cite{Ra}, 
Graham-Rothschild-Spencer \cite{GRS}, \S 5, Theorem 5),
but it seems interesting to
reverse the roles here and derive it from an extreme amenability result.
\medskip

{\bf (B) Hypergraphs}

Let $L_0=\{R_i\}_{i\in I}$ be a finite relational signature with $R_i$ of
arity $n(i)\geq 2$.  A {\it hypergraph of type} $L_0$ is a structure
$\bfa_0=\langle A_0,\{R^{\bfa_0}_i\}_{i\in I}\rangle$ for $L_0$ for which
each $R^{\bfa_0}_i$ is irreflexive and symmetric, i.e., $(a_1,\dots ,
a_{n(i)})\in R^{\bfa_0}_i\Rightarrow a_1,\dots ,a_{n(i)}$ are distinct, and
for any permutation $\pi$ of $\{1,\dots ,n(i)\},
(a_1,\dots ,a_{n(i)})\in
R^{\bfa_0}_i\Rightarrow (a_{\pi (i)},\dots ,a_{\pi (n(i))})\in R^{\bfa_0}_i$.
Thus, in essence, $R^{\bfa_0}_i\subseteq [A_0]^{n(i)}=$ the set of subsets
of $A_0$ of cardinality $n(i)$.

Let $\h_{L_0}$ be the class of finite hypergraphs of type $L_0$, let $L=L_0
\cup\{<\}$, and let $\o\h_{L_0}=\h_{L_0}*\l O$ be the class of structures in
$L$ which are finite ordered hypergraphs, i.e., of the form $\bfa =\langle
\bfa_0,<^\bfa\rangle$, with $\bfa_0$ a finite hypergraph, and $<^{\bfa}$ a linear ordering. It is easy to
check that $\o\h_{L_0}$ is a reasonable Fra\"{\i}ss\'{e} order class (note that
$\h_{L_0}$ satisfies strong amalgamation and joint embedding).  The 
Fra\"{\i}ss\'{e} limit Flim$(\h_{L_0})=\bfh_{L_0}$ is called the {\it random
hypergraph of type} $L_0$, so we call Flim$(\o\h_{L_0})=\bfo\bfh_{L_0}$
the {\it random ordered hypergraph of type} $L_0$.  We have $\bfo\bfh_{L_0}=
\langle\bfh_{L_0},<^{\bfo\bfh_{L_0}}\rangle$, where $<^{\bfo\bfh_{L_0}}$ is an
appropriate ordering on $H_{L_0}$ isomorphic to the rationals.  Now 
Ne\v set\v ril-R\"odl \cite{NR77}, \cite{NR83} and Abramson-Harrington \cite{AH78}
(see also Ne\v set\v ril \cite{N95}) have shown
that $\o\h_{L_0}$ satisfies the Ramsey property, so we have

\begin{theorem}
\label{6.9}
The automorphism group of the random
ordered hypergraph of type $L_0$ is extremely amenable.
\qed
\end{theorem}

In case $L_0=\{E\}$, $E$ a binary relation, $\bfo\bfh_{L_0}=\bfo\bfr$, the random
ordered graph, so this generalizes \ref{6.2}.  As another special case, consider
$L_0=\{R_i\}_{i\in I}$, where each $n(i)=1$.  Then $\bfo\bfh_{L_0}\cong\langle
\bbQ ,\{A_i\}_{i\in I},<\rangle$, where $<$ is the usual ordering on $\bbQ$
and each $A_i\subseteq\bbQ$ is dense and co-dense in $\bbQ$.  Thus the 
automorphism group of the rationals with the usual ordering and a finite
family of subsets each of which is dense and co-dense is extremely amenable.

The preceding example can be further generalized.

\begin{definition}
\label{6.10}
Let $L_0=\{R_i\}_{i\in I}$ be a finite
relational signature with $n(i)\geq 2$.  
A hypergraph $\bfa_0=\langle A_0,\{R^{\bfa_0}_i\}_{i
\in I}\rangle$ is called {\rm irreducible} if $A_0$ has at least two
elements and for every $a\neq b$ in $A_0$ there is $i\in I$ 
and $c_1,\dots ,c_{n(i)-2}\in A_0$ such that $(a,b,c_1,\dots ,c_{n(i)-2})\in
R^{\bfa_0}_i$.
\end{definition}

Let $\a$ be a class of finite irreducible hypergraphs of type $L_0$.  Then
${\mathcal Forb} (\a )$ is the class of all finite hypergraphs of type $L_0$
which omit $\a$, i.e., do not contain a substructure isomorphic to a member
of $\a$.

For example, for $L_0=\{E\},\ E$ binary, $\a =\{K_n\},\ n=3,4,\dots ,
{\mathcal Forb} (\a )={\mathcal Forb}(K_n)$, the class of finite graphs
that do not contain the complete graph of $n$ elements as a substructure.

For $\a$ a class of finite irreducible hypergraphs of type $L_0$, we denote by
$\o\forb (\a )$ the class of finite ordered hypergraphs of
type $L_0$ that omit $\a$.  It is again easy to see that $\o\forb
(\a )$ is a reasonable Fra\"{\i}ss\'{e} class.  We can call again Flim${(\mathcal 
Forb} (\a ))$ the {\it random $\a$-free ordered hypergraph of type} $L_0$.
Ne\v set\v ril-R\"odl \cite{NR77}, \cite{NR83} (see also Ne\v set\v ril 
\cite{N95}) proved that $\o\forb (\a )$ satisfies the Ramsey property.  
So we have

\begin{theorem}
\label{6.11}
For each class $\a$ of finite irreducible
hypergraphs of type $L_0$, the automorphism group of the random $\a$-free
ordered hypergraph of type $L_0$ is extremely amenable.
\qed
\end{theorem}

{\bf (C) Vector spaces}

We will now consider an example of a different type.  Fix a finite field $F$
and consider the signature $L_0=\{+\}\cup\{f_\alpha\}_{\alpha\in F}$ with $+$ a
binary function symbol and $f_\alpha$ a unary function symbol.  Vector spaces
over $F$ can be viewed as structures in this signature.  Let $\calv_F$ be the
class of finite vector spaces over $F$.  This is a Fra\"{\i}ss\'{e} class.  Let
$L=L_0\cup\{<\}$, and consider the following order class defined in Thomas
\cite{Th}.  Fix an ordering on $F$ such that the 0 of the field $F$ is the
least element in that ordering.  If $\bfv_0$ is a finite-dimensional vector
space over $F$ of dimension $n$ and $B$ is a basis for $\bfv_0$, then
every ordering $b_1<\dots <b_n$ of $B$ gives an ordering on $V_0$ by
$$\alpha_1b_1+\dots +\alpha_nb_n<_{a\ell}\beta_1b_1+\dots +\beta_nb_n
\Leftrightarrow$$ $$(\alpha_n<\beta_n)\; {\rm or}\; (\alpha_n=\beta_n \;{\rm and}\;\alpha_{n-1}
<\beta_{n-1})\; {\rm or} \;\dots $$ i.e., $<_{a\ell}$ is the antilexicographical 
ordering induced by the ordering of $B$.  A {\it natural ordering} of $V_0$
is one induced this way by an ordering of a basis.  Let $\o\calv_F$ be the
order class of all $\bfv =\langle\bfv_0,<^\bfv\rangle$, such that $\bfv_0$
is finite-dimensional vector space and $<^\bfv$ a natural ordering on $V_0$.  Thomas
\cite{Th} shows that this is a Fra\"{\i}ss\'{e} class.  Next it is easy to check that
$\o\calv_F$ is reasonable.  Now the Fra\"{\i}ss\'{e} limit $\bfv_F$ of $\calv_F$ is 
easily seen to be the vector space over $F$ of countably infinite dimension,
so if $\bfo\bfv_F$ is the Fra\"{\i}ss\'{e} limit of $\o\calv_F$, then $\bfo\bfv_F=
\langle\bfv_F,<^{\bfo\bfv_F}\rangle$, where $<^{\bfo\bfv_F}$ is an
appropriate linear order on $V_F$.  Let us call $\bfo\bfv_F$ the 
$\aleph_0$-{\it dimensional vector space over $F$ with the canonical 
ordering}.

Finally $\calv_F$ has the Ramsey property as shown in Graham-Leeb-Rothschild
\cite{GLR}.  It is easy to see though that $\o\calv_F$ is order forgetful, 
according to Definition \ref{5.5}.  Thus, by \ref{5.6}, $\bfo\bfv_F$ has the Ramsey
property too.

Thus we have:

\begin{theorem}
\label{6.12}
The automorphism group of the
$\aleph_0$-dimensional vector space over a finite field with the canonical
ordering is extremely amenable.
\qed
\end{theorem}

Of course the automorphism group of this vector space is not extremely
amenable, as it cannot preserve an ordering.
\medskip

{\bf (D) Boolean algebras}

Let now $L_0=\{0, 1, -,\wedge ,\vee\}$, where $0,1$
have arity 0, $-$ has arity 1 and $\wedge,\vee $ have arity 2.  
Boolean algebras
are structures in $L_0$.  Let $\b\a$ be the class of finite Boolean algebras.
Then it is not hard to check that $\b\a$ is a Fra\"{\i}ss\'{e} class and its
Fra\"{\i}ss\'{e} limit is $\bfb_\infty$, the countable atomless Boolean algebra.

We will next define natural orderings on finite Boolean algebras similar
to example {\bf (C)}.  Let $\bfb_0$ be a finite Boolean algebra and $A$ its
set of atoms.  Then every ordering $a_1<\dots <a_n$ of $A$ gives an
ordering of $B_0$ as follows:  Given $x,y\in B_0$, we can write them 
uniquely as $x=\delta_1a_1\vee\dots\vee\delta_na_n,\ y=\epsilon_1a_1\vee\dots
\vee\epsilon_na_n$, where $\delta_i,\epsilon_i\in\{0,1\}$, and for $\epsilon
\in\{0,1\},\ b\in B_0$,
\[
\epsilon b=\begin{cases} b,\text{ if }\epsilon =1,\\ 0^\bfb ,\text{ if }\epsilon =0.
\end{cases}
\]
(Here and below, we simply write $\vee$ instead of $\vee^{\bfb}$, when the Boolean algebra $\bfb$ is understood.) Then put
\[
x<_{a\ell}y\Leftrightarrow (\delta_n<\epsilon_n)\;\text{ or }\;(\delta_n=
\epsilon_n\; {\rm and}\;\delta_{n-1}<\epsilon_{n-1})\;\text{ or }\;\dots
\]
i.e., $<_{a\ell}$ is the antilexicographical ordering induced by the
ordering $<$ of the atoms.  Again a {\it natural ordering} of $\bfb$
is one induced this way from an ordering of the set of atoms.

Let $\o\b\a$ be the order class of all $\bfb =\langle\bfb_0, <^{\bfb}\rangle$,
such that $\bfb_0$ is a finite Boolean algebra and $<^\bfb$ is a natural
ordering of $B_0$.  We now have:

\begin{proposition}
\label{6.13}
$\o\b\a$ is a reasonable Fra\"{\i}ss\'{e}
order class.
\end{proposition}

\begin{proof}
First we check that $\o\b\a$ is reasonable (see Definition
\ref{5.1}).  Let $\bfb_1,\bfb_2$ be two finite Boolean algebras and let $\pi :
\bfb_1\rightarrow\bfb_2$ be an embedding.  Let $<_{a\ell}$ be a natural
ordering on $\bfb_1$ induced by an ordering $a_1<a_2<\dots <a_n$ of the
atoms of $B_1$.  Let $\{c_1,\dots ,c_k\}$ be the atoms of $\bfb_2$.
Then $\pi (a_1)=\bigvee^{k_1}_{i=1}c_{1i},\dots ,\pi (a_n)=
\bigvee^{k_n}_{i=1}c_{ni}$, where $\{c_{1i}\}^{k_1}_{i=1},\dots , \{
c_{ni}\}^{k_n}_{i=1}$ is a partition of $\{c_1,\dots ,c_k\}$.  Order
then the atoms of $\bfb_2$ as follows:
\[
c_{11}\prec c_{12}\prec \dots\prec c_{1k_1}\prec c_{21}\prec\dots\prec
c_{2k_2}\prec\dots\prec c_{n1}\prec\dots\prec c_{nk_n},
\]
and let $\prec_{a\ell}$ be the induced antilexicographical ordering on 
$\bfb_2$.  Then clearly $\pi :\langle\bfb_1,<_{a\ell}\rangle\rightarrow
\langle\bfb_2,\prec_{a\ell}\rangle$ is still an embedding.

We next check that $\o\b\a$ is hereditary.  To see this, let $\bfb_2$ be a
finite Boolean algebra and $\bfb_1$ a subalgebra.  Let $<_{a\ell}$ be a
natural ordering on $\bfb_2$ induced by an ordering $a_1<\dots <a_n$ of the
atoms of $\bfb_2$.  Let now $b_1<_{a\ell}\dots <_{a\ell}b_k$ be the atoms
of $\bfb_1$.  Write $b_i=c_{i1}\vee\dots\vee c_{ik_i}$, where $c_{i1},\dots
,c_{ik_1}$ are atoms of $\bfb_2$ and $c_{i1}<\dots <c_{ik_i}$.  Then $c_{ik_i}
<c_{jk_j}$ if $i<j$.  From this it easily follows that $<_{a\ell}|B_1=$
the antilexicographical ordering induced by $<_{a\ell}|\{b_1,\dots ,b_k\}$,
so a substructure of an element of $\o\b\a$ is also an element of $\o\b\a$.

Finally, we check that $\o\b\a$ satisfies the amalgamation property (from
which JEP also follows, since the two element Boolean algebra embeds
in any Boolean algebra).

Suppose $\bfb$ is a finite Boolean algebra and $b_1<^\bfb\dots<^\bfb b_k$
is an ordering of the atoms of $\bfb$ with induced antilexicographical
ordering $<^\bfb_{a\ell}$.  Let also $\bfc ,\bfd$ be finite Boolean
algebras with orderings $c_1<^\bfc\dots <^\bfc c_l,\ d_1<^\bfd\dots <^\bfd
d_m$ of their atoms, and corresponding induced antilexicographical orderings $<^\bfc_{a\ell},
<^\bfd_{a\ell}$.  Suppose we have embeddings $$f:\langle\bfb ,<^\bfb_{a\ell}
\rangle\rightarrow\langle\bfc ,<^\bfc_{a\ell}\rangle ,\ g:\langle\bfb ,
<^\bfb_{a\ell}\rangle\rightarrow\langle\bfd ,<^\bfd_{a\ell}\rangle.$$We
will find a Boolean algebra $\bfe$ with $m+\ell -k$ atoms and an ordering
$<^\bfe$ on these atoms, so that, if $<^\bfe_{a\ell}$ is the induced
antilexicographical ordering, then there are embeddings $r:\langle\bfc ,
<^\bfd_{a\ell}\rangle\rightarrow \langle\bfe ,\ <^\bfe_{a\ell}\rangle$, $s:\langle\bfd ,
<^\bfc_{a\ell}\rangle\rightarrow <\bfe ,\ <^\bfe_{a\ell}\rangle$, such
that $r\circ f=s\circ g$.  To specify $r,s$, it is of course enough to define
where the atoms of $\bfc ,\bfd$ go.

Let $f(b_i)=c_{i1}\vee\dots\vee c_{ik_i}$, with $c_{i1}<^\bfc\dots <^\bfc
c_{ik_i}$ atoms in $\bfc$.  Then $c_{ik_i}<c_{jk_j}$, if $i<j$.  Similarly,
$g(b_i)=d_{i1}\vee\dots\vee d_{i\ell_i}$, where $d_{i1}<^\bfd\dots <^\bfd
d_{i\ell_i}$ are atoms in $\bfd$ with $d_{i\ell_i}<d_{j\ell_j}$, if $i<j$.

The Boolean algebra $\bfe$ will have atoms $\{\bar c_{ij}\}_{1\leq i\leq k,\ 
1\leq j\leq k_i},\ \{\bar d_{ij}\}_{1\leq i\leq k,\ 1\leq j\leq \ell_i}$ all
distinct, except that
\[
\bar c_{ik_i}=\bar d_{i\ell_i},\ 1\leq i\leq k.
\]
We will now define the ordering $<^\bfe$ on these atoms and decide where
the atoms of $\bfc ,\bfd$ go.

We first order $\{\bar c_{11},\dots ,\bar c_{1k_1}\}\cup\{\bar d_{11},\dots
,\bar d_{1\ell_1}\}$ as follows:
\[
\bar c_{11}<^\bfe\dots <^\bfe\bar c_{1k_1},\bar d_{11}<^\bfe\dots <^\bfe
\bar d_{1\ell_1}(=\bar c_{1k_1})
\]
and extend $<^\bfe$ on the rest in an arbitrary way.  Using the notation
$(a,b]=\{x:a<x\leq b\},\ (-\infty ,a]=\{x:x\leq a\}$ in an arbitrary
ordering, we now define $$r(c_{11})=\bigvee (-\infty ,\bar c_{11}],\ r(c_{12})
=\bigvee (\bar c_{11},\bar c_{12}],\dots ,r(c_{1k_1})=\bigvee (\bar c_{1,k_1-1},
\bar c_{1k_1}],$$ $$s(d_{11})=\bigvee (-\infty ,\bar d_{11}],\dots ,r(d_{1\ell_1})
=\bigvee (\bar d_{1,\ell_1-1},\bar d_{1,\ell_1}],$$
where, for example $\bigvee
(\bar c_{11},\bar c_{12}]$ means $a_1\vee\dots\vee a_p$, where $a_1,\dots 
,a_p$ are the elements of $\{\bar c_{11},\dots ,\bar c_{1k_1}\}\cup
\{\bar d_{11},\dots ,\bar d_{1\ell_1}\}$ in the interval $(\bar c_{11},
\bar c_{12}]$ according to $<^\bfe$.  Clearly $r,s$ are order preserving
from $\{ c_{11} ,\dots , c_{1k_1}\}, \{ d_{11},\cdots d_{1\ell_1} \}$, resp., to $<^\bfe_{a\ell}$ and $r\circ 
f(b_1)=s\circ g(b_1)=\bar c_{11}\vee\dots\vee\bar c_{1k_1}\vee\bar d_{11}
\vee\dots\vee\bar d_{1\ell_1}$.  

Next we extend $<^\bfe$ to $\{\bar c_{11},
\dots ,\bar c_{1k_1}\}\cup\{\bar c_{21},\dots ,\bar c_{2k_2}\}\cup\{
\bar d_{11},\dots ,\bar d_{1\ell_1}\}\cup\{\bar d_{21},\dots ,\bar d_{2\ell_2}
\}$ and define $r(c_{2i}),s(d_{2i})$.  We simply do that by requiring that
$c_{ij}\mapsto\bar c_{ij}\ (i=1,2,\ 1\leq j\leq k_i)$ is order preserving,
and also $d_{ij}\mapsto\bar d_{ij} \;(i=1,2,\ 1\leq j\leq\ell_i)$ is order
preserving, and define it arbitrarily otherwise.  Notice that this guarantees
that $\bar c_{2k_2}=\bar d_{2\ell_2}$ is the largest element, in particular
$\bar c_{2k_2}>\bar c_{1k_1},\bar d_{2\ell_2}>\bar d_{1\ell_1}$.   We then
extend $r,s$ by defining $$r(c_{21})=\bigvee (-\infty ,\bar c_{21}],\dots ,
r(c_{2k_2})=\bigvee (\bar c_{2,k_2-1},\bar c_{2k_2}],$$ $$ s(d_{21})=\bigvee (-\infty ,
\bar d_{21}],\dots ,s(d_{2\ell_2})=\bigvee (\bar d_{2,\ell_2-1},\bar d_{2\ell_2}
],$$where these intervals now refer to the ordering $<^\bfe$ restricted
to $$\{c_{21},\dots ,c_{2k_2}\}\cup\{d_{21},\dots ,d_{2\ell_2}\}.$$Then
$r,s$ are still order preserving and $r\circ f(b_2)=s\circ g(b_2)=\bar c_{21}
\vee\dots\vee\bar c_{2k_2}\vee \bar d_{21}\vee\dots\vee \bar d_{2\ell_2}$.
Proceeding this way, we define $<^\bfe$ on all the atoms of $\bfe$ and $r,s$
on all the atoms of $\bfc ,\bfd$, resp., so that $r,s$ are order preserving
on the atoms and $r\circ f(b)=s\circ g(b)$, for any atom $b$ of $\bfb$. Then
$r,s$ extend uniquely to embeddings from $\langle\bfc ,<^\bfc_{a\ell}
\rangle$ to $\langle\bfe ,<^\bfe_{a\ell}\rangle$ and $\langle\bfd ,
<^\bfd_{a\ell}\rangle$ to $\langle\bfe ,<^\bfe_{a\ell}\rangle$, resp., and
$r\circ f=s\circ g$.
\end{proof}

Finally, it is clear that $\o\b\a$ is order forgetful and, since $\b\a$ satisfies
the Ramsey property by Graham-Rothschild \cite{GR} (the Dual Ramsey Theorem), it
follows that so does $\o\b\a$.  Let $\bfo\bfb_\infty =\langle\bfb_\infty
,<^{\bfo B_\infty}\rangle$ be the Fra\"{\i}ss\'{e} limit of $\o\b\a$, which we
call the {\it countable atomless Boolean algebra with the canonical
ordering}.  Then we have:

\begin{theorem}
\label{6.14}
The automorphism group of the countable
atomless Boolean algebra with the canonical ordering is extremely amenable.
\qed
\end{theorem}

And we conclude this section by providing the following 
characterization of the group Aut$(\langle
\bbQ ,<\rangle )$ in terms of extreme amenability.

\begin{proposition}
\label{6.15}
Let $G\leq S_\infty$ be a closed
subgroup of $S_\infty$ which acts transitively on $[\bbN ]^n=\{A\subseteq
\bbN :{\rm card} (A)=n\} , n=1,2,3\dots$. If $G$ is extremely amenable, then
there is an ordering $\prec$ on $\bbN$ with $\langle\bbN ,\prec\rangle
\cong\langle\bbQ ,<\rangle$ and $G={\rm Aut} (\langle\bbN ,\prec\rangle )$.
\end{proposition}

\begin{proof}
Since $G$ preserves an ordering, it follows from 3.11 of \cite{Ca}
that there is an ordering $\prec$ on $\bbN$ with $\langle\bbN ,\prec
\rangle\cong\langle\bbQ ,<\rangle$, such that $G\leq{\rm Aut}(\langle\bbN ,
\prec\rangle )$.  Since, for each $n$, $G$ acts transitively on increasing
$n$-tuples in $\prec$, $G$ is dense in Aut$(\langle\bbN ,\prec\rangle )$,
so $G={\rm Aut} (\langle\bbN ,\prec\rangle )$.
\end{proof}
\medskip

{\bf (E) Metric spaces}

We can view metric spaces $(X,d)$ as structures for the language $L_0=
\{R_q\}_{q\in\bbQ},\ R_q$ binary, identifying $(X,d)$ with $\bfx =\langle
X,\{R^\bfx_q\}_{q\in\bbQ}\rangle$, where 
$(x,y)\in R^\bfx_q\Leftrightarrow d(x,y)<q$.  
Let $\m_\bbQ$ be the class of finite metric spaces with rational
distances.  Then it is not hard to check that $\m_\bbQ$ is a Fra\"{\i}ss\'{e}
class (see, e.g., \cite{Bo}).  Its Fra\"{\i}ss\'{e} limit is $\bfu_0$, 
originally constructed in Urysohn \cite{U}, and which we will call the {\it
rational Urysohn space}.  Let also $\o\m_\bbQ =\m_\bbQ *\l O$ be the class
of finite ordered metric spaces with rational distances.  Since actually
$\m_\bbQ$ satisfies strong amalgamation and joint embedding, it is easy
to verify that $\o\m_\bbQ$ is a reasonable Fra\"{\i}ss\'{e} order class.  Its
Fra\"{\i}ss\'{e} limit Flim$(\o\m_\bbQ )=\bfo\bfu_0$ will be called the {\it
ordered rational Urysohn space}.   In response to an
inquiry of the authors, Ne\v set\v ril \cite{N03} has announced 
that $\o\m_\bbQ$ satisfies the Ramsey
property. So we have:

\begin{theorem}
\label{6.16}
The automorphism group of the ordered
rational Urysohn space is extremely amenable.
\qed
\end{theorem}

This result has an interesting application, which actually was our
motivation for looking at the case of metric spaces.

Let $\bfu$ be the so-called {\it Urysohn space}, see Urysohn \cite{U}. This is
the unique, up to isometry, complete separable metric space which contains
(up to isometry) all finite metric spaces and is ultrahomogeneous, for
isometries.  Uspenskij \cite{Usp90} showed that Iso$(\bfu )$, with the pointwise
convergence topology, is a universal Polish group, i.e., contains up to
isomorphism any Polish group.  Note that the topology of the group Iso$(\bfu
)$ is that of pointwise convergence on the space $\bfu$ equipped with the
metric topology, not the discrete one, unlike the case of Aut$(\bfo\bfu_0)$.
Pestov \cite{P02}, using quite different techniques
than the ones used in our paper, showed that Iso$(\bfu )$ is extremely
amenable. This result has several applications.  We now use \ref{6.16} to provide
a quite different proof of this theorem.

\begin{theorem}[Pestov \cite{P02}]
\label{6.17}
The group of isometries
Iso$(\bfu )$ of the Urysohn space $\bfu$, with the pointwise convergence
topology, is extremely amenable.
\end{theorem}

\begin{proof}  
We start with the following standard fact.

\begin{lemma}
\label{6.18}
Let $G,H$ be topological groups and $\pi :
G\rightarrow H$ a continuous homomorphism with $\pi (G)$ dense in $H$.
If $G$ is extremely amenable, so is $H$.
\end{lemma}

\begin{proof}
Let $X$ be an $H$-flow.  Denote by $\alpha :H\times X\rightarrow
X$ the action.  Define now $\tilde\alpha :G\times X\rightarrow X$ by $\tilde
\alpha (g,x)=\alpha (\pi (g),x)$.  This turns $X$ into a $G$-flow, so there
is a fixed point $x_0\in X$.  Clearly $x_0$ is a fixed point for the $H$-flow,
since $\pi (G)$ is dense in $H$.
\end{proof}

Now denote by $\langle\bfu_0,\prec\rangle$ the ordered rational Urysohn space
(so that $\bfu_0$ is the rational Urysohn space).  Already Urysohn \cite{U}
showed that the completion of $\bfu_0$ is $\bfu$, so we view $\bfu_0$ as a
dense subspace of $\bfu$.  Thus if $g\in{\rm Iso} (\bfu_0)$, there is a
unique extension $\bar g\in{\rm Iso}(\bfu )$.  Since every $g\in{\rm Aut}
(\langle\bfu_0,\prec\rangle )$ is in particular an isometry of $\bfu_0$,
the map $g\mapsto\bar g$ is 1--1 from Aut$(\langle\bfu_0,\prec
\rangle )$ into Iso$(\bfu )$ and it is easy to check that it is continuous.
It only remains to show that its range is dense in Iso$(\bfu )$ and then
use \ref{6.18} and \ref{6.16}.

\begin{lemma}
\label{6.19}
Let $D\subseteq{\rm Iso}(\bfu )$.  Let $d$
be the metric on $\bfu$.  Then $D$ is dense, if the following holds:
\[
(*)\ \forall\epsilon >0\forall x_1,\dots ,x_n\in U\forall h\in {\rm Iso}
(\bfu )
\]
\[
\exists x'_1,\dots ,x'_n,y'_1,\dots ,y'_n\in U\exists g\in D\
\]
\[
(d(x_i,x'_i)<\epsilon ,\ d(h(x_i),y'_i)<\epsilon ,\ g(x'_i)=y'_i,\ 
i=1,\dots ,n).
\]
\end{lemma}

\begin{proof}
To check that $D$ is dense, fix $\epsilon >0,\ h\in{\rm Iso}
(\bfu ),x_1,\dots ,x_n\in U$, in order to find $g\in D$ with $d(g(x_i),
h(x_i))<\epsilon$.

By $(*)$ find $x'_1,\dots ,x'_n,y'_1,\dots y'_n$ and $g\in D$ for
$\epsilon /2$.  Then
\[
\begin{array}{ll}
d(g(x_i),h(x_i))&\leq  d(g (x_i),g(x'_i))+d(g(x'_i), h(x_i))\\
&=d(x_i,x'_i)+d(y'_i,h(x_i))\\
&<\epsilon
\end{array}
\]
\end{proof}

So to check that $\{\bar g:g\in{\rm Aut} (\langle\bfu_0,\prec\rangle )\}$
is dense in Iso$(\bfu )$, it is enough to show the following.

\begin{lemma}
\label{6.20}
Given $x_1,\dots ,x_n,y_1,\dots ,y_n\in U$
such that $x_i\mapsto y_i,\ i=1,\dots n$, is an isometry, and given $\epsilon
>0$, there are $x'_1,\dots ,x'_n,y'_1,\dots y'_n\in U_0$ so that $x'_i
\mapsto y'_i$ is an order preserving (with respect to $\prec$) isometry and
\[
d(x'_i,x_i)<\epsilon ,\ d(y'_i,y_i)<\epsilon ,\ i=1,\dots ,n.
\]
\end{lemma}

\begin{proof}  
By induction on $n$.

$\boldsymbol{n=1}$: Simply choose $x'_1,y'_1\in U_0$ with $d(x'_1,x_1)<
\epsilon,\ d(y'_1,y_1)<\epsilon$.

$\boldsymbol{n\rightarrow n+1}$:  Suppose $x_1,\dots ,x_n,x_{n+1},y_1,\dots
,y_n,y_{n+1}\in U$ are given so that $x_i\mapsto y_i$ is an isometry.
By induction hypothesis, find $x'_1,\dots x'_n,y'_1,\dots ,y'_n\in U_0$,
so that $x'_i\mapsto y'_i$ is an order preserving isometry and $d(x_i,
x'_i)<\epsilon /2,\ d(y_i,y'_i)<\epsilon /2,\ i=1,\dots , n$.  Let $x^0_{n+1},
y^0_{n+1}\in U_0$ be such that $$d(x^0_{n+1},x_{n+1})\break <\epsilon /2,\ 
d(y^0_{n+1},y_{n+1})<\epsilon /2.$$Put $d(x^0_{n+1},x'_i)=d_i,\ d(y^0_{n+1},
y'_i)=d'_i,\ 1\leq i\leq n$.  We can of course assume that $\epsilon$ is
small enough so that $d_i,d'_i>\epsilon$.

Therefore,
\[
\begin{array}{ll}
|d_i-d(x_{n+1},x_i)|=&|d(x^0_{n+1},x'_i)-d(x_{n+1},x_i)|\\
&\leq d(x^0_{n+1},x_{n+1})+d(x_i,x'_i)<\epsilon
\end{array}
\]
and
\[
|d'_i-d(y_{n+1},y_i)|<\epsilon ,
\]
so
\[
\begin{array}{ll}
|d_i-d'_i|&=|d_i-d(x_{n+1},x_i)+d(x_{n+1},x_i)-d(y_{n+1},y_i)+d(y_{n+1},y_i)
-d'_i|\\
&<2\epsilon .
\end{array}
\]

Put $e_i=\tfrac{d_i+d'_i}{2}$, and consider the ordered metric space
$$\langle\{x'_1,\dots ,x'_n,x^0_{n+1},u\},\break d',\prec'\rangle,$$where $d'
(x'_i,x'_j)=d(x'_i,x'_j),\ d'(x'_i,x^0_{n+1})=d(x'_i,x^0_{n+1}),\ d'
(u,x'_i)=e_i$, and $d'(u,x^0_{n+1})$ is any rational number satisfying
the inequalities
\[
d_i+e_i >2\epsilon >d'(u,x^0_{n+1})\geq |d_i-e_i|,\ i=1,\dots ,n.
\]
Notice here that
\[
d_i+e_i=\frac{3d_i+d'_i}{2}>2\epsilon
\]
and 
\[
|d_i-e_i|=\tfrac{|d_i-d'_i|}{2}<\epsilon,
\]
so such a number exists.  We let $\prec'$ agree with the ordering $\prec$
(of $U_0$) for $x'_1,\dots ,x'_n,\ x^0_{n+1}$ and $x'_i\prec' u,x^0_{n+1}
\prec' u$.  We need of course to verify that $d'$ is indeed a metric:

(i) Since $d'(x^0_{n+1},x'_i)=d_i,\ d'(u,x'_i)=e_i$, we need to check that
\[
|d_i-e_i|\leq d'(u,x^0_{n+1})\leq d_i+e_i,
\]
which is given by the definition of $d'(u,x^0_{n+1})$.

(ii) Let $\alpha_{ij}=d(x'_i,x'_j)$.  We need to verify that
\[
|e_i-e_j|\leq\alpha_{ij}\leq e_i+e_j.
\]
We have
\[
|d_i-d_j|\leq\alpha_{ij}\leq d_i+d_j,
\]
since $d_i=d(x^0_{n+1},x'_i)$.  But also $\alpha_{ij}=d(y'_i,y'_j)$, so we
also have
\[
|d'_i-d'_j|\leq\alpha_{ij}\leq d'_i+d'_j.
\]
Adding and dividing by 2, we get
\[
|e_i-e_j|\leq\alpha_{ij}\leq e_i+e_j.
\]

So by the properties of $\langle\bfu_0,\prec\rangle$, we can find a point
$x'_{n+1}\in U_0$ with $x'_i\prec x'_{n+1},\ i=1,\dots n,\ x^0_{n+1}
\prec x'_{n+1}$, and $d(x'_{n+1},x'_i)=e_i,\ d(x'_{n+1},x^0_{n+1})=d'
(u,x^0_{n+1})<2\epsilon$.   Similarly we can find $y'_{n+1}\in U_0$ with 
$y'_i\prec y'_{n+1},\ i=1,\dots ,n,\ y^0_{n+1}\prec y'_{n+1}$ and
$d(y'_{n+1},y'_i)=e_{i},\ d(y'_{n+1},y^0_{n+1})<2\epsilon$.  Then $x'_i
\mapsto y'_i,\ 1\leq i\leq n+1$, is an order preserving isometry, and
$d(x'_{n+1},x_{n+1})\leq d(x'_{n+1},x^0_{n+1})+d(x^0_{n+1},x_{n+1})<3
\epsilon$ and $d(y'_{n+1},y_{n+1})<3\epsilon$, so the proof is complete.
\end{proof}
\end{proof}

A result similar to \ref{6.16} can be proved for the ordered integer Urysohn
space (where we consider the class of ordered finite metric spaces with
integer distances), since Ne\v set\v ril \cite{N03} has also verified the
corresponding Ramsey property. It is also conceivable that one can push
those ideas to find a new proof of the result of Gromov and Milman \cite{GroM}
that the unitary group of the infinite-dimensional separable Hilbert space
is extremely amenable, as well as a recent strengthening of this result
by Pestov \cite{P02}, who established extreme amenability of the group of
(affine) isometries of the same Hilbert space.

\section{Universal minimal flows and the ordering property}

Consider now a signature $L\supseteq\{<\}$ and put $L_0=L\setminus\{<\}$.
Let $\k$ be a reasonable Fra\"{\i}ss\'{e} order class in $L$ and put $\k_0
=\k |L_0$.  Then by \ref{5.2} we know that $\k_0$ is a Fra\"{\i}ss\'{e} class and if 
$\bff ={\rm Flim}(\k),\ \bff_0={\rm Flim}(\k_0)$, we have $\bff_0=F|L_0$.
Let $<^\bff =\prec_0$.  Put $G_0={\rm Aut}(\bff_0)$ and consider the
logic action of $G_0$ on LO, the compact space of linear orderings
on $F_0=F$ (which of course we can identify, if we want, with $\bbN$).
Let $X_\k$ be the orbit closure of $\prec_0,\ \overline{G_0\cdot\prec_0}
\subseteq$ LO.  We first note the following.

\begin{proposition}
\label{7.1}
A linear ordering $\prec$ is in $X_\k$
iff for every finite substructure $\bfb_0$ of $\bff_0,\ \bfb =\langle
\bfb_0,\prec |B_0\rangle\in\k$.
\end{proposition}

\begin{proof}
Assume $\prec\in X_\k$ and fix a finite substructure $\bfb_0$
of $\bff_0$.  Then as $\prec\in\overline{G_0\cdot\prec_0 }$, there is $g\in G_0
={\rm Aut}(\bff_0)$ such that $\prec |B_0=(g\cdot\prec_0) |B_0$.  So if $g^{-1}
(\bfb_0)=\bfa_0$, a substructure of $\bff_0$, and $\bfa =\langle\bfa_0,
\prec_0|A_0\rangle$, which is in $\k$, we have that $g|A_0:A_0\rightarrow B_0$
is an isomorphism of $\bfa$ with $\bfb = \langle \bfb_0 , \prec | B_0\rangle$, 
so $\bfb\in\k$.

Conversely, assume that for every finite substructure $\bfb_0$ of $\bff_0,\ 
\bfb =\langle\bfb_0,\prec |B_0\rangle\in\k$.  Then there is an embedding
$\pi :\bfb\rightarrow\bff$.  If $\pi (\bfb )=\bfa$, then $\bfa$ is a 
substructure of $\bff$ and $\pi$ is an isomorphism of $\bfb ,\bfa$, and 
thus in particular an isomorphism of $\bfb_0,\bfa_0=\bfa |L_0$.  But
$\bfb_0,\bfa_0$ are finite substructures of $\bff_0$, so, by ultrahomogeneity
of $\bff_0$, there is $g\in {\rm Aut}(\bff_0)=G_0$ extending $\pi^{-1}$, so 
in particular, $\prec |B_0=(g\cdot\prec_0)|B_0$.  Since $\bfb_0$ was 
arbitrary, this shows that $\prec\in X_\k$.
\end{proof}

\begin{definition}
\label{7.2}
We call any linear ordering in $X_\k$ a $\k$-{\rm admissible} ordering.
\end{definition}

Clearly, $X_\k$ is a $G_0$-flow.  We will now derive necessary and sufficient
conditions for $X_\k$ to be a minimal $G_0$-flow.

The following concept plays an important role in the Ramsey theory of graphs
and hypergraphs, see Ne\v set\v ril-R\"odl \cite{NR78}, Ne\v set\v ril 
\cite{N95}, 5.2.
We formulate it here in a general context.

\begin{definition}
\label{7.3}
Let $L\supseteq\{<\}$ be a signature,
$L_0=L\setminus\{<\},\ \k$ a class of structures in $L$ and
let $\k_0=\k |L_0$.  We say that $\k$ satisfies the {\rm ordering property}
if for every $\bfa_0\in\k_0$, there is $\bfb_0\in\k_0$ such that for every
linear ordering $\prec$ on $A_0$ and linear ordering $\prec'$ on $B_0$,
if $\bfa =\langle\bfa_0 ,\prec\rangle\in\k$ and $\bfb =\langle\bfb_0 ,\prec'
\rangle\in\k$, we have $\bfa\leq\bfb$.
\end{definition}

For example, if $\k$ is the class of finite ordered graphs, so that $\k_0=$
the class of finite graphs, then $\k$ satisfies the ordering property by
results of Ne\v set\v ril-R\"odl (see, e.g., Ne\v set\v ril \cite{N95}, 5.2 or
Ne\v set\v ril-R\"odl \cite{NR78}).  On the other hand, it is easy to see that the
class $\k =\o\e\q$ does not have the ordering property:  Let $\bfa_0
=\langle A_0,E\rangle$ be an equivalence relation with two classes $\{a,b\},
\{c\}$ and consider the ordering $\prec$ on $A_0$ given by $a\prec c\prec b$.
Then for any $\bfb_0=\langle B_0,F\rangle$, if $\prec'$ is an ordering on
$B_0$, so that each $F$-class is convex in $\prec'$, clearly $\langle
\bfa_0,\prec\rangle\not\leq\langle\bfb_0,\prec'\rangle$.

Now we have

\begin{theorem}
\label{7.4}
Let $L\supseteq\{<\}$ be a signature, $L_0
=L\setminus\{<\},\ \k$ a reasonable Fra\"{\i}ss\'{e} order class in $L$.  Let
$\k_0=\k |L_0$, and $\bff ={\rm Flim}(\k) ,\ \bff_0={\rm Flim}(\k_0)=\bff
|L_0$.  Let $X_\k$ be the set of linear orderings $\prec$ on $F\;(=F_0)$
which are $\k$-admissible.  Let also $G_0={\rm Aut}(\bff_0)$.  Then the
following are equivalent:
\par
(i) $X_\k$ is a minimal $G_0$-flow.
\par
(ii) $\k$ satisfies the ordering property.
\end{theorem}

\begin{proof}
First we will reformulate (i) in a more explicit form.
Below let $\bff =\langle\bff_0,\prec_0\rangle$.

\medskip\noindent
{\bf Claim}.  Let $\prec$ be a linear ordering on $F_0$.  Then $\prec_0\in
\overline{G_0\cdot\prec}$ iff for every $\bfa\in\k$ there is a finite substructure
$\bfc_0$ of $\bff_0$ such that $\bfc =\langle\bfc_0,\prec |C_0\rangle\cong
\bfa$.
\vskip 5pt
\noindent
{\it Proof of claim}.  Suppose first that $\prec_0\in \overline{G_0\cdot\prec}$, and 
fix $\bfa\in\k$. Find
$\bfd =\langle\bfd_0,\prec_0|D_0\rangle$ a finite substructure of $\bff$
with $\bfd\cong\bfa$.  Then there is $g\in G_0$ such that $(g\cdot\prec )|
D_0=\prec_0|D_0$.  Clearly $\bfc_0=g^{-1}(\bfd_0)$ is a substructure of
$\bff_0$ and $\bfc =\langle\bfc_0,\prec |C_0\rangle\cong\bfd\cong\bfa$.

To prove the converse, it is enough to show that given a finite substructure
$\bfa_0$ of $\bff_0$ we can find $g\in G_0$ such that $(g\cdot\prec )|A_0=
\prec_0|A_0$.  Since $$\langle\bfa_0,\prec_0|A_0\rangle =\bfa\in\k$$
there is a finite substructure $\bfc_0$ of $\bff_0$ such that $\langle\bfc_0,
\prec |\bfc_0\rangle\cong\bfa$, say via the isomorphism $\pi :C_0
\rightarrow A_0$.  In particular, $\pi$ is an isomorphism of $\bfc_0$ with
$\bfa_0$ so, by the ultrahomogeneity of $\bff_0$, there is $g\in G_0$
extending it.  Then clearly $(g\cdot\prec )|A_0=\prec_0|A_0$ and the 
proof is complete.\hfill$\dashv$

\medskip

Thus we see that (i) is equivalent to the following statement:

For every $\prec\in X_\k$ and every $\bfa\in\k$ there is a finite
substructure $\bfc_0$ of $\bff_0$ such that $\bfc =\langle\bfc_0,\prec
|C_0\rangle\cong\bfa$.

We now proceed to prove the equivalence of (i) and (ii).

(ii) $\Rightarrow$ (i):  Fix $\prec\in X_\k ,\ \bfa\in\k$ and let $\bfa_0=
\bfa |L_0$.  By (ii), find $\bfb_0\in\k_0$ as in \ref{7.3}.  We can of course
assume that $\bfb_0$ is a substructure of $\bff_0$.  Then we have, since
$\bfb =\langle\bfb_0,\prec |B_0\rangle\in\k ,\ \bfa\leq\bfb$.  Thus there
is a substructure $\bfc$ of $\bfb$ isomorphic to $\bfa$.  Clearly, if 
$\bfc_0=\bfc|L_0,\ \bfc =\langle \bfc_0,\prec |C_0\rangle\cong \bfa$ and we are done.

(i) $\Rightarrow$ (ii):  Notice first that, in order to verify the ordering
property, it is enough to show that for every $\bfa\in\k$ there is $\bfb_0
\in\k_0$ such that for every linear ordering $\prec'_0$ on $B_0$, if $\bfb=
\langle\bfb_0 ,\prec'_0\rangle\in\k$, then $\bfa\leq\bfb$.  This follows from
the JEP for $\k_0$.

So fix $\bfa\in\k$, and for every finite substructure $\bfc_0$ of $\bff_0$, let
\[
X_{\bfc_0}=\{\prec\in X_\k :\bfa\cong\langle\bfc_0,\prec|C_0\rangle\}.
\]
Then (i) implies that
\[
X_\k=\bigcup_{\bfc_0}X_{\bfc_0},
\]
so, since each $X_{\bfc_0}$ is open, by compactness we have $\bfc^1_0,\dots ,
\bfc^n_0$ with $X_\k=\bigcup^n_{i=1}X_{\bfc^i_0}$, so that $\forall\prec
\in X_\k\exists 1\leq i\leq n (\bfa\cong\langle\bfc^i_0,\prec |C^i_0\rangle
)$.  Let $\bfb_0$ be the (finite) substructure of $\bff_0$ generated by
$\bigcup^n_{i=1}C^i_0$, so that
\[
\forall\prec\in X_\k (\bfa\leq\langle\bfb_0,\prec |B_0\rangle ).
\]
Fix now $\prec'_0$, a linear ordering on $B_0$, such that $\bfb =\langle
\bfb_0 ,\prec'_0\rangle\in\k$.  If we can show that we can extend $\prec'_0$
to a linear ordering $\prec'\in X_\k$, then $\bfa\leq\langle\bfb_0,
\prec' |B_0\rangle =\langle\bfb_0,\prec'_0\rangle =\bfb$, and this verifies
the ordering property.  To find such an extension, note that there
is a finite substructure $\bfd_0$ of $\bff_0$ and an isomorphism $\varphi$
from $\bfb$ to $\bfd =\langle\bfd_0,\prec_0|D_0\rangle$.  In particular,
$\varphi$ is an isomorphism of $\bfb_0$ with $\bfd_0$, so, since $\bff_0$
is ultrahomogeneous, there is $g\in {\rm Aut}(\bff_0)=G_0$ extending $\varphi$.
Then $\prec'=g^{-1}\cdot\prec_0\in X_\k$ and clearly $\prec'$ extends
$\prec'_0$.
\end{proof}

We can finally show that $X_\k$ is the universal minimal flow of $G_0={\rm
Aut}(\bff_0)$, when $\k$ has both the Ramsey and ordering properties.

\begin{theorem}
\label{7.5}
Let $L\supseteq\{<\}$ be a signature, $L_0
=L\setminus\{<\},\ \k$ a reasonable Fra\"{\i}ss\'{e} order
class in $L$, and let $\k_0
=\k |L_0$ and $\bff ={\rm Flim}(\k) ,\bff_0={\rm Flim}(\k_0)=\bff |L_0$.  Let
$G_0={\rm Aut}(\bff_0), G={\rm Aut}(\bff)$, and let $X_\k$ be the 
set of linear orderings of
$F\;(=F_0)$ which are $\k$-admissible. 

(i) If $\k$ has the Ramsey property, the $G_0$-ambit $(X_\k ,\prec_0)$
is the universal $G_0$-ambit with the property that $G$ stabilizes the
distinguished point, i.e., it
can be mapped homomorphically to any $G_0$-ambit $(X,x_0)$
with $G\cdot x_0=\{x_0\}$. Thus any minimal subflow of $X_\k$ is the universal 
minimal flow of $G_0$. In particular, the universal minimal flow of $G_0$ is 
metrizable.

(ii) If $\k$ has the Ramsey and ordering properties, $X_\k$ is the universal 
minimal flow of $G_0$.
\end{theorem}

\begin{proof}
By \ref{4.7}, $G$ is extremely amenable.  Let also $<^\bff =\prec_0$, so that
$\bff =\langle\bff_0,\prec_0\rangle$.  By definition, $g\cdot\prec_0=\prec_0$,
for all $g\in G$.

Let $X$ be a $G_0$-flow and let $x_0\in X$ be such that $g\cdot
x_0=x_0,\ \forall g\in G$.
We will find a homomorphism $\varphi$ of the $G_0$-flow $X_\k$ to the
$G_0$-flow $X$ with $\varphi (\prec_0)=x_0$.  

Let $\Phi$ be the closure of the set
\[
\{(g\cdot\prec_0,g\cdot x_0):g\in G_0\}\subseteq X_\k\times X
\]
in the compact Hausdorff space $X_\k\times X$.  We will show that
$\Phi$ is the graph of a function $\varphi :X_\k\rightarrow X$.  Granting
this, we can easily verify that this $\varphi$ works.  First, since
$(\prec_0,x_0)\in\Phi$, we have that $\varphi (\prec_0)=x_0$.  Next,
since the graph of $\varphi$ is closed, $\varphi$ is continuous.
Finally, it is a $G_0$-map, since if $\varphi (\prec )=x$, then there
is a net $\{g_i\}$ in $G_0$ such that $g_i\cdot\prec_0\rightarrow\prec ,
g_i\cdot x_0\rightarrow x$, so that $gg_i\cdot\prec_0\rightarrow g\cdot
\prec ,gg_i\cdot x_0\rightarrow g\cdot x$ and since $(gg_i\cdot\prec_0,
gg_i\cdot x_0)\in\Phi$, we also have that $(g\cdot\prec ,g\cdot x)\in\Phi$,
i.e., $\varphi (g\cdot\prec )=g\cdot x=g\cdot\varphi (\prec)$.

So it is enough to prove that $\Phi$ is the graph of a function $\varphi :
X_\k\rightarrow X$.  First we notice that for any $\prec\in X_\k$, there
is some $x\in X$ with $(\prec ,x)\in\Phi$.  Indeed, let $\{g_i\}$ be a
net (actually a sequence) such that $g_i\cdot \prec_0\rightarrow\prec$.
Then $\{g_i\cdot x_0\}$ is a net in $X$, so there is a subnet $\{g_{i_j}
\cdot x_0\}$ converging to some $x\in X$.  Then $g_{i_j}\cdot\prec_0
\rightarrow\prec ,g_{i_j}\cdot x_0\rightarrow x$, so $(\prec ,x)\in\Phi$.
Finally, we show that if $(\prec ,x)\in\Phi$ and $(\prec ,x')\in\Phi$,
then $x=x'$. This amounts to proving the following property:

$(*)$ If $\{g_i\},\{h_j\}$ are nets in $G_0$ and $g_i\cdot\prec_0\rightarrow
\prec ,h_j\cdot\prec_0\rightarrow\prec$, and $g_i\cdot\prec_0\rightarrow x,
h_j\cdot\prec_0\rightarrow x'$, then $x=x'$.

Recall that any compact Hausdorff space is regular, so as $(x,x')\not\in
\Delta =\{(y,y):y\in X\}$, there are open nbhds $U, U'$ of $x,x'$, resp.,
and $W$ of $\Delta$ with $(U\times U')\cap W=\emptyset$.  For each $y\in X$,
there is an open nbhd $U_y$ of $y$ with $U_y\times U_y\subseteq W$ and thus,
by regularity again, there is an open nbhd $V_y$ of $y$ with $y\in V_y
\subseteq\overline{V_y}\subseteq U_y$. So, by compactness, we can find
compact sets $K_1,\dots ,K_n$ and open sets $U_1,\dots , U_n$ such that
$X=\bigcup^n_{i=1} K_i$, and $K_i\subseteq U_i, U^2_i\subseteq W, i=1,
\dots ,n$.

Since the action of $G_0$ on $X$ is continuous, the map $g\mapsto \varphi_g$
from $G_0$ to $H(X)$, where $\varphi_g(y)=g\cdot y$, is continuous, where
$H(X)$ has the compact-open topology.  Since $K_i\subseteq U_i$, the set
\[
\bigcap^n_{i=1}\{f\in H(X):f(K_i)\subseteq U_i\}
\]
is an open nbhd of the identity of $H(X)$, so if $d_r$ is a right-invariant
compatible metric for $G_0$, there is $\delta >0$ such that $d_r
(1_{G_0},g)<\delta\Rightarrow g\cdot K_i\subseteq U_i,\forall 1\leq i\leq n
\Rightarrow\forall y\in X(y,g\cdot y)\in W$.  So $d_r(g,h)<\delta
\Rightarrow\forall y(y,gh^{-1}\cdot y)\in W\Rightarrow (g\cdot x_0, h\cdot
x_0)\in W$.

We will now choose a standard right-invariant metric for $G_0$.  Without
loss of generality, we can assume that $F=F_0=\bbN$, so that $G_0$ is a
closed subgroup of $S_\infty$.  The following is then a left-invariant
compatible metric on $S_\infty$ and thus on $G_0:$ For $f\neq g$ in
$S_\infty$,
\begin{center}
$d_\ell (f,g)=2^{-k-1}$, where $k$ is least with $f(k)\neq g(k)$.
\end{center}
Let
\[
d_r(f,g)=d_\ell (f^{-1},g^{-1})
\]
be the corresponding right-invariant compatible metric on $S_\infty$ and 
$G_0$.

Next choose a finite substructure $\bfa_0$ of $\bff_0$ such that for
any $f,g\in{\rm Aut}(\bff_0)=G_0$,
\[
f|A_0=g|A_0\Rightarrow d_\ell(f,g)<\delta .
\]
Since $g_i\cdot x_0\rightarrow x,h_j\cdot x_0\rightarrow x',
g_i\cdot\prec_0\rightarrow\prec ,h_j\cdot\prec_0\rightarrow\prec$, find
$N,M$ large enough so that $g_N\cdot x_0\in U, h_M\cdot x_0\in U'$ and
$$g_N\cdot\prec_0|A_0=h_M\cdot\prec_0|A_0=\prec |A_0.$$
Thus for $a,b\in A_0$,
\[
g^{-1}_N(a)\prec_0g^{-1}_N(b)\Leftrightarrow h^{-1}_M(a)\prec_0 h^{-1}_M
(b).
\]
Let $g^{-1}_N\cdot\bfa_0=\bfb_0,\ h^{-1}_M\cdot\bfa_0=\bfc_0$, so that
$\bfb_0, \bfc_0$ are finite substructures of $\bff_0$.  Let $\bfb =\langle
\bfb_0,\prec_0|B_0\rangle ,\ \bfc =\langle\bfc_0,\prec_0|C_0\rangle$,
so that $\bfb ,\bfc$ are finite substructures of $\bff$.  Thus $\bfb
\cong\bfc$ via the isomorphism $\rho : B\rightarrow C$, given by
\[
\rho (g^{-1}_N(a))=h^{-1}_M(a),\ \ a\in A_0.
\]
Since $\bff$ is ultrahomogeneous, there is $r\in{\rm Aut}(\bff )=G$ extending
$\rho$, i.e., for $a\in A_0, \; r(g^{-1}_N(a))=h^{-1}_M(a)$, so $r\circ
g^{-1}_N|A_0=h^{-1}_M|A_0$, thus $d_\ell(r\circ g^{-1}_N,h^{-1}_M)<\delta$,
$d_r(g_N\circ r^{-1},h_M)<\delta$, and therefore $(g_N\cdot (r^{-1}
\cdot x_0),h_M\cdot x_0)<W$, and since $r^{-1}\cdot x_0=x_0$ (because
$r^{-1}\in G$), we have $(g_N\cdot x_0,h_M\cdot x_0)<W$, a
contradiction.
\end{proof}

The converse of \ref{7.5} (ii) will be proved in \ref{10.8} below.
It follows from the preceding result that if $\k_0$ is a Fra\"\i ss\'e
class in a signature $L_0$ with $\bff_0={\rm Flim}(\k_0), G_0={\rm Aut}
(\bff_0)$, then it is important to understand when there is a reasonable
Fra\"\i ss\'e order class $\k$ in $L=L_0\cup\{<\}$ with $\k_0=\k |L_0$,
which also has the Ramsey property (because then any minimal subflow of
the $G_0$-flow $X_\k$ is the universal minimal flow of $G_0$).  Not every
$\k_0$ admits such a $\k$.  For example, let $G_0$ be an infinite
countable discrete group.  Via its left-regular representation, we can
view $G_0$ as a closed subgroup of $S_\infty$, so let $\bff_0$ be the
induced structure for $G_0$ which is ultrahomogeneous (see Section 2 third
paragraph before the last Remark), and let $\k_0={\rm Age}(\bff_0)$.
Then $\k_0$ does not admit a $\k$ as above with the Ramsey property as
then the universal minimal flow of $G_0$ would be metrizable.  We will
say more about this question in Section 10.

Assuming such a $\k$ with the Ramsey property exists, one can ask further
whether another such $\k'$ exists with both the Ramsey and ordering properties
(in that case $X_{\k'}$ would be the universal minimal flow of $G_0$).
We will see in Section 10 that this is always the case.  Note that by \ref{7.5} any $\k'$
that has both the Ramsey and ordering properties has an important
minimality property among all $\k$ that have the Ramsey property:  The
$G_0$-flow $X_{\k'}$ is (up to isomorphism) a subflow of $X_\k$.  Moreover,
any two such $X_{\k'}$ are isomorphic.  We will exploit further these 
minimality and uniqueness properties of classes that have both the Ramsey and 
ordering properties in Section 9.

\section{Calculating universal minimal flows}

We now apply the results in \S 6, \S 7 to compute the universal minimal
flows of several automorphism groups.

Consider first the classes $\o\g\r$, of finite ordered graphs, $\o\forb
(K_n),\ n=3,4,\dots$ of $K_n$-free finite ordered graphs, $\o\e\q_1$ of
complete finite ordered graphs, $\o\h_{L_0}$ of finite ordered hypergraphs
of type $L_0$, $\o\forb (\a )$ of finite ordered hypergraphs of type
$L_0$ that omit $\a$, where $\a$ is a class of finite irreducible 
hypergraphs of type $L_0$, and $\o\m_{\bbQ}$ of finite ordered metric spaces
with rational distances. 

Each of these classes satisfies the ordering property and this follows
easily from the fact (already used above in \S 6) that each of these
classes satisfies the Ramsey property.  (The case of $\o\e\q_1$ is of 
course trivial.)  This is done via a standard Sierpinski-style of coloring
obtained by comparing two orderings (see, e.g., Ne\v set\v ril 
\cite{N95}, p. 1376). A similar argument will deduce the ordering property
for finite ordered metric spaces with rational distances from the
corresponding Ramsey property.  It should be mentioned, however, that 
typically the ordering property for a given class of structures is a result
considerably easier to prove than the corresponding Ramsey property and
can frequently be proved directly (see Ne\v set\v ril-R\"odl \cite{NR78} and
Ne\v set\v ril \cite{N95},\cite{Nes04a}).

In each one of the above cases, the space of admissible orderings is of course
the space LO of all linear orderings on $\bbN$ (which we identify with
the universe of the Fra\"{\i}ss\'{e} limit of each class).  
Thus, by \ref{7.5}, we have

\begin{theorem}
\label{8.1}
For each one of the groups below, its
universal minimal flow is the space {\rm LO} of linear orderings on $\bbN$,
so in particular it is metrizable:

(i) The automorphism group of the random graph.

(ii) The automorphism group of the random $K_n$-free graph.

(iii) {\rm (Glasner-Weiss \cite{GlW02})}
$S_\infty$, the permutation group of $\bbN$.

(iv) The automorphism group of the random hypergraph of type $L_0$.

(v) The automorphism group of the random $\a$-free hypergraph of type 
$L_0$, where $\a$ is a class of irreducible finite hypergraphs of type $L_0$.

(vi) The isometry group of the rational Urysohn space $\bfu_0$.
\qed
\end{theorem}

Consider now the classes of convexly ordered finite equivalence relations,
naturally ordered finite-dimensional spaces over a finite field $F$, and
naturally ordered finite Boolean algebras.  Each is easily seen to
satisfy the ordering property.  So we have:

\begin{theorem}
\label{8.2}
(i) The automorphism group of the structure
$\langle\bbN,E\rangle$, where $E$ is an equivalence relation on $\bbN$ with 
infinitely many classes, each of which is infinite, has
as universal minimal flow the space of all convex orderings on $\bbN$,
i.e., all orderings on $\bbN$ for which each $E$-class is convex.
\par
(ii) Let $\bfv_F$ be the $\aleph_0$-dimensional vector space over a finite
field $F$.  The universal minimal flow of its automorphism group (i.e.,
${\rm GL}(\bfv_F))$ is the space of all orderings on $V_F$, whose restrictions
to finite-dimensional subspaces are natural.
\par
(iii) Let $\bfb_\infty$ be the countable atomless Boolean algebra.  The
universal minimal flow of its automorphism group is the space of all
orderings on $B_\infty$, whose restrictions to finite subalgebras are
natural.
\qed
\end{theorem}

In particular all these universal minimal flows are metrizable.

The question of whether the universal minimal flow of ${\rm GL}(\bfv_F)$ 
is nontrivial was brought to one of the authors' attention by Pierre de la Harpe.

Note that, by Stone duality, $\bfb_\infty$ can be identified with the
algebra of clopen subsets of $2^\bbN$ and that every $g\in$ Aut$(\bfb_\infty
)$ determines and is uniquely determined by a homeomorphism $\sigma (g)
\in H(2^\bbN)$. In Glasner-Weiss \cite{GlW03} there is
another representation of the universal minimal flow of $H(2^\bbN )$. 
They showed that the space $\Phi (2^\bbN )$ of all maximal chains of
closed subsets of $2^\bbN$, defined by Uspenskij \cite{Usp00}, can serve as the
universal minimal flow of the group $H(2^\bbN )$. The existence of an 
isomorphism between the space $\n (\bfb_\infty )$, of all orderings 
of $\bfb_\infty$ whose restrictions to finite subalgebras are natural, 
and $\Phi (2^\bbN )$ is of course a consequence of the uniqueness of
the universal minimal flow but we exhibit below an explicit one.

\begin{theorem}
\label{8.3}
There exists an (explicit) homeomorphism
$\varphi :\Phi (2^\bbN )\rightarrow\n (\bfb_\infty )$ such that:
$\varphi (\sigma (g)\cdot x)=g\cdot\varphi (x)$, for $x\in\Phi (2^\bbN ),\ 
g\in {\rm Aut} (\bfb_\infty ).$
\end{theorem}

\begin{proof}
Given a maximal chain $\f$ of closed subsets of
$2^\bbN$, for every clopen subset $A$ of $2^\bbN$, let
\[
F_A=\bigcap\{F\in\f :F\cap A\neq\emptyset\}
\]
By the maximality of $\f$, $A \cap F_A$ is a singleton. Note that if $A$ 
is included in $B$, then $F_A$ is included in $F_B$,
though they can also be equal.  Note however that if $A$ and $B$ are
disjoint, then $F_A$ and $F_B$ are different, so for each finite Boolean 
algebra $\bfb$ contained in $\bfb_\infty$, we
have a total ordering of the atoms of $\bfb$ and this induces
the antilexicographical ordering on the Boolean algebra $\bfb$. These 
orderings cohere and produce an ordering $<_\f$
of $\bfb_\infty$ which is in $\n (B_\infty )$.  Let $\varphi (F)=<_\f$.
This defines a homeomorphism $\varphi :\Phi (2^\bbN )\rightarrow\n
(\bfb_\infty )$ having the required property.
\end{proof}

We conclude with another example of a calculation of a universal minimal flow, 
which turns out to be finite in this case.

Let $\langle \bbQ, < , E_n \rangle$ be the structure of the rationals with the 
usual ordering and an equivalence relation $E_n$ with exactly $n$ classes, 
each of which is dense in $\bbQ$. If, as in Section 6, {\bf (A)}, we denote 
by $\e\q_n$ the class of finite equivalence relations with at most $n$ classes 
and let $\o\e\q_n = \e\q_n *\l\o$, then $\langle\bbQ, < , E_n\rangle$ is the 
Fra\"{\i}ss\'{e} limit of $\o\e\q_n$. As we pointed out in Section 6, 
{\bf (A)}, the automorphism group $G_n$ of this structure is not extremely 
amenable, when $n\geq 2$, as it acts on the finite space $\bbQ / E_n$ 
without fixed points. We will calculate below the universal minimal 
flow of $G_n$.

\begin{theorem}
\label{8.4}
For each $n\geq 1$, let $G_n$ be the automorphism group of the structure 
$\langle \bbQ, < , E_n \rangle$, where $<$ is the usual ordering on 
$\bbQ$ and $E_n$ is an equivalence relation with exactly $n$ classes 
each dense in $\bbQ$. Let $A_1, \dots , A_n$ be an enumeration of 
$\bbQ/E_n$ and let $H_n =${\rm Aut}$(\langle \bbQ, < , A_1, \dots , 
A_n \rangle )$, where we view each $A_i$ as a unary relation on $\bbQ$. 
Then $H_n$ is a finite index clopen normal subgroup of $G_n$ and 
(the natural action of $G_n$ on) $G_n/H_n$ is the universal minimal 
flow of $G_n$.
\end{theorem}

\begin{proof}
Put $\bff_n = \langle \bbQ, < , A_1, \dots , A_n \rangle$, so that 
$H_n = \text{Aut} (\bff_n)$. We will show that $H_n$ is extremely 
amenable, from which it is straightforward to see that $G_n/H_n$ is 
the universal minimal flow of $G_n$.

To prove that $H_n$ is extremely amenable, let first $\o\p_n$ denote 
the class of all finite structures of the form 
$\bfa = \langle A, <^{\bfa}, P_1^{\bfa}, \dots , P_n^{\bfa} \rangle $, 
with $<^{\bfa}$ a linear ordering of $A$ and $\{ P_i^{\bfa} \}_{i = 1}^n$ 
a partition of $A$ into disjoint sets. Then it is easy to verify that 
$\o\p_n$ is a Fra\"{\i}ss\'{e} order class and $\bff_n$ is its 
Fra\"{\i}ss\'{e} limit. So, by \ref{6.1}, it is enough to verify 
that $\o\p_n$ has the Ramsey property.

Fix $\bfa\subseteq \bfb$ in $\o\p_n$. We need to find $\bfc\in\o\p_n$ 
such that $\bfb\leq\bfc$ and $\bfc\rightarrow(\bfb)^{\bfa}_2.$ Choose a
partition $\bbN=N_1\cup\cdots\cup N_n$ into $n$ pairwise disjoint
infinite subsets. Let $\bfc_{\bbN}=\langle\bbN,<,N_1,\cdots,N_n\rangle$, 
where $<$ is
the usual ordering of $\bbN.$ As in \ref{4.5}, (i)$\Rightarrow$(ii), it 
suffices to show that
$\bfc_{\bbN}\rightarrow (\bfb)^{\bfa}_2.$ For each $1\leq i\leq n,$ let
$\calu_i$ be a fixed nonprincipal ultrafilter on $N_i.$ We use this
to define another sequence $\vec{\calv}=(\calv_l:l=0,1,\cdots,k-1)$ of
ultrafilters, as follows: let $\calv_l=\calu_i$, where
$i\in\{1,\cdots,n\}$ is determined by $a_l\in P_i^A$ with
$a_0,a_1,\cdots,a_{k-1}$ the $<^{\bfa}-$increasing enumeration of
the universe of $\bfa$. A $\vec{\calv}-${\it tree} is a nonempty subset
of $T\subseteq \bbN^{[\leq k]}$, the set of increasing sequences from $\bbN$ 
of length $\leq k$, closed under the restrictions
$t\rightarrow t| l$, such that, if $|t|$ = length of $t$,  

$$ A_t=\{m\in \bbN: t^{\wedge} m\in T\} \in \calv_{|t|},
\mbox{ for all }t\in T\mbox{ of length } <k.$$
Note that a maximal node (i.e., a node of length $k$) of any
$\vec{\calv}-$tree naturally determines a copy of $\bfa$ inside
$\bfc_\bbN$. Note also that the family of all $\vec{\calv}-$trees forms a
base for an ultrafilter of subsets of 
$\begin{pmatrix}\bfc_\bbN\\ \bfa\end{pmatrix}$. So, given a
coloring $c:\begin{pmatrix}\bfc_\bbN\\ \bfa\end{pmatrix} \rightarrow \{1,2\}$, 
one can find
a $\vec{\calv}-$tree $T$ such that $c$ is constant on the set
$T^{\max}=T\cap\bbN^k$ of maximal nodes of $T$. Assume $\bfb=\langle
\{0,1,\cdots ,m-1\},<,P_1^\bfb,\cdots,P_n^\bfb\rangle $, where $<$ is
the natural ordering of $\{0,1,\cdots,m-1\}$. Recursively on $0\leq
j<m$ we construct a strictly increasing sequence $l_j \;(0\leq j <
m)$ of non-negative integers such that for all $0\leq j<m$:
$$(a) \forall i\in \{ 1,\cdots ,n \} (j\in P_i^\bfb
\Leftrightarrow l_j\in N_i),$$
$$(b) \forall t \in T (\text{range} (t)\subseteq \{l_p: p<j\}\; 
\text{and}\; j\in P_i^{\bfb} \; \text{and}\; a_{|t|} \in P_i^{\bfa} \Rightarrow
l_j \in A_t).$$

\medskip
\noindent Then $\{l_j:0\leq j<m\}$ forms a
substructure $\bfb'$ of $\bfc_\bbN$ isomorphic to $\bfb$ (via $j\mapsto
l_j$) such that every $\bfa'\in \begin{pmatrix}\bfb'\\ \bfa\end{pmatrix}$ 
is determined by
a maximal node of $T$. So $c$ is constant on 
$\begin{pmatrix}\bfb'\\ \bfa\end{pmatrix}$, as required. 
\end{proof}

One can similarly see that the universal minimal flow of the automorphism
group of the rationals with the usual ordering and an equivalene relation with infinitely many classes each of which is dense (see Theorem 6.4) is its action on the space of linear orderings on the set of equivalence classes.

\section{A uniqueness result about the ordering property}

We now return to the context of Theorem \ref{7.5}, in order to exploit another
basic fact of topological dynamics, namely the uniqueness of the universal
minimal flow (see \ref{1.1}).

Recall from Section 5 that if $L_0\subseteq L$ are signatures and $\bfa$
is a structure for $L$, then $\bfa_0=\bfa |L_0$ is its reduct to $L_0$.
In this case we also call $\bfa$ an {\it expansion} of $\bfa_0$ to $L$.
Similarly, if $\k$ is a class of structures in $L$ and $\k_0=\k
|L_0$, then we call $\k$ an {\it expansion} of $\k_0$ to $L$.

Let $\k_0$ be a Fra\"{\i}ss\'e class of structures in a signature
$L_0$.  Theorem \ref{7.5} shows that if $\k_0$ admits a reasonable 
Fra\"{\i}ss\'e order expansion $\k$
in $L=L_0\cup\{<\}$, then $X_\k$
is the univeral minimal flow of $G_0=$Aut(Flim$(\k_0)$), provided 
$\k$ satisfies the Ramsey and ordering properties.  Now, by the 
uniqueness of the universal minimal flow, if $\k'$ is another class with
the same properties, then $X_\k$ and $X_{\k'}$ are isomorphic (as
$G_0$-flows), which might suggest that $\k ,\k'$ are the ``same''
in some sense.  In other words, one concludes that among reasonable
Fra\"{\i}ss\'e order classes that expand $\k_0$, there is ``at most
one'' that satisfies both the Ramsey and ordering properties. As we will see below, there may be quite distinct expansions that satisfy just the Ramsey property, so this illustrates an interesting feature of the
ordering property.

We will first formulate a quite general uniqueness result and then consider
special cases in which it can be strengthened.  We will need to introduce 
some concepts first.

Given a signature $L$, a {\it simple formula} in $L$ is a quantifier-free
formula in the infinitary language $L_{\omega_1\omega}$.  Explicitly,
this means that a simple formula is obtained from the atomic formulas of
$L$ by using only negations, countable conjunctions and disjunctions.
In case we allow only negation and {\it finite} conjunctions and disjunctions,
we call this a {\it first-order simple formula}.  Consider now a 
Fra\"{\i}ss\'e class $\k_0$ in a signature $L_0$ and let $L=L_0\cup
\{<\}$.  If $\k ,\k'$ are reasonable Fra\"{\i}ss\'e order classes
in $L$ which are expansions of $\k$, (i.e., $\k_0=\k |L_0=
\k'|L_0$), then we call $\k ,\k'$ {\it simply bi-definable} if there are
simple formulas $\varphi (x,y),\varphi'(x,y)$ in $L$, each with two
variables, such that for any given $\bfa_0\in\k_0$, $\varphi$ and $\varphi'$
define (uniformly) a bijection between the expansions of $\bfa_0$ in the
signature $L$ that are in $\k$ with those that are in $\k'$.  More precisely,
this means the following:  Consider any $\bfa_0\in\k_0$.  If $\bfa =\langle
\bfa_0,\prec\rangle$ is an expansion of $\bfa_0$ which is in $\k$, let
$\prec'$ be the relation on $A_0$ defined by $\varphi$ over $\bfa$, i.e.,
\[
a\prec'b\Leftrightarrow\bfa\models\varphi [a,b].
\]
Put
\[
\Phi (\bfa )=\langle\bfa_0,\prec'\rangle .
\]
Similarly for any $\bfa'\in\k'$ define $\Phi'(\bfa')$ using $\varphi'$.
Then the above condition means that, for each $\bfa_0, \Phi$ is a bijection
between the expansions of $\bfa_0$ in $\k$ with the expansions of $\bfa_0$
in $\k'$, with inverse $\Phi'$.

The following is easy to verify:

\begin{proposition}
\label{9.1}
Let $\k_0$ be a Fra\"{\i}ss\'e 
class in a signature $L_0$, let 
$L=L_0\cup\{<\}$ and let $\k ,\k'$ be
reasonable Fra\"{\i}ss\'e order classes in $L$ that are expansions of
$\k_0$.  Then, if $\k ,\k'$ are simply bi-definable:

(i) $\k$ satisfies the Ramsey property iff $\k'$ satisfies the Ramsey
property.

(ii) $\k$ satisfies the ordering property iff $\k'$ satisfies the ordering
property.
\qed
\end{proposition}

We can now state the main uniqueness result.

\begin{theorem}
\label{9.2}
Let $\k_0$ be a Fra\"{\i}ss\'e class in
a signature $L_0$, let $L=L_0\cup\{<\}$, and let $\k ,\k'$ be reasonable
Fra\"{\i}ss\'e order classes in $L$ that are expansions of $\k_0$.  If
$\k ,\k'$ satisfy the Ramsey and ordering properties, then $\k ,\k'$
are simply bi-definable.
\end{theorem}

\begin{proof}
Let $\bff_0= $ Flim$(\k_0)$, and 
$\bff =$ Flim$(\k)=\langle\bff_0,\prec_0\rangle,\; \bff' = $ 
Flim$(\k') = \langle\bff_0,\prec'_0\rangle$.
Then, by \ref{7.5}, both $X_\k$ and $X_{\k'}$ are universal minimal flows of
$G_0=$ Aut$(\bff_0)$, so there are isomorphic (as $G_0$-flows), i.e., there
is a homeomorphism $\pi :X_\k\rightarrow X_{\k'}$ such that $\pi (g\cdot
\prec)=g\cdot\pi (\prec )$, for $\prec\in X_\k$.

\medskip
\noindent
{\bf Claim}.  $G_0\cdot\prec_0$ is a dense $G_\delta$ set in $X_\k$ (and
similarly for $G_0\cdot\prec'_0$ in $X_{\k'}$).

\medskip
\noindent
{\it Proof}.  By definition, $G_0\cdot\prec_0$ is dense in $X_\k$.

Assuming, without loss of generality, that $F_0=\bbN$, we will show that $G_0\cdot\prec_0$ is $G_{\delta}$ in $X_\k$. Note that $\prec \in G_0\cdot \prec_0$ iff $\langle \bff_0, \prec \rangle \cong \bff$ iff $\langle \bff_0, \prec \rangle$ has age $\k$ and satisfies the extension property 2.3. It is now easy to verify that these properties
can be expressed in a $G_\delta$ way.\hfill$\dashv$
\medskip

Thus $G_0\cdot\prec_0$ is the unique dense $G_\delta$ orbit of $X_\k$
and similarly $G_0\cdot\prec'_0$ is the unique dense $G_\delta$ orbit of
$X_{\k'}$.  It follows that $\pi (G_0\cdot\prec_0)=G_0\cdot\prec'_0$.  Put
$\pi (\prec_0)=\prec^*_0$.  Then $\bff^*=\langle\bff_0,\prec^*_0\rangle$
is also (isomorphic to) the Fra\"{\i}ss\'e limit of $\k'$ and since 
Aut$(\bff)$ is the stabilizer of $\prec_0$ in the action of $G_0$ on $X_\k$,
and similarly for Aut$(\bff^*)$ and $\prec^*_0$, it follows that Aut$(\bff )$
= Aut$(\bff^*)$.

Thus (at the cost of replacing $\bff'$ by its isomorphic copy $\bff^*$) we 
may as well assume that
\[
\text{Aut}(\bff )=\text{Aut}(\bff').
\]

Consider now the action of $S_\infty$ (and its subgroups) on $\bbN^n, n=1,2,
\dots$:
\[
g\cdot (a_1,\dots ,a_n)=(g(a_1),\dots ,g(a_n)).
\]
\noindent
{\bf Claim}.  If $R\subseteq\bbN^n$ is Aut$(\bff )$-invariant, then $R$ is
definable in $\bff$ by a simple formula $\varphi_R$ in $L$.
\medskip

\noindent
{\it Proof}.  Write $R=\bigcup_{i\in I}R_i$, where $R_i$ are the 
Aut$(\bff )$-orbits on $\bbN^n$ contained in $R$.  (Here $I$ is a countable
index set.)  For each $i\in I$, fix $(a_1,\dots ,a_n)\in R_i$ and let
\[
\varphi_i(x_1,\dots ,x_n)
\]
be the simple formula in $L$ which is the conjunction of all atomic or
negations of atomic formulas $\psi (x_1,\dots ,x_n)$ in $L$ such that
\[
\bff\models\psi [a_1,\dots ,a_n].
\]
Then it is easy to see (using that $\bff$ is ultrahomogeneous) that
\[
(b_1,\dots ,b_n)\in R_i\Leftrightarrow\bff\models\varphi_i[b_1,\dots ,
b_n].
\]

(Note here that if $L$ is a finite relational language, then actually
$\varphi_i$ is first-order simple.)

Take then $\varphi_R$ to be the disjunction of $\varphi_i, i\in I$.
\end{proof}

Now the relation $\prec'_0\subseteq\bbN^2$ is invariant under Aut$(\bff')=$
Aut$(\bff )$, so there is a simple formula in $L, \varphi$, such that
\[
a\prec'_0b\Leftrightarrow\bff\models\varphi[a,b].
\]
Similarly, there is a simple formula $\varphi'$ that defines $\prec_0$ 
in $\bff'$. It is easy now to see that $\varphi ,\varphi'$ witness that $\k ,\k'$
are simply bi-definable.\hfill$\dashv$
\medskip

\noindent
{\it Remark}.  Note that the preceding proof shows that if $L_0$ is finite
and relational, then $\k ,\k'$ are first-order simply bi-definable.

\medskip
In certain special instances, we can  actually assert in \ref{9.2}
that $\k = \k'$.

\begin{proposition}
\label{9.3}
In the context of \ref{9.2}, if $\k ,\k'$ are
comparable under inclusion, i.e., $\k\subseteq\k'$ or $\k'\subseteq\k$, then
$\k =\k'$.  In particular, if $\k =\k_0*\l\o$, then $\k =\k'$.
\end{proposition}

\begin{proof}
Say $\k \subseteq \k'$, but $\k'\setminus\k\not=\emptyset$, so 
that there is $\langle\bfa_0,\prec\rangle\in\k'\setminus\k$. 
Now apply the ordering property to $\bfa_0\in\k'|L_0$ to find 
$\bfb_0\in\k'|L_0 = \k|L_0$ and consider $\prec$ on $A_0$ and 
$\prec'$ on $B_0$, so that $\langle \bfb_0,\prec'\rangle \in \k$. 
Then $\langle\bfa_0,\prec\rangle \leq \langle \bfb_0,\prec'\rangle$, 
so, as $\k$ is hereditary, $\langle\bfa_0,\prec\rangle\in\k$, a 
contradiction.
\end{proof}

\medskip
We see from \ref{9.3} that in all the classes considered in \ref{8.1}, 
i.e., $\k_0=$
(finite) graphs, $K_n$-free graphs, complete graphs, hypergraphs of type
$L_0$ that omit $\a$, and metric spaces with rational distances, $\k =
\k_0*\l\o$ is the unique (reasonable order Fra\"{\i}ss\'e) expansion
that satisfies both the Ramsey and ordering properties.

Note that \ref{9.3} is also a trivial corollary of \ref{9.2}, since 
simple bi-definability implies, in particular, that the cardinality 
of the expansions of any $\bfa_0 \in \k|L_0 = \k' |L_0$ which are in 
$\k$ is the same as that of the expansions which are in $\k'$.

By these simple cardinality considerations, we can also see that \ref{9.2} is not true if one of $\k ,\k'$ fails to satisfy the ordering property.  Consider, for example, the case when $\k_0$ is the class of all $\bfa_0 = \langle A_0, P, Q \rangle$, where $\{ P, Q \}$ is a partition of $A_0$, 
$\k$ is the class of all $\bfa = \langle \bfa_0, < \rangle$, where $<$ is an ordering of $A_0$ with $P < Q $, and $\k'=\k_0*\l\o$.  Then $\k$ satisfies the Ramsey and 
ordering properties, $\k'$ satisfies the Ramsey but not the ordering
property, and $\k ,\k'$ are not simply bi-definable.

In case one of $\k ,\k'$ might not be equal to $\k_0*\l\o$, we do not 
necessarily have the strong uniqueness property of \ref{9.3}. But under certain
conditions on $\k$ we can still strengthen \ref{9.2}.  To motivate what we are 
looking for, let for each $\k$, in the context of \ref{9.2},
\[
\k^*=\{\langle\bfa_0,\prec^*\rangle :\langle A_0,\prec\rangle\in\k\} ,
\]
where
\[
a\prec^*b\Leftrightarrow a\succ b
\]
is the reverse ordering of $\prec$.  Clearly $\k^*$ is simply
bi-definable with $\k$, so $\k$ satisfies the Ramsey
(resp., ordering) property iff $\k^*$ does. We will formulate now a 
condition on $\k$ which
implies that the only $\k'$ simply bi-definable with $\k$ are $\k$ and
$\k^*$.

Given $\bfa =\langle\bfa_0,\prec\rangle\in\k ,\bfb =\langle\bfb_0,\prec'
\rangle\in\k$ and $a,b\in A_0, c,d\in B_0$ with $a\prec b,c\prec d$, we
say that $(a,b),(c,d)$ {\it have the same type if} there is an isomorphism
of the substructure of $\bfa$ generated by $a,b$ with the substructure
of $\bfb$ generated by $c,d$ which sends $a$ to $c$ and $b$ to $d$.
Equivalently, this means that $(a,b),(c,d)$ satisfy (in $\bfa ,\bfb$ resp.)
exactly the same atomic formulas in $L$.  We denote by $tp_{\bfa}(a,b)$
the {\it type of} $(a,b)$, i.e., the set of all atomic formulas satisfied
by $(a,b)$ in $\bfa$.  Thus $(a,b), (c,d)$ have the same type iff 
$tp_{\bfa}(a,b)=tp_{\bfb}(c,d)$.

We now say that $\k$ satisfies the {\it triangle condition} if for any two 
distinct types $\sigma\neq\tau$, there is $\bfa =\langle\bfa_0,\prec
\rangle\in\k$ and $a\prec b\prec c$ in $A_0$ such that
\[
tp_{\bfa}(a,b)=tp_{\bfa}(b,c)=\sigma, tp_{\bfa}(a,c)=\tau
\]
or
\[
tp_{\bfa} (a,b)=tp_{\bfa} (b,c)=\tau ,tp_{\bfa}(a,c)=\sigma
\]
(i.e., two ``sides'' of the triangle $a,b,c$ have one of types $\sigma ,\tau$
and the third one the other).

\begin{corollary}
\label{9.4}
In the context of \ref{9.2}, if $\k$ satisfies
the triangle condition, then $\k'=\k$ or $\k'=\k^*$.
\end{corollary}

\begin{proof}
In the notation of the proof of \ref{9.2}, it is enough to show that
$\prec'_0=\prec_0$ or $\prec'_0=\prec^*_0$.  If this fails, then we can
find $a,b,c,d$ with $a\prec_0b,c\prec_0d,a\prec'_0b,d\prec'_0c$.  Denote by
$\sigma$ the type of $(a,b)$ (in the substructure of $\bff$ generated by
$a,b$) and by $\tau$ the type of $(c,d)$.  Then $\sigma\neq\tau$, since
$\prec'_0$ is Aut$(\bff )$-invariant.  By the triangle condition, there is
$\bfa\in\k$, which without loss of generality we can assume to be
a substructure of $\bff$, and $\bar a,\bar b,\bar c$ in $A$ such that
$\bar a\prec_0\bar b\prec_0\bar c$
and either $tp_{\bfa}(\bar a,\bar b)= tp_{\bfa}(\bar b,\bar c) =\sigma,\; 
tp_{\bfa}(\bar a,\bar c) =\tau$, in which case we have
that $\bar a\prec'_0\bar b\prec'_0\bar c$ but $\bar c\prec'_0\bar a$, a
contradiction, or $tp_{\bfa}(\bar a,\bar b)=tp_{\bfa}(\bar b,\bar c)=
\tau ,\; tp_{\bfa}(\bar a,\bar c)=\sigma$ in which case we have that $\bar c
\prec'_0\bar b\prec'_0\bar a$ but $\bar a\prec'_0\bar c$, also a 
contradiction.
\end{proof}

As an application, let us show, for example, that the class $\k$ of all
convexly ordered finite equivalence relations is the only (reasonable order
Fra\"{\i}ss\'e) expansion of the class $\e\q$ (of finite equivalence
relations) that satisfies both the Ramsey and ordering properties.  Indeed,
if $\k'$ is any other such class then either $\k'\subseteq\k$, and then we
are done by \ref{9.3}, or else $\k'$ contains a structure $\bfa'=\langle A'_0,
E',\prec'\rangle$, where $\prec'$ is such that there are $a,b,c\in A'_0$
with $(a,c)\in E',(a,b)\not\in E',(b,c)\not\in E'$ and $a\prec' b\prec' c$.
Then we can use $\bfa'$ to witness the triangle condition for $\k'$, so
$\k =\k'$ or $\k =(\k')^*$, by \ref{9.4}.  Then $\k'=\k$ or $\k'=(\k )^*=\k$
and we are done.  (Note that $\k$ itself does not satisfy the triangle
condition).

Another interesting example is the following.  Let $\k_0$ = the class of
finite posets and let $\k$ = the class of all finite $\langle A,\sqsubset ,
\prec\rangle$ where $\langle A,\sqsubset\rangle$ is a poset and $\prec$ is
a linear ordering extending $\sqsubset$.  Then (see ``Addendum'' at the end of this paper) $\k$ satisfies the Ramsey
and ordering properties.  Now it is easy to see that $\k$ also satisfies
the triangle condition, so the only expansions that have both of these
properties are $\k$ and $\k^*$.

Finally, let us point out that in other cases, again in the context of \ref{9.2},
there may be many $\k'$ bi-definable with $\k$.  For example, in case
$\k_0$ = all finite-dimensional vector spaces over a finite field $F$,
we can see, by Section 6, $(\bfc )$, that there are different $\k'$
corresponding to different orderings of $F$.  In such cases one seeks to
classify all expansions of a given $\k_0$ that satisfy
the Ramsey and ordering properties.  Of particular interest are, of course, the cases of finite-dimensional vector spaces and Boolean algebras.

For the case of Boolean algebras, we can show that there are exactly
12 reasonable order Fra\"{\i}ss\'e expansions of the class $\k_0=\b\a$ of
finite Boolean algebras.  These are defined as follows:

Let $\pi_1,\pi_2,\dots ,\pi_6$ be the six permutations of the symbols
$0,1,a$.  For each such $\pi_i$, we define an expansion $\k_i$ of $\k_0$
as follows:  Say, for example, $\pi_i=a,1,0$.  Then $\langle\bfa_0,\prec
\rangle\in\k_i$ iff $\bfa_0$ is a finite Boolean algebra, $\prec$ a linear
ordering of $A_0$, and there is a natural ordering $\prec'$ on $A_0$
such that $\prec ,\prec'$ agree on $A_0\setminus\{0^{\bfa_0},1^{\bfa_0}\}$
but $a_0\prec 1^{\bfa_0}\prec 0^{\bfa_0}$, for all $a_0\in A_0\setminus\{0^{\bfa_0},1^{\bfa_0}
\}$. (Note of course that $0^{\bfa_0}\prec' a_0\prec' 1^{\bfa_0}$, for all $a_0\in A_0\setminus\{0^{\bfa_0},1^{\bfa_0}
\}$.)  Clearly $\k_i$ is simply bi-definable with $\o\b\a$ so each $\k_i$
is a reasonable order Fra\"{\i}ss\'e expansion of $\k_0$ satisfying the Ramsey
and ordering properties.  (Note that if, e.g., $\pi_1=0,1,a$, we have
$\k_1=\o\b\a$.)  Thus the 12 classes $\k_i,\k^*_i ,i=1,\dots ,6,$ are distinct
and satisfy both the Ramsey and ordering properties. Then a canonization argument based on the Dual Ramsey Theorem will give us
the following result which is then used to show that every reasonable Fra\"{\i}ss\'e
expansion of $\k_0$ that satisfies the Ramsey and ordering properties must
be in this list:

Given a finite Boolean algebra $\bfa$, there is a finite Boolean algebra $\bfb$,
with $\bfa\leq\bfb$, such that if $\prec$ is any linear ordering on $\bfb$,
extending the partial ordering of $\bfb$, then there is $\bfa'\in
\begin{pmatrix}\bfb\\ \bfa\end{pmatrix}$ such that $\prec |A'$ is natural.

This order canonization theorem should be compared with those of Ne\v set\v ril-Pr\"omel-R\"odl-Voigt \cite{NPRV} (see also Pr\"omel \cite{Pr}) that deal with canonizing orderings of Boolean lattices rather than Boolean algebras. 

Finally, we can also classify all reasonable order Fra\"{\i}ss\'e expansions of
the class $\k_0=\calv_F$ of all finite-dimensional vector spaces over
a fixed finite field $F$.  It turns out that if card$(F)=q$, then there
are exactly $4(q-2)!$ many such classes described as follows:

Let $F$ be a finite field of cardinality $q$, $F^*=F\setminus\{0\}$
its multiplicative group of non-0 elements.  For each ordering $<$ on $F$,
where $0<F^*$, let $<'$ be the ordering of $F$ such that $<'|F^*=\;<|F^*$ but
$F^*<'0$.  (Here 0 is the ``zero'' of the field $F$).

For each such $<$, let $\k_1(<)$ be the class of all $\langle\bfv_0,
\prec\rangle$, where $\bfv_0$ is a finite-dimensional vector space over $F$,
and $\prec$ is an anti-lexicographical ordering on $V_0$ induced from an ordered basis of $V_0$, by using $<$
for ordering $F$.  If we denote by $\bfzero$ the ``zero'' of $V_0$, clearly
$\bfzero\prec V_0\setminus\{\bfzero\}$.  Let also $\k_2(<)$ be the class
of all $\langle\bfv_0,\prec'\rangle$ such that for some $\prec$ with
$\langle\bfv_0,\prec\rangle\in\k_1(<),\prec |V_0\setminus\{\bfzero\}=\;\prec'
|V_0\setminus\{\bfzero\}$ but $V_0\setminus\{\bfzero\}\prec'\bfzero$.  Then
$\k_1(\prec ),\k_2(\prec )$ are simply bi-definable, so they have the Ramsey
and ordering properties.  For each $<'$ as before, define $\l_1(<')$ to be
the class of all $\langle\bfv_0,\prec\rangle$ such that $\prec$ is an 
anti-lexicographical ordering of $V_0$ induced from an ordered basis of $V_0$
by using $<'$ for ordering $F$.  Note now that $V_0\setminus\{\bfzero\}\prec
\bfzero$.  Define, similarly to the above, $\l_2(\prec')$.  Then $$\l_1(<')^*
=\k_1((<')^*),\l_2(<')^*=\k_2((<')^*),$$
and so $\l_1(<'),\l_2(<')$ also have
the Ramsey and ordering properties. Now using an order canonization theorem based on the Ramsey theorem for vector spaces due to 
Graham-Leeb-Rothschild \cite{GLR} we can show that every reasonable Fra\"{\i}ss\'e order
expansion of $\k_0$ satisfying the Ramsey and ordering properties is one of
$$\k_1(<),\k_2(<),\l_1(<'),\l_2(<').$$However these $4(q-1)!$ classes are 
not distinct.

Consider the action of $F^*$ on $F$ by multiplication and the induced action
of $F^*$ on the set of linear orderings on $F$ with $0<F^*$.  This action
is free, so there are exactly $\tfrac{(q-1)!}{(q-1)}=(q-2)!$ orbits.  Let
$\sim$ be the equivalence relation induced by this action.  Then one can
see that $$\k_1(<_1)=\k_1(<_2)$$
iff $<_1\sim <_2$, and so we get exactly $(q-2)!$ distinct classes of the form 
$\k_1(<)$.  Similarly we get
exactly $(q-2)!$ classes of each of the forms $\k_2(<),\l_1(<'),
\l_2(<')$.  One can finally show that the four collections $$\{\k_1(<)\},
\{\k_2(<)\},\{\l_1(<')\},\{\l_2(<')\}$$ are pairwise disjoint, so we get
exactly $4(q-2)!$ many distinct classes.

In connection with these classification problems, it is worth pointing out that, 
in the context of \ref{9.2} again, 
when the signature $L_0$ is relational and finite, there are only finitely many 
reasonable Fra\"{\i}ss\'{e} order classes that expand a given class in $L_0$ and 
satisfy both the Ramsey and ordering properties. In fact, more generally, using 
the concept of type introduced before \ref{9.4}, it is easy to see that this is even 
true when there are only finitely many types. Indeed, in the notation of the 
proof of \ref{9.2}, note that $\prec'_0$ is completely determined (up to an action by 
an element of Aut($\bff_0$)) by the set 
$$T(\prec'_0) = \{tp_{\bff} (a,b) : a\prec'_0 b\},$$
where $$\bff = \langle \bff_0, \prec_0\rangle$$ (i.e., the set of all types of 
pairs on which $\prec'_0$ agrees with $\prec_0$). But there are only finitely 
many possibilities for $T(\prec'_0)$, so there are only finitely many possibilities 
for $\bff' = \langle \bff _0, \prec'_0\rangle$ and thus for $\k' =$ Age$(\bff')$.

We finally show that among all possible Ramsey expansions $\k$ of a
Fra\"\i ss\'e class $\k_0$, the ones that also have the ordering property,
satisfy an important minimality property.

\begin{theorem}
\label{9.5}
Let $\k_0$ be a Fra\"\i ss\'e class in a
signature $L_0$, let $L=L_0\cup\{<\}$ and let $\k ,\k'$ be reasonable
Fra\"\i ss\'e order classes in $L$ that are expansions of $\k_0$.  If $\k$
satisfies the Ramsey property and $\k'$ satisfies both the Ramsey and
ordering properties, then $\k'\subseteq\k$ up to simple bi-definability,
i.e., there is a reasonable Fra\"\i ss\'e order class $\k ''$ in $L$ which is an
expansion of $\k_0$ such that $\k ''\subseteq\k$ and $\k ''$ is simply
bi-definable with $\k'$.
\end{theorem}

\begin{proof}
Using the notation of the proof of \ref{9.2}, we have, by \ref{7.5}, that
the $G_0$-space $X_{\k'}$ is isomorphic to a minimal subflow $X$ of $X_{\k}$.
Let $\pi :X_{\k'}\rightarrow X$ be an isomorphism and let $\pi (\prec'_0)=\prec''_0
\in X$.  Then $$\text{Aut}(\langle\bff_0,\prec'_0\rangle )=\text {Aut}(\langle\bff_0,
\prec''_0\rangle ),$$ so that there is a simple formula $\varphi''$ that defines
$\prec''_0$ in $\langle\bff_0,\prec'_0\rangle$ and a simple formula $\varphi'$
that defines $\prec'_0$ in $\langle\bff_0,\prec''_0\rangle$.  From this it easily
follows that $\langle\bff_0,\prec''_0\rangle$ is ultrahomogeneous.  Let $\k''
={\rm Age}(\langle\bff_0,\prec''_0\rangle )$.  Then $\k''\subseteq\k$ is a
reasonable Fra\"\i ss\'e order class in $L$, which is an expansion of $\k_0$.
As in the proof of \ref{9.2}, it is simply bi-definable with $\k'$.
\end{proof}

Thus the expansions of a given $\k_0$ that have both the Ramsey and ordering
properties are the smallest, up to simple bi-definability, expansions that
have the Ramsey property.  The preceding argument also suggests an approach
to showing that, if there is an expansion $\k$ with the Ramsey property, then
there is one with both Ramsey and ordering properties:  Simply pick $\prec'
\in X_{\k}$ such that $\overline{G_0\cdot \prec'}$ is minimal and let $\k'=
{\rm Age}(\langle\bff_0,\prec'\rangle )$.  Then try to show that $\k'$ has both the Ramsey and ordering properties. We will see in the next section that this approach indeed works.

\section{Ramsey degrees}

Consider the class of finite graphs $\g\r$.  Although this 
class does not satisfy the Ramsey property, there is still an important
Ramsey-type result that holds for finite graphs:  For each finite graph
$\bfa_0\in\g\r$, there is a finite number $t$ such that for all $\bfa_0
\leq\bfb_0\in\g\r$, there is $\bfc_0\geq\bfb_0$ in $\g\r$ with the following property:
For any coloring $c:\begin{pmatrix}\bfc_0\\ \bfa_0\end{pmatrix}\rightarrow
\{1,\dots ,k\}$ (of any number of colors), there is $\bfb'_0\in
\begin{pmatrix}\bfc_0\\ \bfb_0\end{pmatrix}$ such that $c$ takes at most $t$
colors on $\begin{pmatrix}\bfb'_0\\ \bfa_0\end{pmatrix}$.  Moreover, one
can explicitly compute the least number $t$ with that property.

Such results for graphs and other classes of finite structures have
existed implicitly or explicitly, in the literature for some time now
(see, e.g., Ne\v set\v ril-R\"odl \cite{NR77}, \cite{NR83}, 
Abramson-Harrington \cite{AH78},
and Fouch\'e \cite{Fo97}, \cite{Fo98}, \cite{Fo99}, 
\cite{Fo99a}).  We will present
here a version of this theory in the general context of our paper and 
notice some interesting connections with the results in Section 9 leading,
in particular, to the proof of the result mentioned at the end of Sections
7 and 9.

Throughout this section, $\k_0$ {\it will be a Fra\"\i ss\'e class in a 
signature $L_0$, and $\k$ a hereditary order class in $L=L_0
\cup\{<\}$ with $\k_0=\k |L$, i.e., $\k$ is an expansion of $\k_0$.}

For $\bfa_0\in\k_0$, let
\[
X^{\bfa_0}_{\k}=\{\prec:\prec\text{ is a linear ordering on }A_0\text{ and }
\langle\bfa_0,\prec\rangle\in\k\}.
\]
Thus $X^{\bfa_0}_{\k}\neq\emptyset$ by assumption and, since $\k$ is closed
under isomorphism, Aut$(\bfa_0)$ acts on $X^{A_0}_\k$ in the
obvious way.  If $\prec\in
X^{\bfa_0}_{\k}$, we call its orbit under this action the $(\bfa_0-\k$--)
{\it pattern} of $\prec$.  Thus the pattern of $\prec$ is the set of all $\prec'$
such that $\langle \bfa_0,\prec'\rangle\cong\langle\bfa_0,\prec\rangle$.
Put also
\[
\begin{array}{ll}
t_\k(\bfa_0)&=\text{the cardinality of the set of }\bfa_0-\k-\text{patterns}\\
&={\rm card} (X^{\bfa_0}_{\k}/{\rm Aut} (\bfa_0))\\
&=\frac{{\rm card} (X^{\bfa_0}_{\k})}{{\rm card}({\rm Aut} (\bfa_0)) },
\end{array}
\]
since clearly Aut$(\bfa_0)$ acts freely on $X^{\bfa_0}_{\k}$, and thus
every orbit has the same cardinality as Aut$(\bfa_0)$.

Finally, for $\bfa_0\leq\bfb_0\leq\bfc_0$ in $\k_0$, let
\[
\bfc_0\rightarrow (\bfb_0)^{\bfa_0}_{k,t}
\]
mean that for every coloring $c:\begin{pmatrix}\bfc_0\\ \bfa_0\end{pmatrix}
\rightarrow\{1,\dots ,k\}$ there is $\bfb'_0\in\begin{pmatrix}\bfc_0\\ 
\bfb_0\end{pmatrix}$ such that $c$ on $\begin{pmatrix}\bfb'_0\\ \bfa_0
\end{pmatrix}$ takes at most $t$ values.  Thus
\[
\bfc_0\rightarrow (\bfb_0)^{\bfa_0}_{k,1}\Leftrightarrow\bfc_0\rightarrow
(\bfb_0)^{\bfa_0}_{k}.
\]

Also for $\bfa_0\leq\bfb_0\in\k_0$, and $\langle\bfb_0,
<_{\bfb_0}\rangle
\in\k$, let
\[
\begin{array}{ll}
t_\k(\bfa_0,\bfb_0, <_{\bfb_0})=&\text{the number of }\bfa_0-\k-\text{patterns
of }\\
&\prec\in\bfx^{\bfa_0}_\k\text{ such that }\langle\bfa_0,\prec\rangle\leq
\langle\bfb_0, <_{\bfb_0}\rangle .
\end{array}
\]
Clearly $t_\k(\bfa_0,\bfb_0,{\boldsymbol <}_{B_0})\leq t_{\k'}(\bfa_0)$, for any
hereditary class $\k'\subseteq\k$ with 
$\k' |L_0=\k_0, \langle\bfb_0, <_{\bfb_0}\rangle \in \k' $.

\begin{proposition}
\label{10.1}
If $\k$ has the Ramsey property, then
for $\bfa_0\leq\bfb_0$ in $\k_0$, and 
$\langle\bfb_0, <_{\bfb_0}\rangle \in \k , k\geq 2$, there is
$\bfc_0\geq\bfb_0$ in $\k_0$
with
\[
\bfc_0\rightarrow (\bfb_0)^{\bfa_0}_{k, t_\k (\bfa_0, \bfb_0, <_{\bfb_0})},
\]
and thus also
\[
\bfc_0\rightarrow (\bfb_0)^{\bfa_0}_{k, t_{\k'} (\bfa_0)},
\]
for any hereditary class $\k'\subseteq\k$ which is an expansion of $\k_0$.
\end{proposition}

\begin{proof}
Put $t=t_\k(\bfa_0, \bfb_0, <_{\bfb_0})$.  Choose representatives $<_1,
\dots ,<_t$ for the $\bfa_0-\k -$patterns realized by $\prec$ with 
$\langle \bfa_0, \prec\rangle \leq \langle \bfb_0, <_{\bfb_0}\rangle$.  
We will define inductively
$\langle \bfd_i,<_{\bfd_i}\rangle\in\k , 0\leq i\leq t$ as follows:

Let $\langle\bfd_0,<_{\bfd_0}\rangle =\langle\bfb_0,
<_{\bfb_0}\rangle$.  Then for $i\geq 1$, let  $\langle\bfd_i,<_{\bfd_i}\rangle\in\k$
be such that
\[
\langle\bfd_i,<_{\bfd_i}\rangle\rightarrow (\langle\bfd_{i-1},<_{\bfd_{i-1}}
\rangle )^{\langle \bfa_0,<_i\rangle}_{k}.
\]
Take now $\bfc_0=\bfd_t$.  We claim that this works:

Let $c:\begin{pmatrix}\bfc_0\\ A_0\end{pmatrix}\rightarrow\{1,\dots ,k\}$
be a coloring.  Then let
\[
c_t:\begin{pmatrix}\langle\bfc_0,<_{\bfd_t}\rangle \\ \langle\bfa_0,<_t
\rangle\end{pmatrix}\rightarrow\{1,\dots ,k\}
\]
be defined by
\[
c_t:(\langle\bfa'_0,<'\rangle )=c(\bfa'_0),
\]
where $\langle\bfa'_0,<'\rangle =\langle\bfa'_0,<_{\bfd_t}|A'_0\rangle
\cong\langle\bfa_0,<_t\rangle$.  There is a copy $\langle\bfd'_{t-1},
\break <_{\bfd'_{t-1}}\rangle$ of $\langle\bfd_{t-1},<_{\bfd_{t-1}}\rangle$
in $\langle\bfd_t,<_{\bfd}\rangle$ and a color $1\leq k_t\leq k$ such that
the $c_t$-color of any copy of $\langle\bfa_0,<_t\rangle$ in $\langle
\bfd'_{t-1},<_{\bfd'_{t-1}}\rangle$ is equal to $k_t$. Iterate now this process
starting with $\langle\bfd'_{t-1},<_{\bfd'_{t-1}}\rangle$ and the coloring $c_{t-1}:
\begin{pmatrix}\langle\bfd'_{t-1},<_{\bfd_{t-1}}\rangle\\ \langle
\bfa_0,<_{t-1}\rangle\end{pmatrix}\rightarrow\{1,\dots ,k\}$, given by $c_{t-1}
(\langle\bfa'_0,<'\rangle )=c(\bfa'_0)$, where $$\langle\bfa'_0,<'\rangle
=\langle\bfa'_0,<_{\bfd'_{t-1}}|A'_0\rangle\cong\langle\bfa_0,<_{t-1}
\rangle,$$to find a copy $\langle\bfd'_{t-2},<_{\bfd'_{t-2}}\rangle$ of
$\langle\bfd_{t-2},<_{\bfd_{t-2}}\rangle$ in $\langle\bfd'_{t-1},
<_{\bfd'_{t-1}}\rangle$, and thus in the structure $\langle\bfd_t,<_{\bfd_t}\rangle$,
and a color $1\leq k_{t-1}\leq k$ such that the $c_{t-1}$-color of any
copy of $\langle\bfa_0,<_{t-1}\rangle$ in $\langle\bfd'_{t-2},<_{\bfd'_{t-2}}
\rangle$ is equal to $k_{t-1}$, etc.  After $t$ steps, we get a copy of
$\langle\bfd_0,<_{\bfd_0}\rangle =\langle\bfb_0,<_{\bfb_0}\rangle$ in
$\langle\bfd'_1,<_{\bfd'_1}\rangle$ and a color $1\leq k_1\leq k$ that
works for copies of $\langle\bfa_0,<_1\rangle$.  Call this copy $\langle
\bfb'_0,<_{\bfb'_0}\rangle$.  Then clearly $c$ on $\begin{pmatrix}\bfb'_0\\
\bfa_0\end{pmatrix}$ takes at most the values $\{1,\dots , k_t\}$, because if
$\bfa'_0\in\begin{pmatrix}\bfb'_0\\ \bfa_0\end{pmatrix}$ and the pattern of 
$<_t|A'_0 = <'_{\bfb_0}|A'_0$ is $<_i$, then $$\langle\bfa'_0,<_t|A'_0\rangle\in\begin{pmatrix}
\langle\bfd'_i,<_{\bfd'_i}\rangle\\ \langle\bfa_0,<_i\rangle\end{pmatrix},$$
so $c_i(\langle\bfa'_0,<'\rangle )=c(\bfa'_0)=k_i$.
\end{proof}

Actually the preceding proof also establishes the following.
\medskip

\noindent{\bf Proposition \ref{10.1}$'$} {\it  If $\k$ has the Ramsey property,
then for $\bfa_0\leq\bfb_0$ in $\k_0 ,\ k\geq 2$, there is $\bfc_0\geq\bfb_0$
and a linear ordering $\prec$ on $\bfc_0$ with the property that  
$\langle\bfc_0,\prec\rangle\in\k$
and for any coloring $c:\begin{pmatrix}\bfc_0\\ \bfa_0\end{pmatrix}
\rightarrow\{1,\dots ,k\}$, there is $\bfb'_0\in\begin{pmatrix}\bfc_0\\
\bfb_0\end{pmatrix}$ so that the color $c(\bfa'_0)$, for $\bfa'_0\in
\begin{pmatrix}\bfb'_0\\ \bfa_0\end{pmatrix}$, depends only on the pattern of
$\prec |A'_0$.} 

\medskip
If now $\k$ has also the ordering property, the number $t_\k (\bfa_0)$ is
best possible.

\begin{proposition}
\label{10.2}
If $\k$ has the Ramsey and ordering
properties, then for any $\bfa_0\in\k_0$, there is $\bfb_0\geq\bfa_0$ in
$\k_0$ such that for all $\bfc_0\geq\bfb_0$ in $\k_0$ there is a coloring
$c:\begin{pmatrix}\bfc_0\\ \bfa_0\end{pmatrix}\rightarrow\{1,\dots ,t_\k
(\bfa_0)\}$, such that for any $\bfb'_0\in\begin{pmatrix}\bfc_0\\ \bfb_0
\end{pmatrix},\ c$ takes all the values $1,\dots , t_\k (\bfa_0)$ on 
$\begin{pmatrix}\bfb'_0\\ \bfa_0\end{pmatrix}$.
\end{proposition}

\begin{proof}
Let, by the ordering property, $\bfb_0\in\k_0$ be such
that for any ordering $\prec$ on $A_0$, and any ordering $\prec'$ on $B_0$ with
$\langle\bfa_0,\prec\rangle ,\langle\bfb_0,\prec'\rangle\in\k$ we have $\langle
\bfa_0,\prec\rangle\leq\langle\bfb_0,\prec'\rangle$.  Another way of saying this
is that as $\bfa'_0$ varies over $\begin{pmatrix}\bfb_0\\ \bfa_0
\end{pmatrix}, \prec'|A_0$ realizes all possible patterns.

Let now $\bfc_0\in\k_0, \bfc_0 \geq \bfb_0,$ and $<_{\bfc_0}$, a linear ordering on $C_0$, be
such that $$\langle\bfc_0,<_{\bfc_0}\rangle\in\k.$$Define $c:\begin{pmatrix}
\bfc_0\\ \bfa_0\end{pmatrix}\rightarrow\{1,\dots ,t_\k (\bfa_0\}$ by
enumerating the set of patterns as $p_1,\dots , p_{t_\k (\bfa_0)}$ and
letting, for $\bfa'_0\in\begin{pmatrix}\bfc_0\\ \bfa_0\end{pmatrix},
c(\bfa'_0)=i$, where the pattern of $\langle\bfa'_0,<_{\bfc_0}|A'_0
\rangle$ is $p_i$.  Then, by the above, if $\bfb'_0\in\begin{pmatrix}
\bfc_0\\ \bfb_0\end{pmatrix}$, clearly $c$ takes all values $1,\dots ,t_\k
(\bfa_0)$ on $\begin{pmatrix}\bfb'_0\\ \bfa_0\end{pmatrix}$.
\end{proof}

\begin{corollary}
\label{10.3}
If $\k$ has the Ramsey and ordering
properties, then for $\bfa_0\in\k_0, t_\k (\bfa_0)$ is the least number 
$t$ such that for any $\bfa_0\leq\bfb_0$ in $\k_0 ,k\geq 2$, there is
$\bfc_0\geq\bfb_0$ in $\k_0$ with $\bfc_0\rightarrow (\bfb_0)^{\bfa_0}_{k,
t}$.
\qed
\end{corollary}

In particular, this shows that $t_\k (\bfa_0)$ and card($X^{\bfa_0}_\k$) are 
independent of $\k$, as long as $\k$ has the Ramsey and ordering properties.
So one has a uniqueness property of expansions of $\k_0$ that have both the
Ramsey and ordering properties.  Theorem \ref{9.2} provides a much stronger 
uniqueness property, which immediately implies this.

Also if $\k$ has the Ramsey property and $\k'$ has both the Ramsey and
ordering properties, then $t_{\k'} (\bfa_0)\leq t_{\k}(\bfa_0)$ and
card$(X^{\bfa_0}_{\k'} )\leq {\rm card}(X^{\bfa_0}_{\k})$. Of course Theorem \ref{9.5} provides a stronger minimality property of such $\k'$. In particular,
if it happens, as in many examples that we have seen before, that $\k =\k_0
*\l\o$ has both the Ramsey and ordering properties, then clearly $\k$ is
also the unique expansion of $\k_0$ that has the Ramsey property.

For each class $\k_0$ and $\bfa_0 \in \k_0$, let $t(\bfa_0, \k_0)$ be the least $t$, if it exists, such that for any $\bfa_0\leq\bfb_0$ in $\k_0 ,k\geq 2$, there is
$\bfc_0\geq\bfb_0$ in $\k_0$ with 
$$\bfc_0\rightarrow (\bfb_0)^{\bfa_0}_{k,t}.$$
Otherwise, let $t(\bfa_0, \k_0)=\infty $. Following Fouch\'{e} 
\cite{Fo98}, call $t(\bfa_0, \k_0)$ the {\it Ramsey degree} of $\bfa_0$ (in $\k_0$). We have seen that if $\k_0$ admits an expansion $\k$ with the Ramsey and ordering properties, then $t(\bfa_0, \k_0) =t_{\k}(\bfa_0).$

For example, if $\k_0 =\g\r$, then $t(\bfa_0, \k_0)=\frac{{\rm card}(\bfa_0)!}
{{\rm card}({\rm Aut}(\bfa_0))}$.  For instance, if $\bfa_0=K_n$, the
complete graph on $n$ vertices, or $\bfa_0=$ the complement of $K_n$,
then $t(\bfa_0, \k_0)=1$, i.e., 
$\bfa_0$ satisfies Ramsey's Theorem, and moreover the
$K_n$ and their complements are the only graphs that have this property
(see Ne\v set\v ril-R\"odl \cite{NR77}, \cite{NR83}).  For other calculations of this
sort, see \cite{Fo97},\cite{Fo98},\cite{Fo99},\cite{Fo99a}.

Notice also that $t_\k(\bfa_0) =1$ for all $\bfa_0 \in \k_0$ iff $\k$ is order 
forgetful (see \ref{5.5}).

Actually the preceding results admit appropriate converses.  First, for
Proposition \ref{10.1} we have:

\begin{proposition}
\label{10.4}
Suppose $\k$ is reasonable and has the ordering property
and for any $\bfa_0\leq\bfb_0$ in $\k_0 ,k\geq 2$, there is $\bfc_0\geq\bfb_0$
in $\k_0$ with
\[
\bfc_0\rightarrow (\bfb_0)^{\bfa_0}_{k,t_\k (\bfa_0)}.
\]
Then $\k$ has the Ramsey property. 
In particular, any reasonable $\k$ with the ordering property, that is 
contained in an expansion of $\k_0$ with the Ramsey property, also has 
the Ramsey property.
\end{proposition}

\begin{proof}
Suppose $\langle\bfa_0,\prec\rangle ,\langle\bfb_0,\prec'\rangle\in
\k$ are given with $$\langle\bfa_0,\prec\rangle\leq\langle\bfb_0,\prec'\rangle.$$
We need to find $\langle\bfc_0,\prec''\rangle\geq\langle\bfb_0,\prec'\rangle$ in
$\k$ such that
\[
\langle\bfc_0,\prec''\rangle\rightarrow (\langle\bfb_0,\prec'\rangle ))^{\langle\bfa_0,
\prec\rangle}_2 .
\]

First, by the ordering property, we can find $\bfb_1\in\k_0$ such that for
any ordering $\prec_0$ on $B_0$ and any ordering $\prec_1$ on $B_1$ with $\langle
\bfb_0,\prec_0\rangle ,\langle\bfb_1,\prec_1\rangle\in\k$, we have $\langle\bfb_0,
\prec_0\rangle\leq\langle\bfb_1,\prec_1\rangle$. Since $\k$ is reasonable, this
implies that for any ordering $\prec_0$ on $A_0$ and any ordering $\prec_1$ on $B_1$
we also have $\langle\bfa_0,\prec_0\rangle\leq\langle\bfb_1,\prec_1\rangle$.

By hypothesis, there is $\bfc_0\geq\bfb_1$ in $\k_0$ such that
\[
\bfc_0\rightarrow (\bfb_1)^{\bfa_0}_{t_\k (\bfa_0)+1,t_\k (\bfa_0)}.
\]
Fix an ordering $\prec''$ on $\bfc_0$ with $\langle\bfc_0,\prec''\rangle\in\k$.
We claim that $\langle\bfc_0,\prec''\rangle$
works.  Clearly, $\langle\bfb_0,\prec'\rangle\leq\langle\bfc_0,\prec''
\rangle$.

Consider now a coloring
\[
c:\begin{pmatrix}\langle\bfc_0,\prec''\rangle\\ \langle\bfa_0,\prec\rangle
\end{pmatrix}\rightarrow\{1,2\}.
\]
Use this to define a coloring
\[
c':\begin{pmatrix}\bfc_0\\ \bfa_0\end{pmatrix}\rightarrow
\{0\}\cup \{p_1,\dots , p_{t_\k (\bfa_0)}\},
\]
where $p_1,\dots ,p_{t_\k (A_0)}$ enumerate the $\bfa_0-\k -$patterns, as
follows:  Let $\bfa'_0\in\begin{pmatrix}\bfc_0\\ \bfa_0\end{pmatrix}$.
Then
\[
\begin{array}{ll}
c'(\bfa'_0)&=\text{the pattern of }\prec''|A'_0,\text{if this pattern
is different from that of}\\
&\prec ,\text{ or if }c(\langle\bfa'_0,\prec''|A'_0\rangle )=1,\\
&=0,\text{ otherwise}
\end{array}
\]

Let now $\bfb'_1\in\begin{pmatrix}\bfc_0\\ \bfb_1\end{pmatrix}$ be such
that $c'$ takes at most $t_\k(\bfa_0)$ values on $\begin{pmatrix}\bfb'_1\\
\bfa_0\end{pmatrix}$.  Since any pattern is realized among the $\prec''|A'_0$,
where $\bfa'_0\in\begin{pmatrix}\bfb'_1\\ \bfa_0\end{pmatrix}$, it follows
that $c$ is constant on $\begin{pmatrix}\langle\bfb'_1,\prec''|B'_1
\rangle )\\ \langle\bfa_0,\prec\rangle\end{pmatrix}$.  Since $$\langle
B_0,\prec'\rangle\leq\langle\bfb'_1,\prec'' |B'_1\rangle,$$ we are done.

The last assertion follows by also using \ref{10.1}.
\end{proof}

To formulate a converse to Corollary \ref{10.3}, let us first define a local version
of the ordering property.

We say that $\k$ satisfies the {\it ordering property at} $\bfa_0\in\k_0$,
if there is $\bfb_0\geq\bfa_0$ in $\k$ such that for every orderings $\prec,\prec'$
on $A_0,B_0$, resp., with $\langle\bfa_0,\prec\rangle ,\langle\bfb_0,\prec'\rangle\in\k$,
we have $\langle\bfa_0,\prec\rangle\leq\langle B_0,\prec'\rangle$.  Thus $\k$ has
the ordering property iff $\k$ has the ordering property at each $\bfa_0
\in\k_0$.

We now have the following converse to Corollary \ref{10.3}.

\begin{proposition}
\label{10.5}
Assume $\k$ has the Ramsey property.
Then the following are equivalent:

(i) $\k$ has the ordering property at $\bfa_0\in\k$.

(ii) $t(\bfa_0,\k_0)=t_\k (\bfa_0)$.
\end{proposition}

\begin{proof}
(i) $\Rightarrow$ (ii):  As in the proof of Proposition \ref{10.2}.

(ii) $\Rightarrow$ (i):  Assume (i) fails, so that for every $\bfb_0\geq
\bfa_0$ in $\k_0$ there are $\prec,\prec'$ with $\langle\bfa_0,\prec\rangle\in\k ,
\langle\bfb_0,\prec'\rangle\in\k$ but $\langle\bfa_0,\prec\rangle\not\leq\langle
\bfb_0,\prec'\rangle$. Let $t=t_\k (\bfa_0)$. Then $t_\k (\bfa_0,\bfb_0,\prec') 
\leq t-1$, so we are done by \ref{10.1}.
\end{proof}

\begin{corollary}
\label{10.6}
If $\k$ has the Ramsey property, then
the following are equivalent:

(i) $\k$ has the ordering property.

(ii) For every $\bfa_0\in\k , t(\bfa_0, \k_0)=t_\k (\bfa_0)$.
\qed
\end{corollary}

We apply these facts to prove the result mentioned in the
last paragraph of Section 9.

\begin{theorem}
\label{10.7}
Let $\k_0$ be a Fra\"\i ss\'e
class in a signature $L_0$, and assume that
$\k$ is a reasonable Fra\"\i ss\'e order class in $L=L_0\cup\{<\}$ which is an expansion
of $\k_0$ and satisfies the Ramsey property.  Then there is a 
reasonable Fra\"\i ss\'e order class $\k'\subseteq\k$, which is an
expansion of $\k_0$, and satisfies both the Ramsey and ordering
properties.
\end{theorem}

\begin{proof}
We will use the notation of the proof of \ref{9.2}.  Let
$X'\subseteq X_\k$ be a minimal $G_0$-subflow of the $G_0$-flow
$X_\k$, and let $\prec'\in X'$, so that $X'=\overline{G_0\cdot
\prec'}$.  Let $\k'=$ Age$(\langle\bff_0,\prec'\rangle )\subseteq
\k$.  We will show that this works.  Clearly $\k'|L_0=\k_0$,
i.e., $\k'$ is an expansion of $\k_0$, and $\k'$ is hereditary
and satisfies JEP.  It is also easy to see that $\k'$ is
reasonable.  Note that $X'=X_{\k'}=\{\prec:\prec$ is a linear ordering
on $F_0$ and for every finite $\bfb_0\subseteq\bff_0 ,\langle
\bfb_0,\prec|B_0\rangle\in\k'\}$.  It follows then from the proof
of \ref{7.4}, (i) $\Rightarrow$ (ii), that $\k'$ satisfies the ordering
property.  Thus to verify that $\k'$ has the Ramsey property
it is enough, by \ref{10.4}, to check that for any $\bfa_0\leq\bfb_0$
in $\k_0, k\geq 2$, there is $\bfc_0\geq\bfb_0$ in $\k_0$ with
$\bfc_0\rightarrow (\bfb_0)^{\bfa_0}_{k,t_{\k'}(\bfa_0)}$.
This follows from \ref{10.1}.  Finally, from Section 3 (second to
last paragraph) it follows that $\k'$ has the amalgamation property,
and this completes the proof.
\end{proof}

We can also use similar ideas to prove the converse of \ref{7.5} (ii).

\begin{theorem}
\label{10.8}
Let $L\supseteq\{<\}$ be a signature,
$L_0=L\setminus\{<\},\k$ a reasonable Fra\"{\i}ss\'e order class in $L$, and
let $\k_0=\k |L_0$ and $\bff ={\rm Flim}(\k),\bff_0={\rm Flim}
(\k_0)=\bff |L_0$.  Let $G_0={\rm Aut}(\bff_0),G={\rm Aut}(\bff )$ and let
$X_\k$ be the set of linear orderings on $F(=F_0)$ which are $\k$-admissible.
Then the following are equivalent:

(i) $\k$ has the Ramsey and ordering properties.

(ii) $X_\k$ is the universal minimal flow of $G_0$.
\end{theorem}

\begin{proof}
(i) $\Rightarrow$ (ii) is \ref{7.5} (ii).

(ii) $\Rightarrow$ (i): Since $X_\k$ is a minimal flow, $\k$ has the
ordering property by \ref{7.4}.  So it is enough to verify the hypothesis of
\ref{10.4}.  By the usual ultrafilter argument, as in the proof of \ref{4.5}, it is
enough to show for given $\bfa_0\leq\bfb_0$ in $\k ,k\geq 2$, and $c:
\begin{pmatrix}\bff_0\\ \bfa_0\end{pmatrix}\rightarrow\{1,\dots ,k\}$, that
there is $\bfb'_0\in\begin{pmatrix}\bff_0\\ \bfb_0\end{pmatrix}$ such that
$c$ on $\begin{pmatrix}\bfb'_0\\ \bfa_0\end{pmatrix}$ obtains at most 
$t_\k (\bfa_0)$ many values.  To see this, consider the $G_0$-flow $\{1,\dots
,k\}^{\begin{pmatrix}\bff_0\\ \bfa_0\end{pmatrix}}$, where $G_0$ acts on
this space in the usual way:  $g\cdot\gamma (\bfa'_0)=\gamma (g^{-1}
(\bfa'_0))$.  Let $X=\overline{G_0\cdot c}$.  Then there is a homomorphism $\pi :X_\k
\rightarrow X$.  Put $\pi (\prec_0)=c_0$, where $\bff =\langle\bff_0,
\prec_0\rangle$.  Since $G$ stabilizes $\prec_0, G$ also stabilizes $c_0$.
>From this it easily follows that the color $c_0(\bfa'_0)$ for any $\bfa'_0
\in\begin{pmatrix}\bff_0\\ \bfa_0\end{pmatrix}$ depends only on the pattern
of $\prec_0|A'_0$, thus $c_0$ takes at most $t_\k (\bfa_0)$ values.  Fix
now $\bfb''_0\in\begin{pmatrix}\bff_0\\ \bfb_0\end{pmatrix}$.  Then there
is $g\in G_0$ such that $g\cdot c\biggl\vert\begin{pmatrix}\bfb''_0\\ \bfa_0
\end{pmatrix}=c_0\biggl\vert\begin{pmatrix}\bfb''_0\\ \bfa_0\end{pmatrix}$.
Let $\bfb'_0=g^{-1}(\bfb''_0)$.  Then $c$ on $\begin{pmatrix}\bfb'_0\\ 
\bfa_0\end{pmatrix}$ obtains at most $t_\k (\bfa_0)$ many values and we
are done.
\end{proof}

There are some natural questions that are suggested by the preceding facts.
First, recall the question that we raised in Section 7 of understanding when $\k_0$ has an expansion $\k$ with the Ramsey property. A necessary condition
for the existence of such a $\k$ is that $t(\bfa_0, \k_0)<\infty$, for all $\bfa_0
\in\k_0$.  Is that actually a necessary and sufficient condition?  We should
point out that we do not know an example of a $\k_0$ for which there
is an $\bfa_0\in\k_0$ with $t(\bfa_0,\k_0)=\infty$, although one should surely exist.

Second, in the context of \ref{10.8}, le us say that a $G_0$-flow is {\it
universal} if it can be mapped homomorphically to any other $G_0$-flow.
Thus the universal minimal flow of $G_0$ is a minimal and universal
$G_0$-flow.  We have seen in \ref{7.4} that $X_\k$ is a minimal $G_0$-flow
iff $\k$ has the ordering property.  Can we strengthen \ref{10.8} by showing that 
$X_\k$ is a universal $G_0$-flow iff $\k$ has the Ramsey property?

Finally, recall that if $\k_0$ has an expansion $\k$ with the Ramsey and
ordering properties, then the automorphism group $G_0={\rm Aut}(\bff_0)$,
where $\bff_0={\rm Flim}(\k_0)$, has universal minimal flow $X_\k$, which 
is the inverse limit of the family $\{X^{\bfa_0}_\k\}$, where $\bfa_0$
varies over finite substructures of $\bff_0$ ordered under inclusion.  Thus
if $X_\k$ is finite, say of cardinality $n$, then clearly card$(X^{\bfa_0}_\k)
\leq n$, so, in particular, card(Aut$(\bfa_0))\leq n$.  It is easy to find
examples where card$(X_\k)=1$, i.e., $G_0$ is extremely amenable.  Take,
for instance, $\k_0$ to be a Fra\"\i ss\'e order class with the Ramsey
and ordering properties in a language $L_0\supseteq\{<\}$ and let $\k$
consist of all structures of the form $\bfa =\langle\bfa_0,\prec'\rangle$, with 
$\prec'=<^{\bfa_0}$.  However, we do not know examples of $\k_0,\k$ as above with
$2\leq{\rm card}(X_\k )<\infty$.  Also note that if $\k$ is closed under
products, then $X_\k$ cannot be finite, as $\sup_{\bfa_0\in\k_0}{\rm card}
( {\rm Aut} (\bfa_0))=\infty$.

\section{Concluding remarks and problems}

{\bf (A)} One of the two main ingredients in our proofs of extreme
amenability of automorphism groups are the results of the
corresponding structural Ramsey theory.  It is therefore natural
to pose the following problem a solution of which could enhance
the already existing tradition of using the methods of topological
dynamics to prove results of Ramsey theory.

\begin{problem}
\label{11.1}
Find alternate proofs (that use the
methods of topological dynamics itself as well as the intrinsic geometry of
the acting groups) that the automorphism groups of
any of the following structures are extremely amenable:

(i) The rationals with the usual ordering.

(ii) The random ordered graph.

(iii) The random $K_n$-tree ordered graph, ($n=3,4,\dots$).

(iv) The random $\a$-free ordered hypergraph of type $L_0$, for any class
$\a$ of irreducible finite hypergraphs of type $L_0$.

(v) The ordered rational Urysohn space.

(vi) The $\aleph_0$-dimensional vector space over a finite field with the
canonical ordering.

(vii) The countable atomless Boolean algebra with the canonical
ordering.
\end{problem}

\noindent
{\bf (B)}
In connection with the extreme amenability of $U(\ell_2)$ and
the group of isometries of $\ell_2$ (see the end of Section 6), we would
like to ask whether a new proof of these results can be found based on
Ramsey theory, as was done in Section 6, ({\bf E}) for the
isometry group of the Urysohn space.
\medskip

\noindent {\bf (C)}  Topological groups $G$ isomorphic to closed
subgroups of the infinite symmetric group $S_\infty$, or, equivalently, 
those $G$ whose open subgroups form a nbhd basis at the identity 
(groups with small open subgroups) have played a
leading role in this article. One can meaningfully extend to
homogeneous spaces of such groups the concept of Ramsey degree as
follows.

For each bounded nonempty set $A$ in a Euclidean space and $\epsilon>0$, let
$N(\epsilon ,A)$ be the {\it covering number,} that is,
the smallest number of sets of diameter $\leq \epsilon $ that can cover $A$.
Similarly, for a bounded function
$f$ from some
set $X$ into a Euclidean space, we define $N(\epsilon ,f) = N(\epsilon ,
{\mathrm{range}} (f))$.

Let $G$ be a topological group and let $H$ be a subgroup. We define the {\it small
oscillation degree} $n(G,H)$ of $G,H$ to be the smallest number $t$ such that:

(*) for every finite subset $F$ of $G$ , $\epsilon >0$, and bounded
left uniformly
continuous $f$ from $G/H$ to some Euclidean space, we can find some $h$ in $G$ with
$N(\epsilon ,f' | hF) \leq t$, where $f'$ is the lift of $f$ to $G$.

If no such $t$ exists, we put $n(G,H) = \infty$. Put $n(G) = n(G, \{1\})$. Note
that
\[n(G,H) \leq n(G,H') \leq n(G),\]
if $H'$ is a subgroup of $H$.

If $H$ is open in $G$, then $n(G,H)) \leq t$ if and only if for
every coloring $c$ of $G/H$ with any finite number of colors and
every finite $F$ included in $G$, there is an $h$ in $G$ such that
$c$ on $hF$ takes at most $t$ colors. One can prove that for a
group $G$ with small open subgroups, $n(G) = \sup\{n(G,H): H
\mbox{ is an open subgroup of }G\}$. Consequently, such a group
$G$ is extremely amenable iff $n(G) = 1$. (In the case of a discrete semigroup $G$ this was established by 
Mitchell \cite{Mit}, while in a general case of a group acting on a 
metric space the latter equivalence is due to Gromov and Milman 
\cite{GroM}, \cite{M88}, where a suitable extension of the condition 
$n(G)=1$ was studied under the name of {\it concentration property.})

Let now $\k$ be a Fra\"{\i}ss\'e class, $\bff$ its Fra\"{\i}ss\'e limit and $G =
{\mathrm{Aut}}(\bff)$. For each finite substructure $\bfa$ of
$\bff$, let $H(\bfa)$ be the (setwise) stabilizer of the domain of
$\bfa$ in the action of $G$ on the finite substructures of $\bff$.
Then one can see that
 \[n(G,H(\bfa)) = t(\bfa ,\k)= \mbox{ the Ramsey degree of $\bfa$ (in the class $\k$)}.\]

\noindent {\bf (D)} We have seen that the extreme amenability of
the automorphism group of an ultrahomogeneous ordered countable structure
is equivalent to a corresponding finite Ramsey-theoretic result.
This leads us to the following natural problem.

\begin{problem}
\label{11.2}
In each of the cases (i)--(vii) of
Problem 11.1, find the topological dynamics analog of a corresponding
infinite Ramsey-theoretic result.
\end{problem}

\medskip
We will sketch some possible approaches to Problem \ref{11.2} in (E) and
(F) below. But let us first explain what we mean by 'the
corresponding {\it infinite} Ramsey-theoretic result'. First of
all recall that in the arrow notation the infinite Ramsey theorem
can be stated as,

$$\bbN\rightarrow(\bbN)^{k}_l$$

\noindent for finite numbers $k$ and $l$. This is what one calls
{\it finite-dimensional} Ramsey theorem for $\bbN$. There is also
an {\it infinite-dimensional} Ramsey theorem for $\bbN$ which
states

$$\bbN\rightarrow_{*}(\bbN)^{\bbN}_l,$$

\noindent where $*$ signifies some restriction on the colorings
such as for example the restriction on Borel colorings in the well
known Galvin-Prikry theorem \cite{GaP}. Recall that
in (i)-(vii) we really deal with groups of the form Aut$(\bff)$,
where $\bff$ is a Fra\"{\i}ss\'e limit of a countable Fra\"{\i}ss\'e class $\k$.
The corresponding {\it finite-dimensional} Ramsey theoretic
results deal with arrow-relations of the form

\[
\bff\rightarrow (\bff)^{\bfa}_{l,t} ~~ {\rm and} ~~
\bff\rightarrow (\bff)^{\bfa}_{l}
\]

\noindent for $\bfa\in\k$. In other words, for $\bfa\in\k$, one is
interested in the existence and computation of the {\it big Ramsey
degree} $T(\bfa, \k)$, the minimal integer $t$ such that
$\bff\rightarrow (\bff)^{\bfa}_{l,t}$ for every positive integer
$l$. Of course one is interested also in analogs of the
infinite-dimensional Ramsey theorem such as, for example, the
arrow-relations of the form $\bff\rightarrow_{*} (\bff)^{\bff}_{l,t}$
but at this stage in our knowledge even the theory of
arrow-relations of the form $\bff\rightarrow (\bff)^{\bfa}_{l,t}$
is far from being fully developed. The theory, however, does have
substantial results of this form. We mention one quite old but not
so widely known result due to D. Devlin \cite{De} (see also
Todorcevic \cite{T}) that deals with the class $\l\o$ of all finite
linear orderings. More precisely, Devlin's theorem says that for
every positive integer $k$ there is a positive integer $t$ such
that

$$\bbQ\rightarrow(\bbQ)^{k}_{l,t}$$

\noindent for all positive integers $l$, and that the minimal
integer $t$ satisfying this arrow-relation for all $l$ is equal to
the $(2k+1)$st {\it tangent number} $T_{2k+1}$ given by the
formula $\tan z= \sum^\infty_{n=0}T_nz^n/n!$. Thus if $\bfa_k$
denotes a linearly ordered set of size $k$ then the sequence
$t_k=T(\bfa_k , \l\o)$ is a well studied sequence of numbers
which starts as $t_1=1,t_2=2,t_3=16,t_4=272, \dots$ (see
Knuth-Buckholtz \cite{KB}). We note that the existence of
$T(\bfa_k , \l\o)$ was known to R. Laver (see Erd\"{o}s-Hajnal \cite{EH})
before Devlin's work and that the existence follows also rather
directly from results of K. Milliken \cite{Mill} on which Devlin's work
was based. For results about other Fra\"{\i}ss\'e structures such as for
example the random graph the reader is referred to  Pouzet-Sauer
\cite{PS} and Sauer \cite{S}.

One can ask similar questions for other kinds of extremely amenable
topological groups not directly covered by the list \ref{11.1} (i)--(vii).
A particularly important example is the unitary group
$U(\ell_2)$ of the Hilbert space equipped with the strong operator topology.
The result by Gromov and Milman \cite{GroM} that the unitary
group $U(\ell_2)$ is extremely amenable implies the following
property:
If $f$ is a uniformly continuous function on the unit sphere
${\mathbb{S}}^\infty$ of $\ell_2$ with values in some $\bbR^n$, then for every
$\epsilon >0$ and every compact subset $K$ of ${\mathbb{S}}^\infty$ there is $u\in U
(\ell_2)$ such that the oscillation of $f$ on $u(K)$ is $<\epsilon$.
(See Milman and Schechtman \cite{MS}, Milman \cite{M88} or Gromov \cite{Gro}.)
The exact infinite-dimensional
analog of the above property of spheres is impossible as
demonstrated by Odell and
Schlumprecht \cite{OS} in their solution of the famous distortion problem
for $\ell_2$: there exists a bounded uniformly continuous
$f$ on ${\mathbb{S}}^\infty$ such that $f$ has oscillation 1 in every unit sphere
of an infinite-dimensional subspace of $\ell_2$.
\medskip

\noindent
{\bf (E)}
It seems that the above phenomena can
be described within the following
framework. Recall that the {\it left uniformity}, ${\mathcal{U}}_L(G)$,
of a topological
group $G$ has as basic entourages of the diagonal the sets
\[V_L=\{(x,y)\in G\times G\colon x^{-1}y\in V\},\]
where $V$ is a neighborhood of identity in $G$.
In particular, every topological
group admits the {\it completion with regard to the left uniformity,} also
known as the {\it Weil completion.} (See, e.g., Chapter 10 in \cite{RD}.)
For example, in the case where $G$ is a metrizable group, the left completion
of $G$ is just the metric completion of $(G,d)$, where $d$ is any
left-invariant compatible metric. We will denote the left completion by
$\hat G^L$. While in many cases --- for instance, when
$G$ is locally compact, or abelian, or has small invariant neighborhoods
--- the left completion $\hat G^L$ is again a topological group, in general
it is not the case (Dieudonn\'e \cite{Di}), and the left completion of a
topological group is only a topological
semigroup (with jointly continuous multiplication), see Proposition \ref{10.2}(a) in
\cite{RD}.

For example, the left completion of the unitary group $U(\ell_2)$ with the strong
operator topology can be
identified with the semigroup of all linear isometries from $\ell_2$ to
its subspaces, with the composition of maps as the semigroup operation and
the strong operator topology. The semigroup
$\widehat{{\mathrm{Aut}}(\langle\bbQ,<\rangle)}^L$ is formed by all order-preserving injections
from $\bbQ$ to itself, equipped with the composition operation and the
topology of pointwise convergence on $\bbQ$ viewed as discrete.

If $H$ is a subgroup of $G$, then the left uniform structure on
$G/H$ is, by definition, the finest uniform structure making the
factor-map $\pi\colon G\to G/H$
uniformly continuous with regard to ${\mathcal{U}}_L(G)$.

Here are two examples. Fix a point $\xi$ in the unit sphere ${\mathbb{S}}^\infty$
of the Hilbert space $\ell_2$ (`north pole'), and denote by
$H = {\mathrm{St}}_\xi$, the isotropy subgroup of $\xi$:
\[{\mathrm{St}}_\xi = \{u\in U(\ell_2)\colon u(\xi)=\xi\}.\]
This is a closed subgroup, isomorphic to $U(\ell_2)$ itself.
There is a natural identification
\[U(\ell_2)/{\mathrm{St}}_\xi\ni u{\mathrm{St}}_\xi \mapsto u(\xi)\in{\mathbb{S}}^\infty,\]
as topological $G$-spaces.
The left uniform structure on ${\mathbb{S}}^\infty$ viewed as a factor-space of
the unitary group is the norm uniformity. In other words, basic
entourages of diagonal in ${\mathcal{U}}_L({\mathbb{S}}^\infty)$ are of the form
\[V_\epsilon = \{(\xi,\zeta)\colon \norm{\xi-\zeta}<\epsilon\}.\]

Similarly,
fix any finite set $F\subseteq\bbQ$, and denote by ${\mathrm{St}}_F$ the isotropy
subgroup of $F$, that is, the set of all bijections $\tau\in{\mathrm{Aut}}(\langle\bbQ,<\rangle)$
that leave $F$ (and therefore each element of $F$) fixed. This is an
open subgroup of ${\mathrm{Aut}}(\langle\bbQ,<\rangle)$.
The factor-space ${\mathrm{Aut}}(\bbQ,\leq)/{\mathrm{St}}_F$ can be identified with the
set $[\bbQ]^n$ of all $n$-subsets of $\bbQ$, where $n=\vert F\vert$, under the
correspondence
\[{\mathrm{Aut}}(\langle\bbQ,<\rangle)/{\mathrm{St}}_F\ni
\tau {\mathrm{St}}_F\mapsto \tau(F)\in [\bbQ]^n.\]
The left uniformity on the factor-space
${\mathrm{Aut}}(\langle\bbQ,<\rangle)/{\mathrm{St}}_F\cong [\bbQ]^n$
is discrete.

If $f$ is a real-valued function on a set $X$,
the {\it oscillation} of $f$ is
\[{\mathrm{Osc}}(f) = \sup_{x,y\in X}|f(x)-f(y)|.\]
The following definition is modeled on a classical concept from geometric
functional analysis, first introduced by Milman \cite{M67} in the language of
non-emptiness of the spectrum ${\mathfrak{S}}(f)$ of a function $f$.

Let $f\colon G\to\bbR$ be a left uniformly
continuous function on a topological group $G$.
Say that $f$ is {\it oscillation stable} if
for every $\epsilon>0$ and every right ideal $\mathcal I$ of
$\hat G^L$ there is a right ideal ${\mathcal{J}}\subseteq {\mathcal{I}}$
with the property
\[{\mathrm{Osc}}( f\mid {\mathcal J})<\epsilon.\]
Here we have denoted by the same letter $f$ the (unique) extension of $f$ by
continuity over the left completion $\hat G^L$.

If $H$ is a subgroup of a topological group $G$, we say that a
left uniformly
continuous function $f\colon G/H\to\bbR$ is {\it oscillation stable} if
the composition $\tilde f = f\circ\pi$ with the factor-map
$\pi\colon G\to G/H$ is oscillation stable.

Say that the pair $G,H$, where $H$ is a topological subgroup of a
topological group $G$, is {\it oscillation stable} if every bounded
left uniformly continuous function $f\colon G/H\to\bbR$ is
oscillation stable.
One can show that a pair $G,H$ is oscillation stable if and only if for every bounded
left uniformly continuous function $f\colon G/H\to\bbR$ and every $\epsilon>0$
there is a right ideal $\mathcal I$ of
$\hat G^L$ such that ${\mathrm{Osc}}\,( f\mid {\mathcal I})<\epsilon$.


If $G=U(\ell_2)$ and $H={\mathrm{St}}_\xi$, then a function
$f$ on the unit sphere
${\mathbb{S}}^\infty\cong U(\ell_2)/{\mathrm{St}}_\xi$ is oscillation stable
in the sense of our definition
if and only if it is oscillation stable in the classical sense,
see e.g. Definition 13.1 in \cite{BL}. The result by
Odell and Schlumprecht \cite{OS} that the infinite-dimensional Hilbert space has
the distortion property is equivalent, in our language, to
saying that the pair $U(\ell_2),{\mathrm{St}}_\xi$ is not oscillation stable
for some (any) $\xi\in{\mathbb{S}}^\infty$.

Also, it follows from Devlin's theorem that the pair
${\mathrm{Aut}}(\langle\bbQ,<\rangle)$, ${\mathrm{St}}_F$ is oscillation stable if
and only if $|F|=1$.


One question that remains unanswered, is: {\it does there exist a
non-trivial oscillation stable topological group,} that is, a topological
group $G\neq\{e\}$ for which the pair $G,\{e\}$ is 
oscillation stable?\footnote{Recently Hjorth informed the authors that he can prove that no non-trivial Polish group is oscillation stable.}

Oscillation stability for topological groups
is a strictly
stronger concept than extreme amenability, as, for instance, both
topological groups
$U(\ell_2)$ and ${\mathrm{Aut}}(\langle\bbQ,<\rangle)$ are extremely amenable
but not oscillation stable.
\medskip

\noindent {\bf (F)} For topological groups with small open
neighborhoods one is able to capture the quantitative, as well as
qualitative, content of results from Ramsey theory similar to
Devlin's theorem considered above in ({\bf C}).

Let $G$ be a topological group, $H$ be a subgroup. Define the {\it
big oscillation degree}  $N(G,H)$ of
$G$ to be the smallest $t$ such that:

(**) for any bounded left uniformly continuous function
$f$ from $G/H$ to some Euclidean space and
any $\epsilon >0$, there is a right ideal $\mathcal I$ in the left completion
$\hat G^L$ of $G$ such that
(denoting also by $f$ the extension of the lift of $f$ to $\hat G^L$ )
$N(\epsilon , f | {\mathcal I}) \leq t$.

If no such $t$ exists, we say again that $N(G)$ is infinite. We have
$n(G,H) \leq N(G,H) \leq N(G)$ ($= N(G,\{1\})$.

If $H$ is a subgroup of $G$, then the condition
$N(G,H) = 1$ is equivalent to the oscillation stability of $G, H$,
and it turns out that a group $G$ with small open subgroups is oscillation stable iff
$N(G,H)=1$, for all open subgroups $H$ of $G$.

Let $G$ be a non-trivial group with small open subgroups. Any such $G$ is the
automorphism group of a Fra\"{\i}ss\'e structure $\bff$. It turns out now
that if the signature of $\bff$ is relational and finite (or even more generally
if there are only finitely many, up to isomorphism, 2-generated
structures in $\k$), then $G$ is not oscillation stable. On the other hand, we have that in
general for any $\k$,$\bff$,and $G$ as above and $\bfa$ in $\k$,
\[N(G, H(\bfa)) = T(\bfa ,\k).\]

\noindent Thus, if $\k$ is the class of finite linear orderings
and if $\bfa_n$ is a finite linear ordering of size $n$, then
applying Devlin's theorem we get that $N(G, H(\bfa_n))$ is equal
to the $n$th odd tangent number $t_n$.

\medskip
{\bf Addendum}.  We have recently received the preprint Ne\v set\v ril \cite{Nes04a}, which discusses many concepts and results of structural Ramsey theory
relevant to our paper.  In particular, some new examples of classes with
the Ramsey property are presented concerning posets, directed graphs, etc.
One very interesting case is the class of all structures $\langle
A,\sqsubset ,\prec\rangle$, where $\langle A,\sqsubset\rangle$ is a
finite poset and $\prec$ is a linear extension of $\sqsubset$ (see Ne\v set\v ril--R\"odl \cite{NR84}, Fouch\'{e} \cite {Fo97}). The
Fra\"{\i}ss\'{e} limit of this class is of the form $\bff =\langle F_0,
\sqsubset_0,\prec_0\rangle$, where $\langle F_0,\sqsubset_0\rangle$ is
the random poset (i.e., the Fra\"{\i}ss\'{e} limit of the class of finite
posets) and $\prec_0$ is an appropriate linear extension of $\sqsubset_0$.
In particular, it follows from our results here that Aut$(\bff )$ is
extremely amenable and that the universal minimal flow of the
automorphism group of the random poset $\bff_0=\langle F_0,\sqsubset_0
\rangle$ is the space of all linear extensions of $\sqsubset_0$.

The proof of the result announced in Ne\v set\v ril \cite{N03} is contained in the recent preprint Ne\v set\v ril \cite{N04}. Finally, Nguyen Van The \cite{Ng} has shown that the class of finite convexly ordered ultrametric spaces (where an ordered metric space is convexly ordered if each metric ball is an interval) has the Ramsey and ordering properties and uses this to compute the universal minimal flow of the isometry group of the Baire space.

\section*{Appendix 1. A new proof of Veech's theorem}
Veech's theorem (Theorem 2.2.1 in \cite{V}) 
is an important result of abstract topological dynamics, asserting
that every locally compact group acts freely on a suitable compact space.
Alternative proofs of this result can be found in \cite{AS}
(in the second countable case) and in
\cite{Py}. The latter author notes that his proof is `really the same,'
but it emphasizes different features of the original idea. 
The same applies to our proposed proof, which is, we hope, more accessible.
\medskip

\noindent{\bf Lemma A1.1.} 
{\it Let $G$ be a locally compact group, and let $g\in G$,
$g\neq e$. There exists a 
right invariant continuous pseudometric $d$ on $G$, bounded by $1$ and 
such that $0<d(e,g)<1$ and the closure of the open ball of unit radius 
is compact.}

\medskip
\noindent{\bf Proof.}
Let $\nu$ be a left-invariant Haar measure on $G$. For a $f\in L^2(G,\nu)$ 
and $h\in H$, define $^hf\in L^2(G,\nu)$ via
$^hf(x):= f(h^{-1}x)$. 

{\it Case 1: $g^2\neq e$.}
Choose a symmetric compact neighborhood of the identity, $V$, in $G$, with
the property $g,g^2\notin V^2$, and a 
function $f\in L^2(G,\nu)$ supported on $V$ and such that the $L^2$-norm
$\norm f = 1$. Let $\phi = f+\,\,^{g^{-1}}f$. Clearly, $\norm\phi =\sqrt 2$. 
For each $x,y\in G$, define
\[\rho(x,y):= \norm{^{x^{-1}}\phi - \,\,^{y^{-1}}\phi}
\equiv \norm{\phi-\,\,^{yx^{-1}}\phi}.\]
This $\rho$ is a right-invariant continuous pseudometric on $G$, bounded by
$2\sqrt 2$. If translates of $\phi$ by $x^{-1}$ and by $y^{-1}$
are orthogonal, then $\rho(x,y)=2$. 
It follows that,
if $h\in G$ and $\rho(e,h)<2$, then 
$(V \cup g^{-1}V) \cap (h^{-1}V \cup h^{-1}gV) \neq\emptyset$.
This is equivalent to
$h \in V^2 \cup gV^2 \cup V^2 g \cup gV^2 g$, and so the open ball 
${\mathcal O}_2(e)$ of radius $2$ has compact closure.
Also,  $\rho(e,g) = \norm{f-\,\,^{g^{-2}}f} = \sqrt 2$,
because the supports of $f$ and $^{g^{-2}}f$ are disjoint.
The required pseudometric $d$ is now defined by 
\[d(x,y) = \frac 1 2\min\{\rho(x,y), 2\}.\]

\medskip

{\it Case 2: $g^2=e$.} Let $V$ be a compact symmetric 
neighborhood of identity with $g\notin V^2$. Let $f$ be
of $L^2$-norm one and supported on $V$, and let
$\phi = f + 2\cdot\,^gf$.
Define the right-invariant continuous pseudometric $\rho$ via
$\rho(x,y):= \norm{^{x^{-1}}\phi - \,\,^{y^{-1}}\phi}
\equiv \norm{\phi-\,\,^{yx^{-1}}\phi}$. 
Similarly to Case 1, 
the closure of the open ball of radius $\sqrt{10}$ is compact.
Finally, set
\[d(x,y) = \frac 1 {\sqrt{10}}\min\{\rho(x,y),\sqrt{10}\}.\]
\qed
\medskip

\noindent{\bf Veech's Theorem.} 
{\it Every locally compact group $G$ acts freely on the greatest ambit 
$S(G)$.}
\medskip

\noindent{\bf Proof.} Let $G$ be a locally compact group. 
Let $g\in G$ and $g\neq e$. We will show that $g$
has no fixed points in the greatest ambit of $G$.
\smallskip

(A) Choose a a pseudometric $d$ on $G$ as in Lemma A1.1.

\smallskip
(B) Choose an $\varepsilon>0$ satisfying 
\begin{equation}
\label{epsilon} 9{\varepsilon} < d(e,g)< 1 - 4{\varepsilon}. \end{equation}
By Zorn's lemma, there exists a maximal subset
$A\subseteq G$ with the property that whenever $a,b\in A$ and $a\neq
b$, one has ${\mathcal{O}}_{\varepsilon}(a)\cap 
{\mathcal{O}}_{\varepsilon}(b)\neq\emptyset$. Such an $A$ is 
a $2{\varepsilon}$-net in $G$.

\smallskip
(C) Define a graph, $\Gamma$, whose vertex set is $A$ and such
that two vertices, $a,b\in A$, are adjacent if and only
if $a\neq b$ and $ab^{-1}\in K^2$. 

Let $\varkappa$ denote the cardinality of an arbitrary finite
family, $\gamma$, of open balls of radius ${\varepsilon}$ covering the
compact set $K^2$. Any family $\delta$ of pairwise disjoint open balls of
radius ${\varepsilon}$ with centers in $K^2$ has cardinality not exceeding
$\varkappa$, because every mapping assigning to every $B\in\delta$ a ball
$B'$ in the family $\gamma$ containing the center of $B$
is an injection. If $a,b_1,\ldots,b_n\in A$ and $a$ is adjacent to
each $b_i$, then $b_ia^{-1}\in K^2$ for all $i$ and the ${\varepsilon}$-balls
centered at $b_ia^{-1}$, $i=1,2,\ldots,n$, are pairwise disjoint
(the metric $d$ is right-invariant). It follows that
$n\leq\varkappa$ and $\Gamma$ has a finite degree $\leq\varkappa$.
\smallskip

(D) As a consequence, the vertices of $\Gamma$ can be colored
with at most $\varkappa+1$ colors in such a way that no two
adjacent vertices have the same color. Let $A =
\bigsqcup_{i=1}^{m}A_i$, $A_i \not= \emptyset$,  be such a coloring, where
$m\leq\varkappa+1$.
\smallskip

(E) For each $i=1,2,\ldots,m$ define a function $f_i\colon G\to
{\mathbb{R}}$ via
\[f_i(x):=d(x,A_i),\]
the distance from $x\in G$ to $A_i$. This $f_i$ is a 1-Lipschitz
function, bounded by $1$. 
\smallskip

(F) Let $a,b\in A$ be such that $d(a,b)<1$. For some
$i=1,2,\ldots,m$, one has $a\in A_i$ and, since $a$ and $b$ are
adjacent, also $b\notin A_i$. Moreover, $a$ is the
only element of $A_i$ at a distance $<1$ from $b$ and thus the
nearest neighbor to $b$ in $A_i$. Indeed, assuming that there is
a $c\in A_i$ with $c\neq a$ and $d(c,b)<1$, one has $ba^{-1}\in
K$, $cb^{-1}\in K$, and thus $ca^{-1} = cb^{-1}ba^{-1}\in K^2$,
meaning that $c$ and $a$ are adjacent, in contradiction with the
choice of the coloring.

Since each function $f_j$, $j=1,2,\ldots,m$ is
$1$-Lipschitz, 
\begin{equation}\label{ellinf}
\max_{j=1}^m |f_j(a) - f_j(b)| = d(a,b)\mbox{ whenever }a,b\in A\mbox{ and
}d(a,b)<1.\end{equation}

(G) Define the mapping $f\colon G\to \ell^\infty(m)$, where $\ell^\infty (m) = 
\mathbb{R}^m$ with the max norm $\lVert . \rVert_\infty$, as
$f=(f_1,f_2,\ldots,f_m)$. Then (\ref{ellinf})
is equivalent to
\begin{equation} \forall a,b\in
A,~~\left(d(a,b)<1\right)\Rightarrow \lVert f(a)-f(b)\rVert_\infty
= d(a,b). \label{norm}
\end{equation}

(H) Let 
$x\in G$ be arbitrary. There are $a,b\in A$ such that $d(x,a)<2{\varepsilon}$
and $d(gx,b)<2{\varepsilon}$. Since $d(gx,x) = d(g,e)$, it follows by the
triangle inequality that $5{\varepsilon} <d(a,b)< 1$, and
Eq. (\ref{norm}) implies that
\[\lVert f(a)-f(b)\rVert_\infty  = d(a,b) >5{\varepsilon}.\]
By the triangle inequality \begin{eqnarray*}\lVert
f(x)-f(gx)\rVert_\infty &\geq& \lVert f(a)-f(b)\rVert_\infty -
4{\varepsilon}\\ & >& {\varepsilon}.
\end{eqnarray*}

(I) The mapping $f\colon G\to {\mathbb{R}}^m$, being right uniformly
continuous and bounded, admits a unique continuous extension,
$\bar f$, over the greatest ambit $S(G)$ of $G$. By
continuity, 
\[\forall x\in S(G),~~\lVert
f(x)-f(g\cdot x)\rVert_\infty \geq {\varepsilon} > 0.\] In particular,
the action by $g$ on the greatest
ambit is fixed point-free.
\qed

\section*{Appendix 2. Non-metrizability of the universal minimal
flow for non-compact locally compact groups}

A subset $A$ of a group $G$ is called ({\it discretely}) {\it left syndetic}
if finitely many left translates of $A$ cover $G$.
Here is a simple and well-known fact from abstract topological dynamics.
\medskip

\noindent{\bf Lemma A2.1.} {\it Let $G$ be a topological group, and let $M$
be a minimal compact $G$-flow. Let $W\subseteq M$ be a non-empty open
subset, and let $\xi\in M$. Then the set
\[\widetilde W:= \{g\in G\colon g\cdot \xi\in W\}\]
is discretely left syndetic in $G$.
}

\medskip
\noindent{\bf Proof}. 
The translates $h\cdot W$,
$h\in G$ form an open cover of $M$, because otherwise there would be a
point $\zeta\in M$ whose $G$-orbit misses $W$, in contradiction with the
assumed minimality of $M$. Choose finitely many
elements, $h_1,h_2,\ldots, h_n\in G$ with the property that $h_i\cdot W$,
$i=1,2,\ldots,n$, cover $M$. It remains to notice that for every 
$h\in G$, $\widetilde{h\cdot W} = h\widetilde W$, and so the left translates
$h_i\widetilde W$, $i=1,2,\ldots,n$, cover $G$. 
\qed
\medskip

\noindent{\bf Theorem A2.2.} {\it The universal minimal flow $M(G)$ of 
a non-compact
locally compact group $G$ is non-metrizable.}

\medskip
\noindent{\bf Proof}.
Let $G$ be a locally compact group. Let $U$ be a neighborhood of identity
whose closure is compact. 
Use Zorn's lemma to choose a maximal subset $X\subseteq G$
with the property that $\{Ux\colon x\in X\}$ is a disjoint family. 

According to Pym \cite{Py} (the Local Structure Theorem on
p. 172), the closure $\overline X$ of $X$ in 
the greatest ambit $S(G)$ is homeomorphic to $\beta X$, the
Stone-\v Cech compactification of the discrete space $X$. Also, if $V$ is
an open subset of $G$ with $\overline V\subseteq U$, then
the subspace $V\cdot\overline X$ is open in $S(G)$ and
homeomorphic with $V\times \overline X$ (and consequently with
$V\times \beta X$) under the map $(v,\xi) \mapsto
v\cdot \xi$, $v\in V$, $\xi\in \overline X$. 
Finally, given any $\xi\in S(G)$, $U$ and $X$ can be chosen
so that $\xi\in \overline X$. 

Denote by $M$ an isomorphic
copy of the universal minimal flow $M(G)$ sitting inside the greatest
ambit $S(G)$.
Assume from now on that $U$ and $X$ as above are chosen in such a way 
that $\overline X\cap M
\neq\emptyset$. 
Let also $V$ and $V_1$ be open neighborhoods of 
identity in $G$ with the property $V\subseteq \overline V \subseteq V_1
\subseteq\overline{V_1}
\subseteq U$. 
Since $V\cdot\overline X$ is open in $S(G)$, it follows that
$(V\cdot\overline X) \cap M$ is (non-empty and) open in $M$. 

Assume now that $M$ is metrizable, in order to deduce that $G$ is compact.

The closed subspace $(\overline V\cdot\overline X) \cap M$ of $M$ is also
metrizable and compact. The second coordinate projection, ${\mathrm{proj}}_2$,
from $V_1\cdot\overline X\cong V_1\times \beta X$
to $\overline X\cong \beta X$ is continuous, and therefore the image
$K={\mathrm{proj}}_2((\overline V\cdot\overline X) \cap M)$ is a compact
metrizable subspace of the extremally disconnected space $\beta X$. 

Since an extremally disconnected space does not contain any nontrivial
convergent sequences (see e.g. \cite{Eng}, Exercise 6.2.G.(a) on p. 456),
it follows that $K$ is finite. Consequently, for each $\kappa\in K$ the subset
$(V\cdot\kappa)\cap M$ is open in $(V\cdot\overline X) \cap M$, and
we conclude that for some $\kappa^\prime\in \overline X$ the set
$W = (V\cdot \kappa^\prime)\cap M$ is non-empty and 
open in $(V\cdot\overline X) \cap M$
and therefore in $M$ itself. 

Let $\xi\in W = (V\cdot \kappa^\prime)\cap M$ be arbitrary.
By Lemma A2.1, the set
\[\widetilde W = \{g\in G\colon g\cdot \xi\in W\}\]
is discretely 
left syndetic in $G$. For some $v\in V$, one has $\xi = v\cdot\kappa^\prime$.
The set 
\[W^\ddag:= \{g\in G\colon g v\cdot\kappa^\prime \in 
V\cdot\kappa^\prime\}\]
 is bigger than $\widetilde W$ and therefore also
discretely left syndetic. Since the action of $G$ on the greatest ambit
$S(G)$ is free by Veech's theorem, the condition
$g v\cdot\kappa^\prime \in  V\cdot\kappa^\prime$ is equivalent to
$gv\in V$ and, in its turn, implies $g\in V^2$. It follows that
$W^\ddag\subseteq V^2$. 

The compact set $(\overline V)^2$ contains $V^2$ and is therefore
discretely left syndetic as well. Consequently, $G$ is compact. 
\qed


\begin{thebibliography}{100}

\bibitem{AH78}
F.G. Abramson and L.A. Harrington, Models without indiscernibles,
{\it J. Symbolic Logic} {\bf 43} (1978), 572--600.

\bibitem{AS} 
S. Adams and G. Stuck, Splitting of nonnegatively curved leaves in
minimal sets of foliations, {\it Duke Math. J.} (1993), {\bf 71 (1)},
71--92.

\bibitem{A} 
J. Auslander, {\it Minimal Flows and Their Extensions}, North Holland,
1988.

\bibitem{BK} 
H. Becker and A.S. Kechris, {\it The Descriptive Set Theory of
Polish Group Actions}, London Math. Soc. Lecture Note Series, {\bf 232},
Cambridge Univ. Press, 1996.

\bibitem{BL} Y. Benyamini and J. Lindenstrauss [00],
{\it Geometric Nonlinear Functional Analysis,} Vol. 1,
Colloquium Publications {\bf 48}, AMS, 
Providence, RI, xi+488 p., 2000.

\bibitem{Be} 
S.K. Berberian, {\it Lectures in Functional Analysis and Operator
Theory}, Springer-Verlag, 1974.

\bibitem{Bo}
S.A. Bogatyi [02], Metrically homogeneous spaces, {\it Russian Math. Surveys} 
{\bf 57(2)} (2002), 221--240.

\bibitem{Ca} 
P.J. Cameron, {\it Oligomorphic Permutation Groups}, London Math.
Society Lecture Note Series, {\bf 152}, 1990.

\bibitem{Ch} 
G.L. Cherlin [98], {\it The Classification of Countable Homogeneous
Directed Graphs and Countable Homogeneous $n$-tournaments}, Memoirs
of the Amer. Math. Soc., {\bf 131}, No. 621, 1998.

\bibitem{Co} 
C. Constantinescu, {\it $C^*$-algebras, Vol. 2:  Banach Algebras
and Compact Operators}, North Holland, 2001.

\bibitem{De}
D. Devlin, {\it Some partition theorems and ultrafilters on $\omega$}, 
Ph.D. Thesis, Dartmouth College, 1979.

\bibitem{dV} 
J. de Vries, {\it Elements of Topological Dynamics}, Kluwer, 1993.

\bibitem{Di} J. Dieudonn\'e, Sur la completion des groupes topologiques,
{\it C. R. Acad. Sci. Paris} {\bf 218} (1944), 774--776.

\bibitem{E60} 
R. Ellis, Universal minimal sets, {\it Proc. Amer. Math. Soc.} {\bf
11} (1960), 540--543.

\bibitem{E69} 
R. Ellis, {\it Lectures on Topological Dynamics}, W.A. Benjamin, 1969.

\bibitem{Eng}
R. Engelking, {\it General Topology}, Math. Monographs, {\bf 60}, 
PWN - Polish Scient. Publishers, Warsaw, 1977.

\bibitem{EH}
P. Erd\"{o}s and A. Hajnal, 
Unsolved and solved problems in set theory, 
{\it Proc. Symp. Pure Math.} {\bf 25} (1974), 269--287.

\bibitem{EHR65}
P. Erd\"{o}s, A. Hajnal and R. Rado, Partition relations for cardinal numbers, {\it Acta Math. Acad. Sci. Hung.}, {\bf 16}, (1965), 93--196.

\bibitem{ER56}
P. Erd\"{o}s and R. Rado, A partition calculus in set theory, {\it Bull. Amer. Math. Soc.}, {\bf 62} (1956), 427--489.

\bibitem{Fo97}
W.L. Fouch\'{e}, Symmetry and the Ramsey degree of posets, {\it Discrete Math.}
{\bf 167/168} (1997), 309-315.

\bibitem{Fo98}
W.L. Fouch\'{e}, Symmetries in Ramsey theory, {\it East-West J.
Math.}, {\bf 1} (1998), 43--60.

\bibitem{Fo99}
W.L. Fouch\'{e}, Symmetry and the Ramsey degrees of finite relational
structures, {\it J. Comb. Theory Ser. A} {\bf 85(2)} (1999), 135--147.

\bibitem{Fo99a}
W.L. Fouch\'{e}, Symmetries and Ramsey properties of trees, {\it
Discrete Math.} {\bf 197/198} (1999), 325--330.

\bibitem{Fr} 
R. Fra\"{\i}ss\'{e}, Sur l'extension aux relations de quelques propriet\'es
des ordres, {\it Ann. Sci. \'Ecole Norm. Sup.} {\bf 71} (1954), 363-388.

\bibitem{GaP} F. Galvin and K. Prikry,
Borel sets and Ramsey's theorem, {\it J. Symb. Logic} {\bf 38} (1973), 
193-198. 

\bibitem{GiP} 
T. Giordano and V. Pestov, Some extremely amenable groups, {\it C.R. Math. 
Acad. Sci. Paris} {\bf 334 (4)} (2002), 273--278.

\bibitem{GiP2} 
T. Giordano and V. Pestov, Some extremely amenable groups related to
operator algebras and ergodic theory, ArXiv e-print math.OA/0405288,
May 2004, 37 pp.

\bibitem{Gl98} 
E. Glasner, On minimal actions of Polish groups, {\it Top. Appl.}
{\bf 85} (1998), 119--125.

\bibitem{Gl00} 
E. Glasner, Structure theory as a tool in topological dynamics, {\it Descriptive
Set Theory and Dynamical Systems}, M. Foreman et al. (Eds.), London Math.
Society Lecture Note Series, {\bf 277}, 173--209, Cambridge Univ. Press, 2000.

\bibitem{GlW02} 
E. Glasner and B. Weiss, Minimal actions of the group $S(\bbZ )$ of
permutations of the integers, {\it Geom. and Funct. Anal.} {\bf 12} (2002),
964--988.

\bibitem{GlW03} 
E. Glasner and B. Weiss, The universal minimal system for the group
of homeomorphisms of the Cantor set, {\it Fundam. Math.}
{\bf 176} (2003), 277-289. 

\bibitem{GLR} 
R.L. Graham, K. Leeb, and B.L. Rothschild, Ramsey's theorem for a
class of categories, {\it Adv. in Math} {\bf 8} (1972), 417--433.

\bibitem{GR} 
R.L. Graham and B.L. Rothschild, Ramsey's theorem for $n$-parameter
sets, {\it Trans. Amer. Math. Soc.} {\bf 159} (1971), 257--292.

\bibitem{GRS}
R.L. Graham, B.L. Rothschild, and J.H. Spencer, {\it Ramsey Theory}, 
2nd Edition, Wiley, 1990.

\bibitem{Gra} E. Granirer, 
Extremely amenable semigroups I,II,
{\it Math. Scand.} {\bf 17} (1965), 177-179; {\bf 20} (1967), 93-113. 

\bibitem{GraL}
E. Granirer and A.T. Lau, Invariant means on locally compact groups,
{\it Ill. J. Math}. {\bf 15} (1971), 249--257.

\bibitem{Gro}
M. Gromov, Filling Riemannian manifolds, {\it J. Diff. Geometry} {\bf 18} (1983), 
1--147.

\bibitem{GroM} 
M. Gromov and V.D. Milman, A topological application of the isoperimetric
inequality, {\it Amer. J. Math.} {\bf 105} (1983), 843--854.

\bibitem{HC}
W. Herer and J.P.R. Christensen, On the existence of pathological 
submeasures and the construction of exotic topological groups, 
{\it Math. Ann.} {\bf 213} (1975), 203--210.

\bibitem{H} 
W. Hodges, {\it Model Theory}, Cambridge Univ. Press, 1993.

\bibitem{KB} 
D.E. Knuth and T.J. Buckholtz, Computation of Tangent, Euler, and Bernoulli 
numbers, {\it  Math. of Computation} {\bf 21 (100)} (1967), 663--688.

\bibitem{LW} 
A.H. Lachlan and R. Woodrow, Countable ultrahomogeneous undirected
graphs, {\it Trans. Amer. Math. Soc.} {\bf 262} (1980), 51--94.

\bibitem{LMP} 
A.T.-M. Lau, P. Milnes, and J. Pym, On the structure of minimal
left ideals in the largest compactification of a locally compact group,
{\it J. London Math. Soc} {\bf 59 (2)} (1999), 133--152.

\bibitem{L}
M. Ledoux, {\it The concentration of measure phenomenon}. Math. Surveys and 
Monographs, {\bf 89}, Amer. Math. Soc., 2001.

\bibitem{Mill}
K. Milliken, A Ramsey theorem for trees, {\it J. Comb. Theory} {\bf 26} (1979), 
215--237.

\bibitem{M67}
 V.D. Milman, Infinite-dimensional geometry of the unit sphere in Banach space,
{\it Sov. Math., Dokl.} {\bf 8} (1967), 1440-1444.

\bibitem{M88}
V.D. Milman, The heritage of P. L\'{e}vy in geometrical functional analysis, 
{\it Ast\'{e}risque} {\bf 157--158} (1988), 273--301.

\bibitem{MS}
V.D. Milman and G. Schechtman, {\it Asymptotic theory of finite-dimensional 
normed spaces (with an Appendix by M. Gromov)}, 
Lecture Notes in Math., {\bf 1200}, Springer, 1986.

\bibitem{Mit} T. Mitchell, 
Fixed points and multiplicative left invariant means,
{\it Trans. Am. Math. Soc.} {\bf 122} (1966), 195-202. 

\bibitem{N89} 
J. Ne\v set\v ril, For graphs there are only four types of
hereditary Ramsey classes, {\it J. Comb. Theory, Series B} {\bf 46 (2)}
(1989), 127--132.

\bibitem{N95} 
J. Ne\v set\v ril, Ramsey Theory, {\it Handbook of Combinatorics},
R. Graham et al. (Eds.), 1331--1403, Elsevier, 1995.

\bibitem{N03}
J. Ne\v set\v ril, Private communication, 2003.

\bibitem{Nes04a}
J. Ne\v set\v ril, Ramsey classes and homogeneous structures, Preprint,
2003.

\bibitem{N04}
J. Ne\v set\v ril, Metric spaces are Ramsey, Preprint, 2004.

\bibitem{NPRV} 
J. Ne\v set\v ril, H.J. Pr\"omel, V. R\"odl and B. Voigt, Canonical
ordering theorems, a first attempt, Proc. of the 10th Winter School on
Abstract Analysis (Srn\'\i , 1982), {\it Rend. Circ. Math. Palermo (2)} (1982),
193--197.

\bibitem{NR77} 
J. Ne\v set\v ril and V. R\"odl, Partitions of finite relational and 
set systems, {\it J. Comb. Theory} {\bf 22 (3)} (1977), 289--312.

\bibitem{NR78} 
J. Ne\v set\v ril and V. R\"odl, On a probabilistic graph-theoretical
method, {\it Proc. Amer. Math. Soc.} {\bf 72 (2)} (1978), 417--421.

\bibitem{NR83} 
J. Ne\v set\v ril and V. R\"odl, Ramsey classes of set systems,
{\it J. Comb. Theory, Series A} {\bf 34 (2)} (1983), 183--201.

\bibitem{NR84} 
J. Ne\v set\v ril and V. R\"odl, Combinatorial partitions of finite posets and lattices -- Ramsey lattices, {\it Alg. Univers.} {\bf 19)} (1984), 106--119.

\bibitem{NR90}
J. Ne\v set\v ril and V. R\"odl, {\it Mathematics of Ramsey Theory}, 
Springer, 1990.

\bibitem{Ng}
L. Nguyen Van The, Ramsey degrees of finite ultrametric spaces, preprint, 2004

\bibitem{OS}
E. Odell and T. Schlumprecht, The distortion problem, 
{\it Acta Mathematica} {\bf 173} (1994), 259--281.

\bibitem{P98} 
V. Pestov, Some universal constructions in abstract topological dynamics,
{\it Contemporary Math.} {\bf 215} (1998), 83--99.

\bibitem{P98a} 
V. Pestov, On free actions, minimal flows and a problem by Ellis, {\it
Trans. Amer. Math. Soc.} {\bf 350 (10)} (1998), 4149--4165.

\bibitem{P99} 
V. Pestov, Topological groups:  where to from here?, {\it Topology
Proceed.} {\bf 24} (1999), 421--502.

\bibitem{P02} 
V. Pestov, Ramsey-Milman phenomenon, Urysohn metric spaces, and extremely
amenable groups, {\it Israel J. Math.} {\bf 127} (2002), 317--357.

\bibitem{P02a} 
V. Pestov, Remarks on actions on compacta by some infinite-dimensional
groups, {\it Infinite dimensional Lie groups in geometry and representation 
theory, (Washington, DC 2000)}, World Sci., 145--163, 2002.

\bibitem{P02b} 
V. Pestov, $MM$-spaces and group actions, {\it L' Enseign. Math.}
{\bf 48} (2002), 209--236.

\bibitem{PS}
M. Pouzet and N. Sauer, Edge partitions of the Rado graph, 
{\it Combinatorica} {\bf 16(4)} (1996), 505--520.

\bibitem{Pr}
H.J. Pr\"omel, Some remarks on natural orders for combinatorial
cubes, {\it Discrete Math}. {\bf 73} (1989), 189--198.

\bibitem{Py} 
J. Pym, A note on $G^{\l\calu\calc}$ and Veech's Theorem, {\it Semigroup
Forum} {\bf 59} (1999), 171--174.

\bibitem{Ra}
R. Rado, Direct decomposition of partitions, {\it J. London Math. Soc.}
{\bf 29} (1954), 71--83.

\bibitem{RD} W. Roelcke and S. Dierolf, {\it Uniform Structures on Topological
Groups and Their Quotients,} McGraw-Hill, 1981. 

\bibitem{Ru} 
W. Rudin, {\it Functional Analysis}, McGraw-Hill, 1973.

\bibitem{S}
N. Sauer, Coloring finite substructures of countable structures,
Paul Erd\"{o}s and his Mathematics, II, {\it Bolyai Math. Studies,} {\bf 11}
(2002), 525--553.

\bibitem{Th} 
S. Thomas, Groups acting on infinite dimensional projective spaces,
{\it J. London Math. Soc.} {\bf 34 (2)} (1986), 265--273.

\bibitem{T}
S. Todorcevic, {\it Ramsey spaces}, to appear, 2004.

\bibitem{Tu} 
S. Turek, Universal minimal dynamical system for reals, {\it Comment. 
Math. Univ. Carolin.} {\bf 36 (2)} (1995), 371--375.

\bibitem{U} 
P. Urysohn, Sur un espace m\'etrique universel, {\it Bull. Sci.
Math.} {\bf 51} (1927), 43--64, 74--90.

\bibitem{Usp90} 
V. Uspenskij, On the group of isometries of the Urysohn universal 
metric space, {\it Comment. Univ. Carolinae} {\bf 31(1)} (1990), 181-182.

\bibitem{Usp00} 
V. Uspenskij, On universal minimal compact $G$-spaces, {\it Topology
Proc.} {\bf 25} (2000), 301--308.

\bibitem{Usp02} 
V. Uspenskij, Compactification of topological groups, {\it Proc. of the 9th
Prague Topol. Symp. (Prague 2001)}, 331--346 (electronic), Topol. Atlas,
Toronto, 2002.

\bibitem{V} 
W. Veech, Topological dynamics, {\it Bull. Amer. Math. Soc.}
{\bf 83 (5)} (1977), 775--830.

\end{thebibliography}
\end{document}